\documentclass[12,reqno]{amsart}

\usepackage{amsfonts}
\usepackage[active]{srcltx}
\usepackage{amsmath}
\usepackage{amssymb}
\usepackage{amsthm}
\usepackage{mathtools}
\usepackage{bm} 
\usepackage{empheq}
\usepackage{mathptmx}
\usepackage{ulem}

\usepackage{enumitem}
\usepackage{wasysym}
\usepackage{verbatim}
\usepackage{graphicx}
\usepackage[
bookmarks=true,         
bookmarksnumbered=true, 
colorlinks=true, pdfstartview=FitV, linkcolor=blue, citecolor=blue,
urlcolor=blue]{hyperref}
\usepackage{microtype}
\theoremstyle{plain}\newtheorem{theorem}{Theorem}[section]

\newtheorem{lemma}[theorem]{Lemma}
\newtheorem{problem}[theorem]{Problem}

\renewenvironment{proof}[1][Proof]{\textbf{#1.} }{\ \rule{0.5em}{0.5em} \par }
\theoremstyle{remark}

\theoremstyle{definition}
\newtheorem{remark}[theorem]{Remark}

\def\PP{\mathbb{P}}
\def\SS{\mathbb{S}}
\def\RR{\mathbb{R}}
\def\EE{\mathbb{E}}

\def\FF{\mathbb{F}}

\def\BB{\mathbb{B}}
\def\CC{\mathbb{C}}

\def\cB{{\mathcal B}}

\def\cC{{\mathcal C}}

\def\cU{{\mathcal U}}

\def\la{{\lambda}}

\def\cF{{\mathcal F}}

\def\Om{{\Omega}}

\def\al{{\alpha}}

\def\Ga{{\Gamma}}
\def\ga{{\gamma}}

\def\tr{{ \hbox{ Tr} }}

\def\la{{\lambda}}
\def\La{{\Lambda}}
\def\vare{{\varepsilon}}

\def\al{{\alpha}}

\newcommand{\laa}{\Lambda}

\newcommand{\vp}{\varphi}

\renewcommand{\SS}{\mathcal S}

\let\Section=\section
\def\section{\setcounter{equation}{0}\Section}

\begin{document}

\title[]{A linear-quadratic partially observed Stackelberg stochastic differential game with multiple followers and its application to multi-agent formation control}
  \author{\sc Yichun Li }
\address{Key Laboratory of Smart Manufacturing in Energy Chemical Process, Ministry of Education, East China University of Science and Technology, Shanghai 200237, China}
\email{yichunli953@gmail.com}
  \author{\sc Yaozhong Hu }
\address{Department of Mathematical and Statistical Sciences, University of Alberta at Edmonton, Edmonton, T6G 2G1, Canada}
\email{yaozhong@ualberta.ca}
  \author{\sc Jingtao Shi}
\address{School of Mathematics, Shandong University, Jinan, Shandong 250100, China}\email{shijingtao@sdu.edu.cn}
  \author{\sc Yueyang Zheng}
\address{School of Mathematical Sciences, Fudan University, Shanghai 200433, China}\email{zhengyueyang0106@163.com}

\thanks{Y. Li is supported by the National Natural Science Foundation of Shanghai (25ZR1402103) and China Postdoctoral Science Foundation (2024M750904).}
\thanks{Y. Hu is supported by an NSERC discovery fund and a centennial fund of University of Alberta.}
\thanks{J. Shi is supported by National Key R\&D Program of China (2022YFA1006104), National Natural Science Foundations of China (12471419, 12271304), and Shandong Provincial Natural Science Foundation (ZR2024ZD35).}
\thanks{Y. Zheng is supported by National Natural Science Foundations of China (12401577), the Postdoctoral Fellowship Program (Grade B) of China Postdoctoral Science Foundation (GZB20230168), and the China Postdoctoral Science Foundation (2025T180851, 2024M750488).}

\date{}
\maketitle

\begin{abstract} In this paper, we study a linear-quadratic partially observed Stackelberg stochastic differential game problem in which a single leader and multiple followers are involved. We consider more practical formulation for partial information that none of them can observed the complete information and the followers know more than the leader. Some completely different methods including a novel state decomposition and orthogonal decomposition are applied to overcome the difficulties caused by partially observability which improves the tools and relaxes the constraint condition imposed on admissible control in the existing literature. More precisely, the followers encounter the standard linear-quadratic partially observed optimal control problems, however, a kind of forward-backward indefinite linear-quadratic partially observed optimal control problem is considered by the leader. Instead of maximum principle of forward-backward control systems, inspired by the existing work related to definite case and classical forward control system, some distinct forward-backward linear-quadratic decoupling techniques including the method of completion of squares are applied to solve the leader's problem. More interestingly, we develop the deterministic formation control in multi-agent system with a framework of Stackelberg differential game and extend it to the stochastic case. The optimal strategies are obtained by our theoretical result suitably.
\end{abstract}

\maketitle

\section{Introduction}

Inspired by some natural phenomenon appearing in the swarming behavior of living beings, such as flocks of birds, schools of fish, colonies of bacteria and herds of wildebeest, the introduction of formation control of multi-agent systems gives a corresponding theoretical characteristic. The studies of formation control now becomes a hot frontier research topic because of its widespread applications to mobile robots, satellite control, unmanned surface vessels and unmanned aerial vehicles. Therefore, the aim of formation control is to design a controller  to drive all the agents to maintain the formation along the desired trajectories. For solving the formation control, there are now several approaches such as behavior-based approach, virtual-structure approach, leader-follower approach, consensus and optimal control and differential game. However, it is worth noting that most of investigations are related to the deterministic formation control in which only deterministic differential game approaches are applied. As far as we know, there is little literature that studies the stochastic case where necessary probability theories closed to forward backward stochastic differential equations (FBSDEs for short) are the key to model the stochastic framework of formation control. Meanwhile, the stochastic characteristic in Stackelberg differential game shows the mutual feedback between leader and follower. Therefore, instead of Nash differential game approach, the Stackelberg game including leader-follower mechanism is adapted to formation control well when consider one of the agent leading the formation direction and tracking the desired trajectory as the leader while other agents keeping the formation by following the leader with a fixed distance as the followers when they can exchange information such as relative positions and velocities in order to maintain a formation.

The Stackelberg game, also known as the leader-follower game, was introduced firstly by Stackelberg \cite{Stack52} where its hierarchical feature was first shown under the economic background that certain companies have advantages of domination over others in some markets environment. In the classical Stackelberg game, two agents played by asymmetric roles are usually considered in which the leader claims his strategy in advance, and the follower, considering the given leader's strategy, makes an instantaneous reaction (that can be regarded as a map from or feedback to the leader's control) by optimizing his cost functional. Then the leader would like to choose an optimal strategy to optimize his cost functional by taking the rational response of the follower into account. Since its meaningful structure and background, the leader-follower's feature mainly plays an important role in many aspects, such as principle-agent problems \cite{CZ13}, newsvendor/wholesaler problems \cite{OSU13}, optimal reinsurance problems \cite{CS18} and operations management and marketing channel problems \cite{LS17}. Now in our paper, the stochastic Stackelberg differential game can be applied in the formation control problems where the uncertainties coming from Brownian motions are considered and the stochastic model is studied in the martingale's environment.

When it comes to the Stackelberg game in the past decades, there have been a great deal of significant works among which we only mention a few closely related. In the stochastic case, Yong \cite{Yong02} studied an indefinite linear quadratic (LQ for short) leader-follower game with random coefficients and control-dependent diffusion, where FBSDEs and Riccati equations are applied to obtain the state feedback representation of open-loop Stackelberg equilibrium. Since then, Stackelberg differential games has been studied extensively in the literature such as in the fields of mean-field \cite{FW24, HWX21, LJZ19, MB18}, delay \cite{XSZ18, MS22}, jump \cite{Moon21} and regime switching \cite{Moon23, LXZ23}. In terms of the agents' different information sets, Bensoussan et al. \cite{BCS15} introduced several solution concepts, and obtained the maximum principle for the leader's problem under the adapted closed-loop memoryless information structure. Li and Yu \cite{LY18} studied the solvabilities of coupled  FBSDEs with a multilevel self-similar domination-monotonicity structure which was applied to solve a LQ generalized multilevel hierarchical leader-follower game and characterize its corresponding equilibrium in a closed form. Then it was extended to the jump-diffusion \cite{LXY21} and time-delay systems \cite{LW24}.

Different from the complete information background in previous literature, some of the agents with asymmetry of roles in reality can only get the partial information usually due to the delay in information transmission, inequality in market competition, private information in different class and limitation policy by the government. In particular, the major difficulty in the formation control problem shows that each agent is only able to exchange information with other agents according to the communication topology, which means the constrainted information limits many game strategy design approaches that usually consider individual agents to have global information. In this case, the formation control problem would be formulated in such a way that individual agents try to minimize their partially observed formation errors and to solve it from a stochastic Stackelberg differential game theory point of view. Based on the above analysis, therefore, the partially observed setting matches the leader-follower feature well, such problems have aroused intense interest \cite{SWX16, ZS22}, and there are also lots of works on this issue. Li et al. \cite{LMFZ23} investigated the LQ Gaussian Stackelberg game under a class of nested observation information, where the follower only use the observation data to design its strategy but the leader use the global observation data. Very recently, Wang and Wang \cite{WW24} studied a partially observed LQ stochastic Stackelberg game with two followers where both of followers can not observe the state process but leader can know the full information. Then Si and Shi \cite{SS24arXiv} studied an overlapping information LQ Stackelberg stochastic differential game with two leaders and two followers in which the noisy information available to the leaders and the followers may be asymmetric and have overlapping part. Especially, in the recent work \cite{ZS22}, based on the classical state decomposition technique \cite{Ben92} and backward separation principle developed in the partially observed control problem \cite{WW08, WWX15}, a kind of partially observed LQ Stackelberg game theory was investigated and was applied to a practical dynamic advertising problem. However, the conditions satisfied by control process are relatively strong and not natural enough in \cite{ZS22}. Moreover, the stochastic LQ optimal control of fully coupled FBSDEs would be encountered by the leader in most of the Stackelberg literature in which the maximum principle is the main tool to obtain the maximum condition satisfied by optimal control, then get the optimal state feedback by decoupling the Hamiltonian system. Few work considers solving the fully coupled forward backward stochastic LQ optimal control by the completion of squares methods but the recent work studied by \cite{HJX23}, and there is no any result on the counterpart with partially observable information.

The stochastic Stackelberg game theory mentioned in the above literature are mostly applied to investigate  in the fields of mathematical finance, management science and economy, but relatively few literature consider the multi-agent formation control such as unmanned aerial vehicle (UAV), robots and aircraft carrier formation corrupted by stochastic noise such as Brownian motion, in which the dynamic system is described by stochastic differential equations (SDEs for short) or FBSDEs instead of deterministic model. Different from the deterministic formation control, we introduce the notion of partially observed stochastic formation control whose acquirement need more probability and stochastic analysis technique, on one hand, the observation equation should be introduced and the formation control becomes a stochastic process adapted to the information filtration generated by observable process; On the other hand, stochastic filtering theory and solvability of FBSDEs are necessary, which add more difficulties than the deterministic case.

Therefore, the main contributions of our work in this paper can be summarized as follows:
\begin{enumerate}
\item[{\bf (i)}] A class of partially observed stochastic LQ Stackelberg differential game problem is investigated, where one leader and $N$ followers constitute a group of $N+1$ agents all of whom can only have access to noisy observation process related to the state. Motivated by a class of formation control, a distinguished partially observable setting \eqref{assumption_filtration} that the leader knows less than the followers is considered in our paper.

\item[{\bf (ii)}] In the existing literature, the state decomposition technique and backward separation principle has been widely applied to solve the partially observed LQ optimal control of SDEs and FBSDEs and plenty of excellent works are obtained. However, this comes at the cost of requiring stringent double-adaptedness conditions to be satisfied, which are necessary to break the circular dependence between the controlled filtration and the control process. Motivated by the idea in \cite{SX23}, we utilize the orthogonal decomposition technique of the state process to divide the original cost functional into two sub-cost functionals: one is the functional of control processes and filtering processes which can be regarded as the completely observed classical optimal control problem and another is independent of control process. However, in the followers' problem, the filtration is not only dependent on followers' control $u_1$ but also influenced by the leader's control $u_2$. For getting rid of the dependence on $u_2$, a novel state decomposition technique \eqref{decomposition}, \eqref{X^1Y^1} and \eqref{X^0Y^0} is introduced without any additional adaptedness imposed on control, which is distinguished from the existing decomposition in \cite{Ben92, WWX15}. In the new procedure studied in the leader's problem, as far as we know, it is the first time that we make an extension to the fully coupled forward backward stochastic systems by applying the orthogonal decomposition to the forward and backward state process and illustrating the independence between filtering error processes and filtrations. Therefore, we relax the limitation imposed on the control process. On one hand, the new state decomposition technique which divides the state equation into two still controlled equations, actually aims to introduce an $u_2$-independent sub-filtration to overcome the loop; On the other hand, the orthogonal decomposition idea used in our paper mainly focus on the separation of cost functionals instead of filtrations, the core is to first obtain the filtering equation for any fixed control, then determining the optimal control by the completion of squares method in the transformed completely observable case.

\item[{\bf (iii)}] For the leader's problem, instead of applying maximum principle to get a candidate control, we initially study the indefinite partially observed LQ  optimal control of fully coupled FBSDEs with the new methods combining with completion of squares method, idea of four-step scheme and dimension-enlargement motivated by \cite{HJX23} which is further extended by us to the multi-dimensional case with inhomogeneous terms.

\item[{\bf (iv)}] Besides the difference of methods from that applied in \cite{ZS22}, the multi-agent formation control is investigated in the stochastic case for the first time, where the leader and followers play a asymmetrical Stackelberg game and Nash game are played among these finite followers with symmetrical roles. Moreover, there exists some deterministic Stackelberg differential game approach to formation control, but in our work, we add the partially observation element into the background to characterize the asymmetry of leader-follower structure, so the partially observed stochastic formation control is closed to the stochastic filtering theory of FBSDEs. We obtain the optimal state filtering feedback form of Stackelberg equilibrium.
\end{enumerate}

The rest of the article is organized as follows. Section \ref{s2} presents some preliminaries and formulates the partially observed LQ stochastic Stackelberg differential game with $N$ followers. A LQ partially observed optimal control of the follower's problem is studied, the filtering equation of forward state is derived, then the optimal state filtering feedback forms of the followers' optimal control and optimal cost are obtained by introducing two Riccati equations and two backward stochastic differential equations (BSDEs for short) in Section \ref{follower pro}. In Section \ref{leader pro}, we investigate an indefinite LQ partially observed optimal control of fully coupled FBSDEs with conditional mean field terms, where the forward backward stochastic filtering equation is first obtained, then we complete the leader's problem by introducing a series of ordinary differential equations (ODEs for short), Riccati equations and BSDEs and we get the optimal state filtering feedback form of optimal control and its corresponding optimal cost. In Section \ref{s5}, for an application, we consider a class of multi-agent  formation control problem and extend it into the stochastic case, then it can be solved by the obtained partially observed Stackelberg game approach. Section \ref{s6} concludes the paper. In the Appendix, we give the mathematical analysis on the independence between filtering error process and filtration and investigate an LQ optimal control theory of fully coupled FBSDEs for the multi-dimensional case with inhomogeneous terms.

\section{Problem Formulation}\label{s2}

Let $(\Om,\cF,\PP)$ be a complete probability space on which two standard $l_1$-dimensional and $l_2$-dimensional Brownian motions $W^1=(W^{1,1}_t,W^{1,2}_t,\dots,W^{1,l_1}_t)^\top_{0\leq t\leq T}$ and $W^2=(W^{2,1}_t,W^{2,2}_t, $ $\dots,W^{2,l_2}_t)^\top_{0\leq t\leq T}$ are defined. Assume that $\FF=\{\cF_t,0\leq t\leq T\}$ is the $\PP$-augmentation of the natural filtration of $W^1, W^2$, where $\cF_0$ contains all $\PP$-null sets of $\cF$. Denote by $\RR^n$ the $n$-dimensional real Euclidean space, and let $\RR^{n\times m}$ be the set of $n\times m$ real matrixes. Let $\langle\cdot,\cdot\rangle$ (respectively,$\vert\cdot\vert$) denote the usual scalar product (respectively, usual norm) of $\RR^n$ and $\RR^{n\times m}$. The scalar product (respectively, norm) of matrices $M=(m_{ij}),N=(n_{ij})\in\RR^{n\times m}$ is denoted by $\langle M,N\rangle=\tr\{MN^\top\}$ (respectively, $\Vert M\Vert=\sqrt{MM^\top}$), where the superscript $\top$ denotes the transpose of vectors or matrices.

We introduce the following spaces. We denote by $\SS^n$ all $n\times n$ symmetric matrices, by $\SS^n_{+}$ the subspace of all nonnegative definite matrices of $\SS^n$, by $\SS^n_{-}$ the subspace of all negative definite matrices of $\SS^n$, by $L^\infty(0,T;\RR^{n\times m})$ the space of essential bounded measurable $\RR^{n\times m}$-valued functions.

We consider the following multi-dimensional state equations of the leader and the followers:
\begin{equation}\label{state SDE}
\left\{
\begin{aligned}
dX^{u_1,u_2}_t&=\big[A(t)X^{u_1,u_2}_t+B_1^\top(t) u_1(t)+B_2(t)u_2(t)+\alpha_t\big]dt+C_1(t)dW^1_t+C_2(t)dW^2_t,\\
X^{u_1,u_2}_0&=x_0,
\end{aligned}
\right.
\end{equation}
where for $t\in[0,T]$, $X_t\in\RR^n, A(t)\in\RR^{n\times n}$. And we denote $u_1=(u_{11},u_{12},\cdots,u_{1N})^\top\in\RR^{Nm}$ with $u_{1i}\in\RR^m$ for $i=1,\dots, N$, $B_1=(B_{11}^\top,B_{12}^\top,\cdots,B_{1N}^\top)^\top$ with $B_{1i}\in\RR^{n\times m}$ for $i=1,\dots,N$, $B_1^\top(t) u_1(t)=\sum_{i=1}^NB_{1i}(t)u_{1i}(t)$, $B_2(t)\in\RR^{n\times m}, \al(t)\in\RR^n, C_1(t)\in\RR^{n\times l_1}, C_2(t)\in\RR^{n\times l_2}$, and $W^1_t,W^2_t$ are both $\RR^{l_1}$ and $\RR^{l_2}$-valued mutually independent Brownian motions with independent components.

We suppose that the followers and the leader can not observe the above process $X$, instead, the only observed processes by them are described as follows, respectively:
\begin{equation*}
\begin{aligned}
dY^{u_1,u_2}_1(t)&=\big[f_1(t)X^{u_1,u_2}_t+g_1(t)\big]dt+K_1(t)dW^1_t,\ Y^{u_1,u_2}_1(0)=0,\\
dY^{u_1,u_2}_2(t)&=\big[f_2(t)X^{u_1,u_2}_t+g_2(t)\big]dt+K_2(t)dW^2_t,\ Y^{u_1,u_2}_2(0)=0,
\end{aligned}
\end{equation*}
where we assume the processes observed by all the followers are the same, i.e., $Y_{1i}^{u_1,u_2}=Y_1^{u_1,u_2}\in\RR^{l_1}, f_{1i}=f_1\in\RR^{l_1\times n},g_{1i}=g_1\in\RR^{l_1},K_{1i}=K_1\in\RR^{l_1\times l_1}$ for $i=1,2,\dots,N$, $Y_2^{u_1,u_2}\in\RR^{l_2}, f_2\in\RR^{l_2\times n},g_2\in\RR^{l_2},K_2\in\RR^{l_2\times l_2}$. We denote filtrations $\cF_t^{Y_i^{u_1,u_2}}$ by the $\sigma$-algebra generated by $Y_i^{u_1,u_2}$ up to $t$, respectively, for $i=1,2$, for $t\in[0,T]$ and first assume that $\cF_t^{Y_2^{u_1,u_2}}\subset\cF_t^{Y_1^{u_1,u_2}}$, which then will be given in detail in term of inclusion relationship. It should be noted that the drift term can contain the state variable in the leader's observation equation, which is an extension of that in \cite{ZS22} because of the new technique we applied in this paper.

The cost functionals of the followers and the leader are defined as follows:
\begin{equation*}
\begin{aligned}
J_{1i}(x_0;u_1,u_2)&=\EE\bigg\{\int_0^T\bigg[\bigg\langle\begin{pmatrix}Q_{1i}&S_{1i}^\top\\S_{1i}&R_{1i}\end{pmatrix}\begin{pmatrix}X\\u_{1i}\end{pmatrix},\begin{pmatrix}X\\u_{1i}\end{pmatrix}\bigg\rangle +2\bigg\langle\begin{pmatrix}q_{1i}\\r_{1i}\end{pmatrix},\begin{pmatrix}X\\u_{1i}\end{pmatrix}\bigg\rangle \bigg]dt\\
&\qquad +\langle G_{1i}X_T,X_T\rangle+2\langle g_{1i},X_T\rangle\bigg\},\ i=1,2,\dots,N,\\
J_2(x_0;u_1,u_2)&=\EE\bigg\{\int_0^T\bigg[\bigg\langle\begin{pmatrix}Q_2&S_2^\top\\S_2&R_2\end{pmatrix}\begin{pmatrix}X\\u_2\end{pmatrix},\begin{pmatrix}X\\u_2\end{pmatrix}\bigg\rangle +2\bigg\langle\begin{pmatrix}q_2\\r_2\end{pmatrix},\begin{pmatrix}X\\u_2\end{pmatrix}\bigg\rangle \bigg]dt\\
&\qquad +\langle G_2X_T,X_T\rangle+2\langle g_2,X_T\rangle\bigg\},
\end{aligned}
\end{equation*}
which are minimized for the purpose of this paper, by choosing the corresponding optimal controls $\bar{u}_1$ and $\bar{u}_2$ in the admissible control sets $\cU^F_{ad}$ and $\cU^L_{ad}$ defined by
\begin{equation*}
\begin{split}
\cU^F_{ad}=&\bigg\{u_1:[0,T]^N\times\Om\rightarrow\RR^{Nm}\bigg|u_1(t)=(u_{11}(t),\cdots,u_{1N}(t))^\top\in\mathcal{F}_t^{Y_1^{u_1,u_2}}\\
&\qquad\text{and}\ \mathbb{E}\int_0^T|u_1(t)|^2dt<\infty\bigg\},\\
\cU^L_{ad}=&\bigg\{u_2:[0,T]\times\Om\rightarrow\RR\bigg|u_2(t)\in\mathcal{F}_t^{Y_2^{u_1,u_2}}\text{and}\ \mathbb{E}\int_0^T|u_2(t)|^2dt<\infty\bigg\},
\end{split}
\end{equation*}
respectively.

Now we need the following assumptions:
\begin{enumerate}
\item[{\bf (A1)}] $A(\cdot)\in L^\infty(0,T;\RR^{n\times n}),\ B_{1i}(\cdot), B_2(\cdot)\in L^\infty(0,T;\RR^{n\times m}),\ \text{for}\ i=1,\dots, N,\\
 \al(\cdot)\in L^\infty(0,T;\RR^n),\ C_1(\cdot)\in L^\infty(0,T;\RR^{n\times l_1}),\ C_2(\cdot)\in L^\infty(0,T;\RR^{n\times l_2})$;
\item[{\bf (A2)}] $f_1(\cdot)\in L^\infty(0,T;\RR^{l_1\times n}),\ f_2(\cdot)\in L^\infty(0,T;\RR^{l_2\times n}),\ g_1(\cdot)\in L^\infty(0,T;\RR^{l_1}),\\
 g_2(\cdot)\in L^\infty(0,T;\RR^{l_2}),\ K_1(\cdot)\in L^\infty(0,T;\RR^{l_1\times l_1}),\ K_2(\cdot)\in L^\infty(0,T;\RR^{l_2\times l_2})$;
\item[{\bf (A3)}] $Q_{1i}(\cdot)\in L^\infty(0,T;\SS^n_+),\ S_{1i}(\cdot)\in L^\infty(0,T;\RR^{n\times m}),\ R_{1i}(\cdot)\in L^\infty(0,T;\SS^m_+),\\
 q_{1i}(\cdot)\in L^\infty(0,T;\RR^{n}),\ r_{1i}(\cdot)\in L^\infty(0,T;\RR^m),\ \text{for}\ i=1,\dots, N,\\
 Q_2(\cdot)\in L^\infty(0,T;\SS^n_+), S_2(\cdot)\in L^\infty(0,T;\RR^{n\times m}),\ R_2(\cdot)\in L^\infty(0,T;\SS^m_{-}),\\
 q_2(\cdot)\in L^\infty(0,T;\RR^n),\ r_2(\cdot)\in L^\infty(0,T;\RR^m)$.
\end{enumerate}

\section{The Followers' Problem}\label{follower pro}

In this section, we consider a linear-quadratic partially observed optimal control problem by the technique in \cite{SX23}. 
In order to overcome the difficulty that the control process is adapted to the controlled filtration, for a fixed admissible control $u_1\in\cU^F_{ad}$, first, we consider decomposing the state and followers' observation equations as follows:
\begin{equation}\label{decomposition}
X^{u_1,u_2}=X^{1,u_2}+X^{0,u_1},\ Y^{u_1,u_2}_1=Y^{1,u_2}_1+Y^{0,u_1}_1,
\end{equation}
where $X^{1,u_2},X^{0,u_1},Y^{1,u_2}_1$ and $Y^{0,u_1}_1$ satisfiy
\begin{equation}\label{X^1Y^1}
\left\{
\begin{aligned}
dX^{1,u_2}_t&=\big[A(t)X^{1,u_2}_t+B_2(t)u_2(t)+\alpha_t\big]dt,\quad X^{1,u_2}_0=0,\\
dY^{1,u_2}_1(t)&=\big[f_1(t)X^{1,u_2}_t+g_1(t)\big]dt,\quad Y^{1,u_2}_1(0)=0,
\end{aligned}
\right.
\end{equation}
\begin{equation}\label{X^0Y^0}
\left\{
\begin{aligned}
dX^{0,u_1}_t&=\big[A(t)X^{0,u_1}_t+B_1^\top(t) u_1(t) \big]dt+C_1(t)dW^1_t+C_2(t)dW^2_t,\quad X^{0,u_1}_0=x_0,\\
dY^{0,u_1}_1(t)&=f_1(t)X^{0,u_1}_tdt+K_1(t)dW^1_t,\quad Y^{0,u_1}_1(0)=0,
\end{aligned}
\right.
\end{equation}
respectively. Then, we give the following detailed inclusion between the filtrations
\begin{equation}\label{assumption_filtration}
\cF_t^{Y_2^{u_1,u_2}}\subset\Big(\cF_t^{Y_1^{u_1,u_2}}\wedge\cF_t^{Y_1^{0,u_1}}\Big),
\end{equation}
where the filtration $\cF_t^{Y_1^{0,u_1}}$ is denoted by the $\sigma$-algebra generated by $Y_1^{0,u_1}$ up to $t$. Now we will give the auxiliary result for obtaining the filtering equation in the following.
\begin{lemma}\label{equivalent_filtration}
Given $u_2(\cdot)\in\cU^L_{ad}$, for fixed followers' control $u_1(\cdot)\in\cU^F_{ad}$ and let assumption \eqref{assumption_filtration} hold. Then $\cF_t^{Y_1^{u_1,u_2}}=\cF_t^{Y_1^{0,u_1}}$.
\end{lemma}
\begin{proof}
On one hand, since $u_2(\cdot)$ is $\cF_t^{Y_2^{u_1,u_2}}\Big(\subset\cF_t^{Y_1^{0,u_1}}\Big)$-adapted, it follows from the first equation in \eqref{X^1Y^1} that $X^{1,u_2}_{\cdot}$ is $\cF_t^{Y_1^{0,u_1}}$-adapted, so is $Y^{1,u_2}_1(\cdot)$. Then $Y^{u_1,u_2}_1=Y^{1,u_2}_1+Y^{0,u_1}_1$ is $\cF_t^{Y_1^{0,u_1}}$-adapted, i.e., $\cF_t^{Y_1^{u_1,u_2}}\subset\cF_t^{Y_1^{0,u_1}}$.

On the other hand, since $u_2(\cdot)$ is $\cF_t^{Y_2^{u_1,u_2}}\Big(\subset\cF_t^{Y_1^{u_1,u_2}}\Big)$-adapted, it follows from the first equation in \eqref{X^1Y^1} that $X^{1,u_2}_{\cdot}$ is $\cF_t^{Y_1^{u_1,u_2}}$-adapted, so is $Y^{1,u_2}_1(\cdot)$. Then $Y^{0,u_1}_1=Y^{u_1,u_2}_1-Y^{1,u_2}_1$ is $\cF_t^{Y_1^{u_1,u_2}}$-adapted, i.e., $\cF_t^{Y_1^{0,u_1}}\subset\cF_t^{Y_1^{u_1,u_2}}$.
\end{proof}
\begin{remark}
It is necessary to stress the novelty of the decomposition method \eqref{decomposition}, \eqref{X^1Y^1} and \eqref{X^0Y^0}. The state decomposition technique is a critical idea to overcome the circular dependence between the control process and control-dependent filtration, see \cite{Ben92} and \cite{WWX15}. However, the cost is an addition of strong limit that the admissible control have to be adapted to both original filtration and uncontrolled filtration simultaneously. Then the adaptedness requirement are also applied to solve the Stackelberg game case \cite{ZS22} where two control processes are considered. Fortunately, in our new decomposition \eqref{decomposition}, \eqref{X^1Y^1} and \eqref{X^0Y^0}, it is enough to separate the controls $(u_1,u_2)$ into one single sub-system, respectively, instead of putting both controls into the same equation labeled 1 like that in \cite{ZS22}. The advantage is that the equivalence of filtrations in Lemma \ref{equivalent_filtration} can be obtained more naturally with the strong adaptedness condition relaxed. The main reason of this setting counts on the feasibility of partially observed controlled sub-system \eqref{X^0Y^0} when fixing $u_1$ solved by \cite{SX23}, which is an interesting finding in our paper. This is a combination of decomposition technique \cite{Ben92} and orthogonal decomposition \cite{SX23}, which benefits for not only getting rid of the dependence of followers' filtration $\cF_t^{Y_1^{u_1,u_2}}$ on the leader's control $u_2$ but also relaxing the double adaptedness requirement. Moreover, the new decomposition technique can be further applied to study more similar game structure such as the partially observed large-population Stackelberg games.
\end{remark}
Based on the Lemma \ref{equivalent_filtration}, for a fixed admissible control $u_1\in\cU^F_{ad}$, we define the following notations:
\begin{equation*}
\hat{X}^{u_2}_t\triangleq\EE\Big[X^{u_2}_t\Big|\cF_t^{Y_1^{u_1,u_2}}\Big]=\EE\Big[X^{u_2}_t\Big|\cF_t^{Y_1^{0,u_1}}\Big],\ \tilde{X}_t\triangleq X^{u_2}_t-\hat{X}^{u_2}_t
\end{equation*}
in which the filtration $\cF_t^{Y_1^{0,u_1}}$ is independent of leader's control $u_2$ and it is benefit to obtain the filtering equation for $\hat{X}^{u_2}_t$, and define a stochastic process
\begin{equation*}
V_t\triangleq Y_1^{u_2}(t)-\int_0^t\big[f_1(s)\hat{X}^{u_2}_s+g_1(s)\big]ds.
\end{equation*}
Note that the process $\tilde{X}_\cdot$ is indeed independent of $u_2$ due to the $\cF_t^{Y_1^{0,u_1}}$-adaptedness of $u_2$ and the inclusion relation of filtrations, and so is $V$. The detailed proof verification would be provided strictly later. For any given admissible control $u_2$ and a fixed admissible control $u_1\in\cU_{ad}^F$, we can obtain
\begin{equation*}
\check{V}_t\triangleq\int_0^tK_1^{-1}(s)dV_s
\end{equation*}
is a standard $\cF_t^{Y_1^{0,u_1}}$-Brownian motion. Indeed, first, we have $\check{V}$ is a continuous $\cF_t^{Y_1^{0,u_1}}$-adapted, integrable process, and
\begin{equation*}
\begin{aligned}
d\check{V}_t=K^{-1}_1(t)dV_t=K^{-1}_1(t)f_1(t)\tilde{X}_tdt+dW^1_t.
\end{aligned}
\end{equation*}
Thus, by Lemma \ref{equivalent_filtration}, for $0\leq s<t\leq T$,
\begin{equation*}
\begin{aligned}
&\EE\Big[\check{V}_t-\check{V}_s\Big|\cF_s^{Y_1^{u_1,u_2}}\Big]=\EE\Big[W^1_t-W^1_s\Big|\cF_s^{Y_1^{0,u_1}}\Big]+\EE\Big[\int_s^tK^{-1}_1(r)f_1(r)\tilde{X}_rdr\Big|\cF_s^{Y_1^{0,u_1}}\Big]\\
&=\EE\Big[\EE[W^1_t-W^1_s|\cF_s]\Big|\cF_s^{Y_1^{0,u_1}}\Big]+\int_s^tK^{-1}_1(r)f_1(r)\EE\tilde{X}_rdr=0,
\end{aligned}
\end{equation*}
which implies that $\check{V}$ is a $\cF_t^{Y_1^{0,u_1}}$-martingale. Meanwhile,
\begin{equation*}
\begin{aligned}
d\big[\check{V}_t(\check{V}_t)^\top\big]&=K^{-1}_1(t)f_1(t)\tilde{X}_t(\check{V}_t)^\top dt+(dW^1_t)(\check{V}_t)^\top\\
&\quad +\check{V}_tK^{-1}_1(t)f_1(t)(\tilde{X}_t)^\top dt+\check{V}_t(dW^1_t)^\top+I_{l_1}dt.
\end{aligned}
\end{equation*}
Due to
\begin{equation*}
\begin{aligned}
&\EE\Big[\int_s^t\check{V}_r(K_1^{-1}(r)f_1(r)\tilde{X}_r)^\top dr\Big|\cF_s^{Y_1^{u_1,u_2}}\Big]\\
&=\int_s^t\EE\Big[\EE\big[\check{V}_r(\tilde{X}_r)^\top\big|\cF_r^{Y_1^{0,u_1}}\big]\Big|\cF_s^{Y_1^{0,u_1}}\Big]f^\top_1(r)(K_1^{-1}(r))^\top dr\\
&=\int_s^t\EE\Big[\check{V}_r\EE(\tilde{X}_r)^\top\Big|\cF_s^{Y_1^{0,u_1}}\Big]f^\top_1(r)(K_1^{-1}(r))^\top dr=0,
\end{aligned}
\end{equation*}
and
\begin{equation*}
\begin{aligned}
\EE\Big[\int_s^t\check{V}_rd(W^1_r)^\top\Big|\cF_s^{Y_1^{u_1,u_2}}\Big]=\EE\Big[\EE\big[\int_s^t\check{V}_rd(W^1_r)^\top\big|\cF_s\big]\Big|\cF_s^{Y_1^{0,u_1}}\Big]=0,
\end{aligned}
\end{equation*}
we have
\begin{equation*}
\EE\Big[(\check{V}_t)^2-(\check{V}_s)^2\Big|\cF_s^{Y_1^{u_1,u_2}}\Big]=(t-s)I_{l_1},
\end{equation*}
which implies that $\check{V}$ is a standard Brownian motion.

Next, we consider the following process
\begin{equation*}
\Lambda^{1}_t\triangleq\hat{X}^{u_2}_t-x_0-\int_0^t\big[A(s)\hat{X}^{u_2}_s+B_1^\top(s)u_1(s)+B_2(s)u_2(s)+\alpha(s)\big]ds.
\end{equation*}
In the following, we will prove that $\Lambda^{1}\equiv\{\Lambda^{1}_t;0\leq t\leq T\}$ is a square-integrable $\{\cF_t^{Y_1^{u_1,u_2}}\}$-martingale having right continuous with left limits (RCLL) paths, and $\Lambda_0=0$, a.s.. Indeed, for fixed $0\leq s<t\leq T$,
\begin{equation*}
\begin{aligned}
&\EE\Big[\laa^{1}_t-\laa^{1}_s\Big|\cF_s^{Y_1^{u_1,u_2}}\Big]=\EE\Big[\laa^{1}_t-\laa^{1}_s\Big|\cF_s^{Y_1^{0,u_1}}\Big]\\
&=\EE\Big[\EE\big[X^{u_2}_t\big|\cF_t^{Y_1^{0,u_1}}\big]-\EE\big[X^{u_2}_s\big|\cF_s^{Y_1^{0,u_1}}\big]\Big|\cF_s^{Y_1^{0,u_1}}\Big]\\
&\quad -\int_s^t\bigg\{A(r)\EE\Big[\EE\big[X^{u_2}_r\big|\cF_r^{Y_1^{0,u_1}}\big]\Big|\cF_s^{Y_1^{0,u_1}}\Big]+B_1^\top(r)\EE\Big[u_1(r)\Big|\cF_s^{Y_1^{0,u_1}}\Big]\\
&\qquad\qquad +B_2(r)\EE\Big[u_2(r)\Big|\cF_s^{Y_1^{0,u_1}}\Big]+\al(r)\bigg\}dr\\
&=\EE\Big[X^{u_2}_t-X^{u_2}_s-\int_s^t\big\{A(r)X^{u_2}_r+B_1^\top(r)u_1(r)+B_2(r)u_2(r)+\al(r)\big\}dr\Big|\cF_s^{Y_1^{0,u_1}}\Big]\\
&=\EE\Big[\EE\Big[\int_s^tC_1(r)dW^1_r+\int_s^tC_2(r)dW^2_r\Big|\cF_s\Big]\Big|\cF_s^{Y_1^{0,u_1}}\Big]=0.
\end{aligned}
\end{equation*}
Now by the theorem of Fujisaki, Kallianpar and Kunita \cite{RW00}, there exists a square-integrable $\{\cF_t^{Y_1^{0,u_1}}\}$-progressively measurable process $\la^{1}=\{\la^{1}_t;0\leq t\leq T\}\in\RR^{n\times l_1}$ such that
\begin{equation*}
\laa^{1}_t=\int_0^t\la^{1}_s(K_1(s)K_1^\top(s))^{-1} dV_s,\ 0\leq t\leq T.
\end{equation*}
To determine the process $\la^{1}$, we let $\zeta=\{\zeta_t;0\leq t\leq T\}\in\RR^{n\times l_1}$ be a fixed but arbitrary square integrable $\{\cF_t^{Y_1^{0,u_1}}\}$-progressively measurable process. Consider the $\cF_t^{Y_1^{0,u_1}}$-martingale
\begin{equation*}
\eta_t\triangleq\int_0^t\zeta_s(K_1(s)K_1^\top(s))^{-1} dV_s,\ 0\leq t\leq T,
\end{equation*}
and then
\begin{equation*}
\EE\big[\laa^{1}_t(\eta_t)^\top\big]=\EE\int_0^t\la^{1}_s(K_1(s)K_1^\top(s))^{-1}(\zeta_s)^\top ds.
\end{equation*}
On the other hand, we also have
\begin{equation*}
\begin{aligned}
&\EE\big[\laa^{1}_t(\eta_t)^\top\big]=\EE\big[\hat{X}^{u_2}_t(\eta_t)^\top\big]-\EE\int_0^t\big[A(s)\hat{X}^{u_2}_s+B_1^\top(s)u_1(s)+B_2(s)u_2(s)+\alpha(s)\big](\eta_t)^\top ds.
\end{aligned}
\end{equation*}
Since for $0\leq s\leq t\leq T$,
\begin{equation*}
\begin{aligned}
\EE\big[\hat{X}^{u_2}_s(\eta_t)^\top\big]&=\EE\Big[\EE\big[\hat{X}^{u_2}_s(\eta_t)^\top\big|\cF_s^{Y_1^{0,u_1}}\big]\Big]=\EE\big[\hat{X}^{u_2}_s(\eta_s)^\top\big]\\
&=\EE\Big[\EE\big[X^{u_2}_s\big|\cF_s^{Y_1^{0,u_1}}\big](\eta_s)^\top\Big]=\EE\big[X^{u_2}_s(\eta_s)^\top\big],
\end{aligned}
\end{equation*}
and
\begin{equation*}
\begin{aligned}
&\EE\big[\big(B_1^\top(s)u_1(s)+B_2(s)u_2(s)+\al(s)\big)(\eta_t)^\top\big]=\EE\Big[\sum_{i=1}^NB_{1i}(s)\EE\big[u_{1i}(s)(\eta_t)^\top\big|\cF_s^{Y_1^{0,u_1}}\big]\Big]\\
&\qquad+\EE\Big[B_2(s)\EE\big[u_2(s)(\eta_t)^\top\big|\cF_s^{Y_1^{0,u_1}}\big]\Big]+\EE\Big[\al(s)\EE\big[(\eta_t)^\top\big|\cF_s^{Y_1^{0,u_1}}\big]\Big]\\
&\quad=\EE\big[B_1^\top(s)u_1(s)(\eta_s)^\top+B_2(s)u_2(s)(\eta_s)^\top+\al(s)(\eta_s)^\top\big],
\end{aligned}
\end{equation*}
then we obtain
\begin{equation}\label{Lambda eta}
\begin{aligned}
&\EE\big[\laa^{1}_t(\eta_t)^\top\big]=\EE\big[X^{u_2}_t(\eta_t)^\top\big]\\
&\quad -\EE\int_0^t\big[A(s)X^{u_2}_s+B_1^\top(s)u_1(s)+B_2(s)u_2(s)+\alpha(s)\big](\eta_s)^\top ds.
\end{aligned}
\end{equation}
Observing that
\begin{equation*}
d\eta_t=\zeta_t(K_1(t)K_1^\top(t))^{-1} f_1(t)\tilde{X}_tdt+\zeta_tK_1^{-1}(t)dW^1_t,
\end{equation*}
and using It\^o's formula, it has
\begin{equation*}
\begin{aligned}
&\EE\big[X^{,u_2}_t(\eta_t)^\top\big]=\EE\int_0^t\Big\{X^{u_2}_s\big[\zeta_s(K_1(s)K_1^\top(s))^{-1} f_1(s)\tilde{X}_s\big]^\top\\
&+\big[A(s)X^{u_2}_s+B_1^\top(s)u_1(s) +B_2(s)u_2(s)+\al(s)\big](\eta_s)^\top+C_1(s)\big(\zeta_sK_1^{-1}(s)\big)^\top\Big\}ds,
\end{aligned}
\end{equation*}
where
\begin{equation}\label{deduction_1}
\begin{aligned}
&\EE\int_0^tX^{u_2}_s(\tilde{X}_s)^\top f^\top_1(s)(K_1(s)K_1^\top(s))^{-1}(\zeta_s)^\top ds\\
&=\EE\int_0^t\hat{X}^{u_2}_s(\tilde{X}_s)^\top f^\top_1(s)(K_1(s)K_1^\top(s))^{-1}(\zeta_s)^\top ds\\
&\quad +\EE\int_0^t\tilde{X}_s(\tilde{X}_s)^\top f^\top_1(s)(K_1(s)K_1^\top(s))^{-1}(\zeta_s)^\top ds\\
&=\EE\int_0^t\EE\Big[\hat{X}^{u_2}_s(\tilde{X}_s)^\top\Big|\cF_s^{Y_1^{0,u_1}}\Big] f^\top_1(s)(K_1(s)(K_1(s)K_1^\top(s))^{-1}(\zeta_s)^\top ds\\
&\quad +\EE\int_0^t\EE\Big[\tilde{X}_s(\tilde{X}_s)^\top\Big|\cF_s^{Y_1^{0,u_1}}\Big] f^\top_1(s)(K_1(s)K_1^\top(s))^{-1}(\zeta_s)^\top ds\\
&=\EE\int_0^t\EE\big[\tilde{X}_s(\tilde{X}_s)^\top\big] f^\top_1(s)(K_1(s)K_1^\top(s))^{-1}(\zeta_s)^\top ds\\
&=\EE\int_0^t\Sigma_sf^\top_1(s)(K_1(s)K_1^\top(s))^{-1}(\zeta_s)^\top ds,
\end{aligned}
\end{equation}
and we set $\Sigma_s\triangleq \EE\big[\tilde{X}_s(\tilde{X}_s)^\top\big]$, which implies that
\begin{equation*}
\begin{aligned}
&\EE\big[X^{u_2}_t(\eta_t)^\top\big]=\EE\int_0^t\big[A(s)X^{u_2}_s+B_1^\top(s)u_1(s)+B_2(s)u_2(s)+\alpha(s)\big](\eta_s)^\top ds\\
&\quad +\EE\int_0^t\big[\Sigma_sf^\top_1(s)+C_1(s)K_1^\top(s)\big](K_1(s)K_1^\top(s))^{-1}(\zeta_s)^\top ds.
\end{aligned}
\end{equation*}
The independence of $\tilde{X}_t$ and $\cF^{Y^{u_1,u_2}_1}_t$ has been applied in the above deduction \eqref{deduction_1}. For detailed proof, please see the Appendix \ref{Appendix_A}.

Then making the comparison with \eqref{Lambda eta} we get
\begin{equation*}
\EE\big[\laa^{1}_t(\eta_t)^\top\big]=\EE\int_0^t\big[\Sigma_sf^\top_1(s)+C_1(s)K_1^\top(s)\big](K_1(s)K_1^\top(s))^{-1}(\zeta_s)^\top ds.
\end{equation*}
So we obtain
\begin{equation}\label{lambda}
\la^{1}_s=\Sigma_sf^\top_1(s)+C_1(s)K_1^\top(s),\ 0\leq s\leq T.
\end{equation}
The filtering process $\hat{X}^{u_2}$ satisfies
\begin{equation}\label{hat X}
\begin{aligned}
d\hat{X}^{u_2}_t&=\big[A(t)\hat{X}^{u_2}_t+B_1^\top(t)u_1(t)+B_2(t)u_2(t)+\al(t)\big]dt\\
&\qquad +\big[\Sigma_tf^\top_1(t)+C_1(t)K_1^\top(t)\big](K_1(t)K_1^\top(t))^{-1}dV_t\\
&=\big[A(t)\hat{X}^{u_2}_t+B_1^\top(t)u_1(t)+B_2(t)u_2(t)+\al(t)\big]dt\\
&\qquad +\big[\Sigma_tf^\top_1(t)+C_1(t)K_1^\top(t)\big](K_1(t)K^\top_1(t))^{-1}f_1(t)\tilde{X}_tdt\\
&\qquad +\big[\Sigma_tf^\top_1(t)+C_1(t)K_1^\top(t)\big]K_1^{-1}(t)dW^1_t,
\end{aligned}
\end{equation}
and $\tilde{X}$ satisfies
\begin{equation}\label{tilde X}
\begin{aligned}
d\tilde{X}_t&=\big\{A(t)-\big[\Sigma_tf^\top_1(t)+C_1(t)K_1^\top(t)\big](K_1(t)K_1^\top(t))^{-1}f_1(t)\big\}\tilde{X}_tdt\\
&\quad -\Sigma_tf^\top_1(t)K_1^{-1}(t)dW^1_t+C_2(t)dW^2_t.
\end{aligned}
\end{equation}
Then apply It\^o's formula to $\tilde{X}(\tilde{X})^\top$, we have
\begin{equation}\label{Sigma}
\begin{aligned}
\Sigma_t&=\int_0^t\Big[\big(A(s)-C_1(s)K_1^{-1}(s)f_1(s)\big)\Sigma_s+\Sigma_s\big(A(s)-C_1(s)K_1^{-1}(s)f_1(s)\big)^\top\\
&\qquad -\Sigma_sf^\top_1(s)(K_1(s)K_1^\top(s))^{-1}f_1(s)\Sigma_s+C_2(s)C_2^\top(s)\Big]ds.
\end{aligned}
\end{equation}

Next, based on the orthogonal decomposition of the state process \cite{SX23}, we divide the cost functional $J_{1i}(x_0;u_1,u_2)$ for the $i$th follower as follows:
\begin{equation*}
\begin{aligned}
J_{1i}(x_0;u_1,u_2)
&=\EE\bigg[\int_0^T\Big(\big\langle Q_{1i}(t)\tilde{X}_t,\tilde{X}_t\rangle+2\langle q_{1i}(t),\tilde{X}_t\big\rangle\Big) dt+\big\langle G_{1i}\tilde{X}_T,\tilde{X}_T\big\rangle+2\big\langle g_{1i},\tilde{X}_T\big\rangle\bigg]\\
&\quad+\EE\bigg[\int_0^T\bigg[\bigg\langle\begin{pmatrix}Q_{1i}&S_{1i}^\top\\S_{1i}&R_{1i}\end{pmatrix}
  \begin{pmatrix}\hat{X}^{u_2}\\u_{1i}\end{pmatrix},
  \begin{pmatrix}\hat{X}^{u_2}\\u_{1i}\end{pmatrix}\bigg\rangle\\
&\qquad +2\bigg\langle\begin{pmatrix}q_{1i}\\r_{1i}\end{pmatrix},\begin{pmatrix}\hat{X}^{u_2}\\u_{1i}\end{pmatrix}\bigg\rangle \bigg]dt
  +\big\langle G_{1i}\hat{X}^{u_2}_T,\hat{X}^{u_2}_T\big\rangle+2\big\langle g_{1i},\hat{X}^{u_2}_T\big\rangle\bigg]\\
&\triangleq\tilde{J}_{1i}(x_0;u_1,u_2)+\hat{J}_{1i}(x_0;u_1,u_2).
\end{aligned}
\end{equation*}
We introduce Riccati equations:
\begin{equation}\label{P ^ 1i}
\left\{
\begin{aligned}
&\dot{P}^{1i}_t+A^\top(t)P^{1i}_t+P^{1i}_tA(t)-\big(P^{1i}_tB_{1i}(t)+S^\top_{1i}(t)\big)R_{1i}^{-1}(t)\big(B^\top_{1i}(t)P^{1i}_t+S_{1i}(t)\big)+Q_{1i}(t)=0,\\
&P^{1i}_T=G_{1i},
\end{aligned}
\right.
\end{equation}
and BSDEs:
\begin{equation}\label{e.3.14}
\left\{
\begin{aligned}
-d\vp^{1i,u_2}_t&=\Big\{\big[A(t)-B_{1i}(t)R_{1i}^{-1}(t)\big(B^\top_{1i}(t) P^{1i}_t+S_{1i}(t)\big)\big]^\top\vp^{1i,u_2}_t+P^{1i}_t\al(t)\\
&\qquad +q_{1i}(t)-\big[R_{1i}^{-1}(t)\big(B^\top_{1i}(t)P^{1i}_t+S_{1i}(t)\big)\big]^\top r_{1i}(t)+P^{1i}_tB_2(t)u_2(t)\Big\}dt\\
&\quad -\la^{1i,u_2}_td\check{V}_t,\\
\vp^{1i,u_2}_T&=g_{1i},\quad i=1,\cdots,N.
\end{aligned}
\right.
\end{equation}
Apply It\^o's formula to $\langle P^{1i}_t\hat{X}^{u_2}_t,\hat{X}^{u_2}_t\rangle$ and $\langle\vp^{1i,u_2}_t,\hat{X}^{u_2}_t\rangle$ on $[0,T]$, respectively, and insert them into $\hat{J}_{1i}(x_0;u_1,u_2)$, then we have (We usually omit the time variable $t$ in the following for simplicity.)
\begin{equation*}
\begin{aligned}
&\hat{J}_{1i}(x_0;u_1,u_2)=\EE\big[\langle P^{1i}_0x_0+2\vp^{1i,u_2}_0,x_0\rangle\big]+\EE\int_0^T\Big[\big\langle R_{1i}\big[u_{1i}+R_{1i}^{-1}(B^\top_{1i}P^{1i}+S_{1i})\hat{X}^{u_2}\\
&\quad +R_{1i}^{-1}(B_{1i}^\top\vp^{1i,u_2}+r_{1i})\big],\big[u_{1i}+R_{1i}^{-1}(B^\top_{1i}P^{1i}+S_{1i})\hat{X}^{u_2}+R_{1i}^{-1}(B_{1i}^\top\vp^{1i,u_2}+r_{1i})]\big\rangle\\
&\quad -\Big\langle R_{1i}^{-1}(B_{1i}^\top\vp^{1i,u_2}+r_{1i}),(B_{1i}^\top\vp^{1i,u_2}+r_{1i})\Big\rangle+\sum_{j=1}^{l_1}\big\langle P^{1i}(\Delta_j+C_{1j}),\Delta_j+C_{1j}\big\rangle\\
&\quad +2\big\langle B_2u_2,\vp^{1i,u_2}\big\rangle+2\big\langle\vp^{1i,u_2},\al\big\rangle+2\big\langle\Delta+C_1,\la^{1i,u_2}\big\rangle\Big] dt,
\end{aligned}
\end{equation*}
where $\Delta(t)\equiv(\Delta_1(t),\dots,\Delta_{l_1}(t))\triangleq\Sigma_t \big(K_1^{-1}(t)f_1(t)\big)^\top$.
Similarly, we also introduce ODEs:
\begin{equation*}
\left\{
\begin{aligned}
&\dot{\Pi}^{1i}_t+\big[A-(\Sigma f^\top_1+C_1K_1^\top)(K_1K_1^\top)^{-1}f_1)\big]^\top\Pi^{1i}_t\\
&\qquad+\Pi^{1i}_t\big[A-(\Sigma f^\top_1+C_1K_1^\top)(K_1K_1^\top)^{-1}f_1\big]+Q_{1i}=0,\ i=1,\cdots,N\\
&\Pi^{1i}_T=G_{1i},
\end{aligned}
\right.
\end{equation*}
and BSDEs:
\begin{equation*}
\left\{
\begin{aligned}
-d\pi_t^{1i}&=\Big\{\big[A-(\Sigma f^\top_1+C_1K_1^\top)(K_1K_1^\top)^{-1}f_1\big]^\top\pi^{1i}_t+q_{1i}\Big\}dt-\beta^1_tdW^1_t-\beta^2_tdW^2_t,\\
\pi_T^{1i}&=g_{1i},\quad i=1,\cdots,N.
\end{aligned}
\right.
\end{equation*}
Apply It\^o's formula to $\langle\Pi_t^{1i}\tilde{X}_t,\tilde{X}_t\rangle$ and $\langle\pi^{1i}_t,\tilde{X}_t\rangle$ on $[0,T]$, respectively, and insert them into $\tilde{J}_{1i}(x_0;u_1,u_2)$, then we have
\begin{equation*}
\begin{aligned}
\tilde{J}_{1i}(x_0;u_1,u_2)&=\EE\int_0^T\bigg[\sum_{j=1}^{l_1}\langle\Pi^{1i}\Delta_j,\Delta_j\rangle+\sum_{j=1}^{l_2}\langle\Pi^{1i}C_{2j},C_{2j}\rangle\\
&\qquad\qquad -2\langle\Sigma f^\top_1(K_1^\top)^{-1},\beta^1\rangle+2\langle C_2,\beta^2\rangle\bigg]dt.
\end{aligned}
\end{equation*}

Now, the optimal filtering feedback form of the $i$th follower's optimal control is
\begin{equation}\label{bar u1}
\bar{u}_{1i}(t)=-R_{1i}^{-1}(t)\big[B_{1i}^\top(t) P^{1i}_t+S_{1i}(t)\big]\hat{X}^{u_2}_t-R_{1i}^{-1}(t)\big[B_{1i}(t)^\top\vp^{1i,u_2}_t+r_{1i}(t)\big],
\end{equation}
where $\vp^{1i,u_2}$ satisfies \eqref{e.3.14} and $\hat{X}^{u_2}$ satisfies
\begin{equation*}
\begin{aligned}
d\hat{X}^{u_2}_t&=\bigg[\bigg(A(t)-\sum_{i=1}^NB_{1i}(t)R_{1i}^{-1}(t)\big[B_{1i}(t)^\top P^{1i}_t+S_{1i}(t)\big]\bigg)\hat{X}^{u_2}_t\\
&\qquad -\sum_{i=1}^NB_{1i}(t)R_{1i}^{-1}(t)\big[B_{1i}(t)^\top\vp^{1i,u_2}_t+r_{1i}(t)\big]+B_2(t)u_2(t)+\al(t)\bigg]dt\\
&\quad +\big[\Sigma_tf^\top_1(t)+C_1(t)K_1^\top(t)\big](K_1(t)K_1^\top(t))^{-1}dV_t.
\end{aligned}
\end{equation*}
Now, for any fixed $u_2$, the optimal cost of the $i$th follower is given by
\begin{equation*}
\hspace{-2mm}\begin{aligned}
J_{1i}&(x_0;\bar{u}_1,u_2)=\EE\big[\langle P^{1i}_0x_0+2\vp^{1i,u_2}_0,x_0\rangle\big]+\EE\int_0^T\bigg\{\sum_{j=1}^{l_1}\big\langle\Pi^{1i}\Delta_j,\Delta_j\big\rangle +\sum_{j=1}^{l_2}\big\langle\Pi^{1i}C_{2j},C_{2j}\big\rangle\\
&-2\langle\Delta,\beta^1\rangle-\big\langle R_{1i}^{-1}(B_{1i}^\top\vp^{1i,u_2}+r_{1i}),B_{1i}^\top\vp^{1i,u_2}+r_{1i}\big\rangle +\sum_{j=1}^{l_1}\big\langle P^{1i}(\Delta_j+C_{1j}),\Delta_j+C_{1j}\big\rangle\\
& +2\langle B_2u_2,\vp^{1i,u_2}\rangle+2\langle\al,\vp^{1i,u_2}\rangle+2\langle C_2,\beta^2\rangle+2\big\langle(\Delta+C_1),\la^{1i,u_2}\big\rangle\bigg\}dt.
\end{aligned}
\end{equation*}
\begin{remark}
It is acknowledged that the separation principle introduced in \cite{Wh68} is the key step that separate the optimal control and state estimate in the classical partially observed stochastic control problem of SDEs. Then a backward separation principle is also introduced in \cite{WW08} for solving the partially observed FBSDEs. However, the necessary state decomposition technique in the above two work plays an important role, therefore, the restriction imposed on the admissible control requires additional condition in terms of the filtration generated by observation process. In their procedure, different from that in \cite{WW08} in which they show the ``backward" idea that first addresses the control problem and the filtering equation is then obtained, the classical idea in \cite{Wh68} is to investigate the filtering equation first, which is similar to the steps in our work. Based on the orthogonal decomposition, we focus on the decomposition of cost functional alternatively and then transform the original partially observed control into a class of ``completely observed" control problem based on the known information filtration, so we do not need the restriction introduced by state decomposition and relax the definition of admissible control. As far as we know, the similar technique is also applied to investigate the partially observed LQ stochastic control of mean-field SDEs recently in \cite{MB24}, but there is no any related result for forward backward stochastic system yet, which has been now completed in the following section for leader's problem.
\end{remark}

\section{The Leader's Problem}\label{leader pro}

In this section, the leader encounters a kind of partially observed stochastic optimal control of non-standard FBSDEs with conditional mean-field terms. So we extend the result in \cite{SX23} to the case of forward-backward partially observed stochastic system.

Since the state equation (\ref{state SDE}) of $X^{u_2}\equiv X^{\bar{u}_1,u_2}$, now, contains $\vp^{1i,u_2}$ and $\hat{X}^{u_2}\equiv\hat{X}^{\bar{u}_1,u_2}$ by $\bar{u}_{1i}$ in (\ref{bar u1}), thus the state equation of the leader becomes an FBSDE with conditional mean field terms. In the following part, for obtaining the explicit form, we consider the case that $f_1(\cdot)=g_1(\cdot)=0$, i.e., $\cF_t^{Y_1^{u_2}}\equiv\cF_t^{Y_1^{\bar{u}_1,u_2}}=\cF_t^{W^1}$, then $\check{V}_t\equiv\check{V}_t=W^1_t$. And thus
\begin{equation}\label{leader_state}
\left\{
\begin{aligned}
dX_t^{u_2}&=\bigg\{A(t)X_t^{u_2}-\sum_{i=1}^NB_{1i}(t)R_{1i}^{-1}(t)\big[B_{1i}^\top(t) P^{1i}_t+S_{1i}(t)\big]\hat{X}^{u_2}_t+B_2(t)u_2(t)+\al(t)\\
&\qquad -\sum_{i=1}^NB_{1i}(t)R_{1i}^{-1}(t)\big[B_{1i}^\top(t)\vp^{1i,u_2}_t+r_{1i}(t)\big] \bigg\}dt+C_1(t)dW^1_t+C_2(t)dW^2_t,\\
-d\vp^{1i,u_2}_t&=\Big\{\big[A(t)-B_{1i}(t)R_{1i}^{-1}(t)\big(B_{1i}^\top(t) P^{1i}_t+S_{1i}(t)\big)\big]\vp^{1i,u_2}_t+P^{1i}_t\al(t)+q_{1i}(t)\\
&\qquad -R_{1i}^{-1}(t)\big[B_{1i}^\top(t) P^{1i}_t+S_{1i}(t)\big]r_{1i}(t) +P^{1i}_tB_2(t)u_2(t)\Big\}dt-\la^{1i,u_2}_tdW^1_t,\\
X_0^{u_2}&=x_0,\ \vp^{1i,u_2}_T=g_{1i},\quad i=1,2,\dots,N.
\end{aligned}
\right.
\end{equation}
Now we have $\cF^{Y_2^{u_2}}_t\subset\cF_t^{Y_1^{u_2}}=\cF_t^{W^1}$, where we have denoted by $Y_2^{u_2}\equiv Y_2^{\bar{u}_1,u_2}$. Set, for a fixed but arbitrary admissible control $u_2\in\cU^L_{ad}$,
\begin{equation*}
\begin{aligned}
\check{X}_t&\triangleq\EE\Big[X_t\Big|\cF_t^{Y_2^{u_2}}\Big],\quad \tilde{\tilde{X}}_t=X_t-\check{X}_t,\quad \check{\vp}^{1i}_t\triangleq\EE\Big[\vp^{1i}_t\Big|\cF_t^{Y_2^{u_2}}\Big],\\
\check{\la}^{1i}_t&\triangleq\EE\Big[\la^{1i}_t\Big|\cF_t^{Y_2^{u_2}}\Big],\ i=1,2,\dots,N,\quad U_t\triangleq Y_2(t)-\int_0^t\big[f_2(s)\check{X}_s+g_2(s)\big]ds,
\end{aligned}
\end{equation*}
and $\check{U}_t\triangleq\int_0^tK_2^{-1}(s)dU_s$ is an $\cF^{Y_2^{u_2}}_t$-Brownian motion. First, we consider the process
\begin{equation*}
\begin{aligned}
\laa^{2}_t\triangleq\check{X}_t-x_0-\int_0^t\bigg\{&\bigg[A-\sum_{i=1}^NB_{1i}R_{1i}^{-1}\big(B_{1i}^\top P^{1i}+S_{1i}\big)\bigg]\check{X}-\sum_{i=1}^NB_{1i}R_{1i}^{-1}B_{1i}^\top\check{\vp}^{1i}\\
&\quad +B_2u_2-\sum_{i=1}^NB_{1i}R_{1i}^{-1}r_{1i}+\al\bigg\}ds,
\end{aligned}
\end{equation*}
which is $\cF_t^{Y^{u_2}_2}$-adapted and RCLL, $\laa^{2}_0=0$. For fixed $0\leq s<t\leq T$, we have
\begin{equation*}
\begin{aligned}
\EE\Big[\laa^{2}_t-\laa^{2}_s\Big|\cF_s^{Y^{u_2}_2}\Big]&=\EE\bigg[\check{X}_t-\check{X}_s-\int_s^t\Big[\Big(A-\sum_{i=1}^NB_{1i}R_{1i}^{-1}\big(B_{1i}^\top P^{1i}+S_{1i}\big)\Big)\check{X}\\
&\quad-
\sum_{i=1}^NB_{1i}R_{1i}^{-1}B_{1i}^\top\check{\vp}^{1i}+B_2u_2-\sum_{i=1}^NB_{1i}R_{1i}^{-1}r_{1i}+\al\Big]dr\bigg|\cF_s^{Y^{u_2}_2}\bigg],
\end{aligned}
\end{equation*}
where
\begin{equation*}
\begin{aligned}
&\EE\Big[\check{X}_t-\check{X}_s\Big|\cF_s^{Y^{u_2}_2}\Big]=\EE\Big[\EE\Big[X_t\Big|\cF_t^{Y^{u_2}_2}\Big]\Big|\cF_s^{Y^{u_2}_2}\Big]
-\EE\Big[X_s\Big|\cF_s^{Y^{u_2}_2}\Big]=\EE\Big[X_t-X_s\Big|\cF_s^{Y^{u_2}_2}\Big],\\
&\EE\bigg[\int_s^t\bigg(A-\sum_{i=1}^NB_{1i}R_{1i}^{-1}\big(B_{1i}^\top P^{1i}+S_{1i}\big)\bigg)\check{X}_rdr\bigg|\cF_s^{Y^{u_2}_2}\bigg]\\
&=\int_s^t\bigg(A-\sum_{i=1}^NB_{1i}R_{1i}^{-1}\big(B_{1i}^\top P^{1i}+S_{1i}\big)\bigg)\EE\Big[\check{X}_r\Big|\cF_s^{Y^{u_2}_2}\Big]dr\\
&=\int_s^tA\EE\Big[\EE\Big[\check{X}_r\Big|\cF_r^{Y^{u_2}_2}\Big]\Big|\cF_s^{Y^{u_2}_2}\Big]dr\\
&\quad -\int_s^t\sum_{i=1}^NB_{1i}R_{1i}^{-1}\big(B_{1i}^\top P^{1i}+S_{1i}\big)\EE\Big[\EE\Big[\check{X}_r\Big|\cF_r^{Y^{u_2}_2}\Big]\Big|\cF_s^{Y^{u_2}_2}\Big]dr\\
&=\int_s^tA\EE\Big[X_r\Big|\cF_s^{Y^{u_2}_2}\Big]dr\\
&\quad -\int_s^t\sum_{i=1}^NB_{1i}R_{1i}^{-1}\big(B_{1i}^\top P^{1i}+S_{1i}\big)\EE\Big[\EE\Big[\EE\Big[X_r\Big|\cF_r^{W^1}\Big]\Big|\cF_r^{Y^{u_2}_2}\Big]\Big|\cF_s^{Y^{u_2}_2}\Big]dr\\
&=\EE\bigg[\int_s^tAX_rdr\bigg|\cF_s^{Y^{u_2}_2}\bigg]
-\int_s^t\sum_{i=1}^NB_{1i}R_{1i}^{-1}\big(B_{1i}^\top P^{1i}+S_{1i}\big)\EE\Big[\EE\Big[\hat{X}_r\Big|\cF_r^{Y^{u_2}_2}\Big]\Big|\cF_s^{Y^{u_2}_2}\Big]dr\\
&=\EE\bigg[\int_s^tAX_rdr\bigg|\cF_s^{Y^{u_2}_2}\bigg]
-\int_s^t\sum_{i=1}^NB_{1i}R_{1i}^{-1}\big(B_{1i}^\top P^{1i}+S_{1i}\big)\EE\Big[\hat{X}_r\Big|\cF_s^{Y^{u_2}_2}\Big]dr\\
&=\EE\bigg[\int_s^tAX_r-\sum_{i=1}^NB_{1i}R_{1i}^{-1}\big(B_{1i}^\top P^{1i}+S_{1i}\big)\hat{X}_rdr\bigg|\cF_s^{Y^{u_2}_2}\bigg],
\end{aligned}
\end{equation*}
and
\begin{equation*}
\begin{aligned}
\EE\bigg[\int_s^tB_{1i}R_{1i}^{-1}B_{1i}^\top\check{\vp}_r^{1i}dr\bigg|\cF_s^{Y^{u_2}_2}\bigg]&=\int_s^tB_{1i}R_{1i}^{-1}B_{1i}^\top\EE\Big[\EE\Big[\vp_r^{1i}\Big|\cF_r^{Y^{u_2}_2}\Big]\Big|\cF_s^{Y^{u_2}_2}\Big]dr\\
&=\EE\bigg[\int_s^tB_{1i}R_{1i}^{-1}B_{1i}^\top\vp_r^{1i}dr\bigg|\cF_s^{Y^{u_2}_2}\bigg],
\end{aligned}
\end{equation*}
which implies that
\begin{equation*}
\begin{aligned}
\EE\Big[\laa^{2}_t-\laa^{2}_s\Big|\cF_s^{Y^{u_2}_2}\Big]&=\EE\bigg[X_t-X_s-\int_s^t\Big[AX-\sum_{i=1}^NB_{1i}R_{1i}^{-1}\big(B_{1i}^\top P^{1i}+S_{1i}\big)\hat{X}\\
&\qquad -\sum_{i=1}^NB_{1i}R_{1i}^{-1}B_{1i}^\top\vp^{1i}+B_2u_2-\sum_{i=1}^NB_{1i}R_{1i}^{-1}r_{1i}+\al\Big]dr\bigg|\cF_s^{Y^{u_2}_2}\bigg]\\
&=\EE\bigg[\int_s^tC_1dW^1+\int_s^tC_2dW^2\bigg|\cF_s^{Y^{u_2}_2}\bigg]\\
&=\EE\bigg[\EE\bigg[\int_s^tC_1dW^1+\int_s^tC_2dW^2\bigg|\cF_s\bigg]\bigg|\cF_s^{Y^{u_2}_2}\bigg]=0.
\end{aligned}
\end{equation*}
Thus, $\laa^{2}_t$ is an $\cF_t^{Y^{u_2}_2}$-martingale. Then similarly, there exists $\cF_s^{Y^{u_2}_2}$-progressively measurable $\la^{2}$, such that $\laa^{2}_t=\int_0^t\la^{2}_s(K_2(s)K_2^\top(s))^{-1}dU_s$. Let $\tilde{\zeta}=\{\tilde{\zeta}_t;0\leq t\leq T\}$ be a fixed but arbitrary $\RR^{n\times l_2}$-valued square-integrable $\cF_t^{Y^{u_2}_2}$-progressively measurable process. Consider the $\cF_t^{Y^{u_2}_2}$-martingale:
\begin{equation*}
\tilde{\vp}_t\triangleq\int_0^t\tilde{\zeta}_s(K_2(s)K_2^\top(s))^{-1}dU_s.
\end{equation*}
Applying It\^o's formula, we have
\begin{equation*}
\EE\big[\laa^{2}_t(\tilde{\vp}_t)^\top\big]=\EE\int_0^t\la^{2}_s(K_2(s)K_2^\top(s))^{-1}(\tilde{\zeta}_s)^\top ds,\ \forall\ 0\leq t\leq T.
\end{equation*}
Meanwhile, we can get
\begin{equation*}
\begin{aligned}
\laa^{2}_t(\tilde{\vp}_t)^\top&=\check{X}_t(\tilde{\vp}_t)^\top-x_0(\tilde{\vp}_0)^\top-\int_0^t\bigg\{\Big[A-\sum_{i=1}^NB_{1i}R_{1i}^{-1}\big(B_{1i}^\top P^{1i}+S_{1i}\big)\Big]\check{X}_s\\
&\quad -\sum_{i=1}^NB_{1i}R_{1i}^{-1}B_{1i}^\top\check{\vp}^{1i}_s+B_2u_2-\sum_{i=1}^NB_{1i}R_{1i}^{-1}r_{1i}+\al\bigg\}(\tilde{\vp}_s)^\top ds.
\end{aligned}
\end{equation*}
Taking expectation on both sides, and noting that
\begin{equation*}
\begin{aligned}
\EE\big[\check{X}_t(\tilde{\vp}_t)^\top\big]&=\EE\Big[\EE\Big[X_t\Big|\cF_t^{Y^{u_2}_2}\Big](\tilde{\vp}_t)^\top\Big]=\EE\big[X_t(\tilde{\vp}_t)^\top\big],\\
\EE\big[\check{X}_s(\tilde{\vp}_t)^\top\big]&=\EE\Big[\EE\Big[\check{X}_s(\tilde{\vp}_t)^\top\Big|\cF_s^{Y^{u_2}_2}\Big]\Big]
=\EE\Big[\check{X}_s\EE\Big[(\tilde{\vp}_t)^\top\Big|\cF_s^{Y^{u_2}_2}\Big]\Big]\\
&=\EE\big[\check{X}_s(\tilde{\vp}_s)^\top\big]=\EE\big[X_s(\tilde{\vp}_s)^\top\big],\\
\end{aligned}
\end{equation*}
\begin{equation*}
\begin{aligned}
\EE\big[\check{X}_s(\tilde{\vp}_t)^\top\big]&=\EE\Big[\EE\Big[X_s\Big|\cF_s^{Y^{u_2}_2}\Big](\tilde{\vp}_t)^\top\Big]
=\EE\Big[\EE\Big[\EE\Big[X_s\Big|\cF_s^{W^1}\Big]\Big|\cF_s^{Y^{u_2}_2}\Big](\tilde{\vp}_t)^\top\Big]\\
&=\EE\Big[\EE\Big[\hat{X}_s\Big|\cF_s^{Y^{u_2}_2}\Big](\tilde{\vp}_t)^\top\Big]=\EE\Big[\EE\Big[\EE\Big[\hat{X}_s\Big|\cF_s^{Y^{u_2}_2}\Big](\tilde{\vp}_t)\Big|\cF_s^{Y^{u_2}_2}\Big]\\
&=\EE\Big[\EE\Big[\hat{X}_s\Big|\cF_s^{Y^{u_2}_2}\Big](\tilde{\vp}_s)^\top\Big]=\EE\big[\hat{X}_s(\tilde{\vp}_s)^\top\big],\\
\EE\big[\check{\vp}^{1i}_s(\tilde{\vp}_t)^\top\big]&=\EE\Big[\EE\Big[\check{\vp}^{1i}_s(\tilde{\vp}_t)^\top\Big|\cF_s^{Y^{u_2}_2}]\Big]=\EE\Big[\check{\vp}^{1i}_s(\tilde{\vp}_s)^\top\Big]\\
&=\EE\Big[\EE\Big[\vp^{1i}_s\Big|\cF_s^{Y^{u_2}_2}\Big](\tilde{\vp}_s)^\top\Big]=\EE\big[\vp^{1i}_s(\tilde{\vp}_s)^\top\big],\\
\EE\big[u_2(s)(\tilde{\vp}_t)^\top\big]&=\EE\Big[\EE\Big[u_2(s)(\tilde{\vp}_t)^\top\Big|\cF_s^{Y^{u_2}_2}\Big]\Big]=\EE\big[u_2(s)(\tilde{\vp}_s)^\top\big],
\end{aligned}
\end{equation*}
we get
\begin{equation*}
\begin{aligned}
\EE\big[\laa^{2}_t(\tilde{\vp}_t)^\top\big]&=\EE\big[X_t(\tilde{\vp}_t)^\top\big]-\EE\bigg[\int_0^t\bigg\{AX_s
-\sum_{i=1}^NB_{1i}R_{1i}^{-1}\big(B_{1i}^\top P^{1i}+S_{1i}\big)\hat{X}_s\\
&\quad -\sum_{i=1}^NB_{1i}R_{1i}^{-1}B_{1i}^\top\vp^{1i}_s+B_2u_2-\sum_{i=1}^NB_{1i}R_{1i}^{-1}r_{1i}+\al\Big](\tilde{\vp}_s)^\top ds\bigg].
\end{aligned}
\end{equation*}
Noting
\begin{equation*}
\begin{aligned}
d\tilde{\vp}_t&=\tilde{\zeta}_t(K_2(t)K_2^\top(t))^{-1}dU_t=\tilde{\zeta}_t(K_2(t)K_2^\top(t))^{-1}\big[f_2(t)(X_t-\check{X}_t)dt+K_2(t)dW^2_t\big]\\
&=\tilde{\zeta}_t(K_2(t)K_2^\top(t))^{-1}f_2(t)\tilde{\tilde{X}}_tdt+\tilde{\zeta}_t(K_2(t)K_2^\top(t))^{-1}K_2(t)dW^2_t,
\end{aligned}
\end{equation*}
and applying It\^o's formula, it follows that
\begin{equation*}
\begin{aligned}
\EE\big[X_t(\tilde{\vp}_t)^\top\big]&=\EE\int_0^t\bigg\{X_s(\tilde{\tilde{X}}_s)^\top f^\top_2(K_2K_2^\top)^{-1}(\tilde{\zeta}_s)^\top\\
&\qquad\quad +\bigg[AX_s-\sum_{i=1}^NB_{1i}R_{1i}^{-1}\big(B_{1i}^\top P^{1i}+S_{1i}\big)\hat{X}_s-\sum_{i=1}^NB_{1i}R_{1i}^{-1}B_{1i}^\top\vp^{1i}_s\\
&\qquad\quad +B_2u_2-\sum_{i=1}^NB_{1i}R_{1i}^{-1}r_{1i}+\al\bigg](\tilde{\vp}_s)^\top+C_2K_2^\top(K_2K_2^\top)^{-1}(\tilde{\zeta}_s)^\top\bigg\}ds,
\end{aligned}
\end{equation*}
among which
\begin{equation*}
\begin{aligned}
&\EE\big[X_s(\tilde{\tilde{X}}_s)^\top f^\top_2(K_2K_2^\top)^{-1}(\tilde{\zeta}_s)^\top\big]\\
&=\EE\Big[\EE\Big[X_s(\tilde{\tilde{X}}_s)^\top f^\top_2(K_2K_2^\top)^{-1}(\tilde{\zeta}_s)^\top\Big|\cF_s^{Y^{u_2}_2}\Big]\Big]\\
&=\EE\Big[\EE\Big[(\check{X}_s+\tilde{\tilde{X}}_s)(\tilde{\tilde{X}}_s)^\top\Big|\cF_s^{Y^{u_2}_2}\Big]f^\top_2(K_2K_2^\top)^{-1}(\tilde{\zeta}_s)^\top\Big]\\
&=\EE\Big[\EE\Big[\check{X}_s(\tilde{\tilde{X}}_s)^\top+\tilde{\tilde{X}}_s(\tilde{\tilde{X}}_s)^\top\Big|\cF_s^{Y^{u_2}_2}\Big]f^\top_2(K_2K_2^\top)^{-1}(\tilde{\zeta}_s)^\top\Big]\\
&=\EE\Big[\EE\Big[\check{X}_s(\tilde{\tilde{X}}_s)^\top\Big|\cF_s^{Y^{u_2}_2}\Big]f^\top_2(K_2K_2^\top)^{-1}(\tilde{\zeta}_s)^\top\Big]
+\EE\big[\tilde{\Sigma}_sf^\top_2(K_2K_2^\top)^{-1}(\tilde{\zeta}_s)^\top\big]\\
&=\EE\big[\tilde{\Sigma}_sf^\top_2(K_2K_2^\top)^{-1}(\tilde{\zeta}_s)^\top\big],
\end{aligned}
\end{equation*}
where, similarly to the discussion in the counterpart of followers' problem, $\tilde{\Sigma}_s\triangleq\EE\big[\tilde{\tilde{X}}_s(\tilde{\tilde{X}}_s)^\top\big]$, and $\tilde{\tilde{X}}_t$ is independent of $u_2$, then we also apply the fact that for each fixed admissible control $u_2$, $\tilde{\tilde{X}}_t$ is independent of $\cF^{Y^{u_2}_2}_t$. Thus, we have
\begin{equation*}
\begin{aligned}
\EE&\big[X_t(\tilde{\vp}_t)^\top\big]=\EE\int_0^t\bigg\{\bigg[AX_s-\sum_{i=1}^NB_{1i}R_{1i}^{-1}\big(B_{1i}^\top P^{1i}+S_{1i}\big)\hat{X}_s-\sum_{i=1}^NB_{1i}R_{1i}^{-1}B_{1i}^\top\vp^{1i}_s\\
& +B_2u_2-\sum_{i=1}^NB_{1i}R_{1i}^{-1}r_{1i}+\al\bigg](\tilde{\vp}_s)^\top +\big(\tilde{\Sigma}_sf_2^\top+C_2K_2^\top\big)(K_2K_2^\top)^{-1}(\tilde{\zeta}_s)^\top\bigg\}ds,
\end{aligned}
\end{equation*}
and
\begin{equation*}
\EE\big[\laa^{2}_t(\tilde{\vp}_t)^\top\big]=\EE\int_0^t\big(\tilde{\Sigma}_sf_2^\top+C_2K_2^\top\big)(K_2K_2^\top)^{-1}(\tilde{\zeta}_s)^\top ds,
\end{equation*}
Therefore, we obtain $\la^{2}_t=\tilde{\Sigma}_tf_2^\top(t)+C_2(t)K_2^\top(t),0\leq t\leq T$, and the equation of $\check{X}$ is
\begin{equation*}
\begin{aligned}
\check{X}_t&=x_0+\int_0^t\bigg\{\Big[A-\sum_{i=1}^NB_{1i}R_{1i}^{-1}\big(B_{1i}^\top P^{1i}+S_{1i}\big)\Big]\check{X}-\sum_{i=1}^NB_{1i}R_{1i}^{-1}B_{1i}^\top\check{\vp}^{1i}+B_2u_2\\
&\quad-\sum_{i=1}^NB_{1i}R_{1i}^{-1}r_{1i}+\al\bigg\}ds+\int_0^t\big(\tilde{\Sigma}_sf_2^\top+C_2K_2^\top\big)(K_2K_2^\top)^{-1}dU_s.
\end{aligned}
\end{equation*}

In the following, we derive, for a fixed admissible control $u_2\in\cU^L_{ad}$, the filtering equation for $\check{\vp}^{1i}$. We introduce
\begin{equation*}
\begin{aligned}
\Ga_t^{i}&\triangleq\check{\vp}^{1i}_t-g_{1i}-\int_t^T\Big\{\big[A(s)-B_{1i}(s)R_{1i}^{-1}(s)\big(B_{1i}(s)^\top P^{1i}_s+S_{1i}(s)\big)\big]\check{\vp}^{1i}_s+P^{1i}_s\al(s)\\
&\quad +q_{1i}(s)-R_{1i}^{-1}(s)\big(B_{1i}(s)^\top P^{1i}_s+S_{1i}(s)\big)r_{1i}(s)+P^{1i}_sB_2(s)u_2(s)\Big\}ds,\quad i=1,2,\dots,N.
\end{aligned}
\end{equation*}
For fixed $0\leq s<t\leq T$, we have
\begin{equation*}
\begin{aligned}
&\EE\Big[\Ga_t^{i}-\Ga_s^{i}\Big|\cF_s^{Y^{u_2}_2}\Big]=\EE\bigg[\check{\vp}^{1i}_t-\check{\vp}^{1i}_s+\int_s^t\Big\{\big[A-B_{1i}R_{1i}^{-1}\big(B_{1i}^\top P^{1i}+S_{1i}\big)\big]\check{\vp}^{1i}\\
&\qquad +P^{1i}\al+q_{1i}-R_{1i}^{-1}(B_{1i}^\top P^{1i}+S_{1i})r_{1i}+P^{1i}B_2u_2\Big\}dr\bigg|\cF_s^{Y^{u_2}_2}\bigg],
\end{aligned}
\end{equation*}
where
\begin{equation*}
\begin{aligned}
&\EE\Big[\check{\vp}^{1i}_t-\check{\vp}^{1i}_s\Big|\cF^{Y^{u_2}_2}_s\Big]\\
&=\EE\Big[\EE\Big[\vp^{1i}_t\Big|\cF_t^{Y^{u_2}_2}\Big]\Big|\cF_s^{Y^{u_2}_2}\Big]-\EE\Big[\vp_s^{1i}\Big|\cF_s^{Y^{u_2}_2}\Big]=\EE\Big[\vp^{1i}_t-\vp_s^{1i}\Big|\cF_s^{Y^{u_2}_2}\Big],
\end{aligned}
\end{equation*}
and
\begin{equation*}
\begin{aligned}
&\EE\bigg[\int_s^t\big[A-B_{1i}R_{1i}^{-1}\big(B_{1i}^\top P^{1i}+S_{1i}\big)\big]\check{\vp}^{1i}_rdr\bigg|\cF_s^{Y^{u_2}_2}\bigg]\\
&=\int_s^t\big[A-B_{1i}R_{1i}^{-1}\big(B_{1i}^\top P^{1i}+S_{1i}\big)\big]\EE\Big[\EE\Big[\vp^{1i}_r\Big|\cF^{Y^{u_2}_2}_r\Big]\Big|\cF_s^{Y^{u_2}_2}\Big]dr\\
&=\EE\bigg[\int_s^t\big[A-B_{1i}R_{1i}^{-1}\big(B_{1i}^\top P^{1i}+S_{1i}\big)\big]\vp^{1i}_rdr\bigg|\cF_s^{Y^{u_2}_2}\bigg],
\end{aligned}
\end{equation*}
which imply that
\begin{equation*}
\begin{aligned}
&\EE\Big[\Ga_t^{i}-\Ga_s^{i}\Big|\cF_s^{Y^{u_2}_2}\Big]=\EE\bigg[\vp^{1i}_t-\vp^{1i}_s+\int_s^t\Big\{\big[A-B_{1i}R_{1i}^{-1}\big(B_{1i}^\top P^{1i}+S_{1i}\big)\big]\vp^{1i}\\
&\qquad +P^{1i}\al+q_{1i}-R_{1i}^{-1}\big(B_{1i}^\top P^{1i}+S_{1i}\big)r_{1i}+P^{1i}B_2u_2\Big\}dr\bigg|\cF_s^{Y^{u_2}_2}\bigg]\\
&=\EE\bigg[\int_s^t\la_r^{1}dW^1\bigg|\cF_s^{Y^{u_2}_2}\bigg]=\EE\bigg[\EE\bigg[\int_s^t\la_r^{1}dW^1\bigg|\cF_s^{W^1}\bigg]\bigg|\cF_s^{Y^{u_2}_2}\bigg]=0,
\end{aligned}
\end{equation*}
so $\Ga_t^{i}$ is an $\cF^{Y^{u_2}_2}_t$-martingale. Similarly, we have
\begin{equation*}
\Ga_t^{i}=\int_0^t\tilde{\eta}_s^{i}(K_2K_2^\top)^{-1}dU_s,\ \forall\ 0\leq t\leq T,
\end{equation*}
where $\tilde{\eta}^{i}=\{\tilde{\eta}^{i}_t;0\leq t\leq T\}$ is to be determined, $i=1,2,\dots,N$. Let $\tilde{\beta}^{i}=\{\tilde{\beta}^{i}_t;0\leq t\leq T\}$ be a fixed but arbitrary $\RR^{n\times l_2}$-valued square-integrable $\cF_t^{Y^{u_2}_2}$-progressively process. Consider the $\{\cF_t^{Y^{u_2}_2}\}$-martingale
\begin{equation*}
\varepsilon^{i}_t\triangleq\int_0^t\tilde{\beta}^{i2}_s(K_2(s)K_2^\top(s))^{-1}dU_s,
\end{equation*}
and apply It\^o's formula, then we get
\begin{equation*}
\EE\big[\Ga_t^{i}(\varepsilon_t^{i})^\top\big]=\EE\int_0^t\tilde{\eta}^{i}_s(K_2(s)K_2^\top(s))^{-1}(\tilde{\beta}^{i}_s)^\top ds.
\end{equation*}
Moreover, we have
\begin{equation*}
\begin{aligned}
&\EE\big[\Ga_t^{i}(\varepsilon_T^{i})^\top\big]\\
&=\EE\big[\check{\vp}^{1i}_t(\varepsilon^{i}_T)^\top\big]-\EE\int_t^T\Big\{\big[A(s)-B_{1i}(s)R_{1i}^{-1}(s)\big(B_{1i}(s)^\top P^{1i}_s+S_{1i}(s)\big)\big]\check{\vp}^{1i}_s\\
&\quad +P^{1i}_s\al(s)+q_{1i}(s)-R_{1i}^{-1}(s)\big(B_{1i}(s)^\top P^{1i}_s+S_{1i}(s)\big)r_{1i}(s)+P^{1i}_sB_2(s)u_2(s)\Big\}(\varepsilon^{i}_T)^\top ds,
\end{aligned}
\end{equation*}
where
\begin{equation*}
\begin{aligned}
\EE\big[\Ga_t^{i}(\varepsilon_T^{i})^\top\big]&=\EE\Big[\EE\Big[\Ga_t^{i}(\varepsilon_T^{i})^\top\Big|\cF_t^{Y^{u_2}_2}\Big]\Big]=\EE\big[\Ga_t^{i}(\varepsilon_t^{i})^\top\big],\\
\EE\big[\check{\vp}^{1i}_t(\varepsilon^{i}_T)^\top\big]&=\EE\Big[\EE\Big[\check{\vp}^{1i}_t(\varepsilon^{i}_T)^\top\Big|\cF_t^{Y^{u_2}_2}\Big]\Big]=\EE\big[\check{\vp}^{1i}_t(\varepsilon^{i}_t)^\top\big]\\
&=\EE\Big[\EE\Big[\vp^{1i}_t\Big|\cF_t^{Y^{u_2}_2}\Big](\varepsilon^{i}_t)^\top\Big]=\EE\big[\vp^{1i}_t(\varepsilon^{i}_t)^\top\big],\\
\EE\big[\check{\vp}^{1i}_s(\varepsilon^{i}_T)^\top\big]&=\EE\Big[\EE\Big[\check{\vp}^{1i}_s(\varepsilon^{i}_T)^\top\Big|\cF_t^{Y^{u_2}_2}\Big]\Big]=\EE\big[\check{\vp}^{1i}_s(\varepsilon^{i}_s)^\top\big]\\
&=\EE\Big[\EE\Big[\vp^{1i}_s\Big|\cF_t^{Y^{u_2}_2}\Big](\varepsilon^{i}_s)^\top\Big]=\EE\big[\vp^{1i}_s(\varepsilon^{i}_s)^\top\big],\\
\EE\big[u_2(s)(\varepsilon^{i}_T)^\top\big]&=\EE\Big[\EE\Big[u_2(s)(\varepsilon^{i}_T)^\top\Big|\cF_s^{Y^{u_2}_2}\Big]\Big]=\EE\big[u_2(s)(\varepsilon^{i}_s)^\top\big],
\end{aligned}
\end{equation*}
which implies that
\begin{equation*}
\begin{aligned}
&\EE\big[\Ga_t^{i}(\varepsilon_t^{i})^\top\big]=\EE\big[\vp^{1i}_t(\varepsilon^{i}_t)^\top\big]-\EE\int_t^T\Big\{\big[A(s)-B_{1i}(s)R_{1i}^{-1}(s)\big(B_{1i}^\top(s)P^{1i}_s+S_{1i}(s)\big)\big]\vp^{1i}_s\\
&\qquad +P^{1i}_s\al(s)+q_{1i}(s)-R_{1i}^{-1}(s)\big(B_{1i}^\top(s)P^{1i}_s+S_{1i}(s)\big)r_{1i}(s)+P^{1i}_sB_2(s)u_2(s)\Big\}(\varepsilon^{i}_s)^\top ds.
\end{aligned}
\end{equation*}
Now, we also get
\begin{equation*}
d\varepsilon^{i}_t=\tilde{\beta}^{i}_t(K_2K_2^\top)^{-1}dU_t=\tilde{\beta}^{i}_t(K_2K_2^\top)^{-1}f_2\tilde{X}^{u_2}_tdt+\tilde{\beta}^{i}_t(K_2K_2^\top)^{-1}K_2dW^2.
\end{equation*}
By It\^o's formula, we have
\begin{equation*}
\begin{aligned}
\EE\big[g_{1i}(\varepsilon^{i}_T)^\top\big]&=\EE\big[\vp^{1i}_t(\varepsilon^{i}_t)^\top\big]+\EE\int_t^T\Big\{\vp^{1i}_s(\tilde{X}_s)^\top f^\top_2(K_2K_2^\top)^{-1}(\tilde{\beta}^{i}_s)^\top\\
&\qquad -\big[\big[A-B_{1i}R_{1i}^{-1}\big(B_{1i}^\top P^{1i}+S_{1i}\big)\big]\vp^{1i}_s+P^{1i}\al+q_{1i}\\
&\qquad -R^{-1}_{1i}(B_{1i}^\top P^{1i}+S_{1i})r_{1i}+P^{1i}B_2u_2\big](\varepsilon^{i}_s)^\top\Big\} ds,
\end{aligned}
\end{equation*}
among which
\begin{equation*}
\begin{aligned}
&\EE\big[\vp^{1i}_s(\tilde{\tilde{X}}_s)^\top f^\top_2(K_2K_2^\top)^{-1}(\tilde{\beta}^i_s)^\top\big]=\EE\Big[\EE\Big[\vp^{1i}_s(\tilde{\tilde{X}}_s)^\top f^\top_2(K_2K_2^\top)^{-1}(\tilde{\beta}^i_s)^\top\Big|\cF_s^{Y^{u_2}_2}\Big]\Big]\\
&\qquad=\EE\Big[\EE\Big[\check{\vp}^{1i}_s(\tilde{\tilde{X}}_s)^\top\Big|\cF_s^{Y^{u_2}_2}\Big]f^\top_2(K_2K_2^\top)^{-1}(\tilde{\beta}^i_s)^\top\Big]\\
&\qquad\quad+\EE\Big[\EE\Big[\big(\vp^{1i}_s-\check{\vp}^{1i}_s\big)(\tilde{\tilde{X}}_s)^\top\Big|\cF_s^{Y^{u_2}_2}\Big]f^\top_2(K_2K_2^\top)^{-1}(\tilde{\beta}^i_s)^\top\Big]\\
&\qquad=\EE\big[\check{\Sigma}^{i}_sf^\top_2(K_2K_2^\top)^{-1}(\tilde{\beta}^i_s)^\top\big],
\end{aligned}
\end{equation*}
where $\check{\Sigma}^{i}_s\triangleq\EE\big[\big(\vp^{1i}_s-\check{\vp}^{1i}_s\big)(\tilde{\tilde{X}}_s)^\top\big]$, similarly, $\vp^{1i}_t-\check{\vp}^{1i}_t$ is independent of $u_2$, and for each fixed admissible control $u_2$, $\vp^{1i}_t-\check{\vp}^{1i}_t$ is independent of $\cF^{Y^{u_2}_2}_t$. Then we have
\begin{equation*}
\EE\big[\Ga_t^{i}(\varepsilon_t^{i})^\top\big]=-\EE\int_t^T\check{\Sigma}^{i}_sf^\top_2(K_2K_2^\top)^{-1}(\tilde{\beta}^i_s)^\top ds.
\end{equation*}
Comparing this with $\EE\big[\Ga_t^{i}(\varepsilon_t^{i})^\top\big]=\EE\int_0^t\tilde{\eta}^{i}_s(K_2K_2^\top)^{-1}(\tilde{\beta}^i_s)^\top ds$, we achieve
\begin{equation*}
0=\EE\big[\Ga_T^{i}(\varepsilon_T^{i})^\top\big]=\EE\int_0^T\tilde{\eta}^{i}_s(K_2K_2^\top)^{-1}(\tilde{\beta}^i_s)^\top ds,
\end{equation*}
and
\begin{equation*}
\EE\big[\Ga_t^{i}(\varepsilon_t^{i})^\top\big]=-\EE\int_t^T\tilde{\eta}^{i}_s(K_2K_2^\top)^{-1}(\tilde{\beta}^i_s)^\top ds,
\end{equation*}
which implies that $\tilde{\eta}^{i}_t=\check{\Sigma}^{i}_tf^\top_2(t),0\leq t\leq T$. Therefore, we obtain the equation of $\check{\vp}^{1i}$:
\begin{equation*}
\left\{
\begin{aligned}
-d\check{\vp}^{1i}_t&=\Big\{\big[A-B_{1i}R_{1i}^{-1}\big(B_{1i}^\top P^{1i}+S_{1i}\big)\big]\check{\vp}^{1i}+P^{1i}\al+q_{1i}\\
&\qquad -R_{1i}^{-1}\big(B_{1i}^\top P^{1i}+S_{1i}\big)r_{1i}+P^{1i}B_2u_2\Big\}dt-\check{\ga}^{1i}_td\check{U}_t,\\
\check{\vp}^{1i}_T&=g_{1i},
\end{aligned}
\right.
\end{equation*}
where $\check{\ga}^{1i}_t\triangleq\check{\Sigma}^{i}_tf^\top_2(t)(K_2(t)K_2^\top(t))^{-1}$.

Next, we mainly give the coupled representations of $\tilde{\Sigma}\equiv\EE\big[(X-\check{X})(X-\check{X})^\top\big]$ and $\check{\Sigma}^{i}\equiv\EE\big[\big(\vp^{1i}-\check{\vp}^{1i}\big)(X-\check{X})^\top\big]$ by the following equations:
\begin{equation}\label{tilde Sigma}
\begin{aligned}
\tilde{\Sigma}_t&=\int_0^t\bigg\{\big[A-\big(\tilde{\Sigma}f_2^\top+C_2K_2^\top\big)(K_2K_2^\top)^{-1}f_2\big]\tilde{\Sigma}+\tilde{\Sigma}\big[A-\big(\tilde{\Sigma}f_2^\top+C_2K_2^\top\big)(K_2K_2^\top)^{-1}f_2\big]^\top\\
&\qquad\quad -\sum_{i=1}^NB_{1i}R_{1i}^{-1}\big(B_{1i}^\top P^{1i}+S_{1i}\big)\Sigma_{11}-\Sigma_{11}\bigg[\sum_{i=1}^NB_{1i}R_{1i}^{-1}\big(B_{1i}^\top P^{1i}+S_{1i}\big)\bigg]^\top
\\
&\qquad\quad -\sum_{i=1}^NB_{1i}R_{1i}^{-1}B_{1i}^\top\check{\Sigma}^{i}-\sum_{i=1}^N\check{\Sigma}^{i}\big(B_{1i}R_{1i}^{-1}B_{1i}^\top\big)^\top+C_1C_1^\top
+\tilde{\Sigma}f^\top_2(K_2K_2^\top)^{-1}f_2\tilde{\Sigma}\bigg\}ds,
\end{aligned}
\end{equation}
and
\begin{equation}\label{check Sigma}
\begin{aligned}
\check{\Sigma}^{i}_t&=\int_0^t\bigg\{\Big[A-\big(\tilde{\Sigma}f_2^\top+C_2K_2^\top\big)(K_2K_2^\top)^{-1}f_2-\Big(A-B_{1i}R_{1i}^{-1}\big(B_{1i}^\top P^{1i}+S_{1i}\big)\Big)\Big]\check{\Sigma}^i\\
&\qquad\quad -B_{1i}R_{1i}^{-1}\big(B_{1i}^\top P^{1i}+S_{1i}\big)\Sigma^i_{13}-B_{1i}R_{1i}^{-1}B_{1i}^\top\Sigma^i_{14}\\
&\qquad\quad -\check{\Sigma}^{i}f_2^\top(K_2K_2^\top)^{-1}f_2\tilde{\Sigma}+\la^{1i}C^\top_1+\tilde{\Sigma}f^\top_2K_2^{-1}(\check{\ga}^{1i})^\top\bigg\}ds,
\end{aligned}
\end{equation}
where $\Sigma_{11}:=\EE[(\hat{X}-\check{X})(\tilde{X})^\top],\Sigma_{12}:=\EE[(\hat{X}-\check{X})(\hat{X}-\check{X})^\top],\Sigma^i_{13}:=\EE[(\hat{X}-\check{X})(\vp^{1i}-\check{\vp}^{1i})^\top],\\\Sigma^i_{14}:=\EE[(\vp^{1i}-\check{\vp}^{1i})(\vp^{1i}-\check{\vp}^{1i})^\top]$ satisfy
\begin{equation}\label{hat check tilde X u2}
\begin{aligned}
&\Sigma_{11}=\int_0^t\bigg\{\Big[A-\big(\tilde{\Sigma}f_2^\top+C_2K_2^\top\big)(K_2K_2^\top)^{-1}f_2\Big]\tilde{\Sigma}-\sum_{i=1}^NB_{1i}R_{1i}^{-1}\big(B_{1i}^\top P^{1i}+S_{1i}\big)\Sigma_{12}\\
&\qquad\quad -\sum_{i=1}^NB_{1i}R_{1i}^{-1}B_{1i}^\top(\Sigma^i_{13}+\check{\Sigma}^{i})+\Big[A-\sum_{i=1}^NB_{1i}R_{1i}^{-1}\big(B_{1i}^\top P^{1i}+S_{1i}\big)\Big]\tilde{\Sigma}+C_1C_1^\top\bigg\}ds,
\end{aligned}
\end{equation}
\begin{equation}\label{hat check X u2}
\hspace{-15mm}\begin{aligned}
&\Sigma_{12}=\int_0^t\bigg\{\Big[A-\sum_{i=1}^NB_{1i}R_{1i}^{-1}\big(B_{1i}^\top P^{1i}+S_{1i}\big)\Big]\Sigma_{12}+\Sigma_{12}\Big[A-\sum_{i=1}^NB_{1i}R_{1i}^{-1}\big(B_{1i}^\top P^{1i}+S_{1i}\big)\Big]^\top\\
&\qquad\quad -\sum_{i=1}^NB_{1i}R_{1i}^{-1}B_{1i}^\top\Sigma^i_{13}-\sum_{i=1}^N\Sigma^i_{13}\big(B_{1i}R_{1i}^{-1}B_{1i}^\top\big)^\top+\big(\tilde{\Sigma}f^\top_2+C_2K_2^\top\big)(K_2K_2^\top)^{-1}\big(\tilde{\Sigma}f^\top_2+C_2K_2^\top\big)^\top\\
&\qquad\quad -\big(\tilde{\Sigma}f^\top_2+C_2K_2^\top\big)(K_2K_2^\top)^{-1}f_2\Sigma_{11}-\Sigma_{11}[(\tilde{\Sigma}f^\top_2+C_2K_2^\top)(K_2K_2^\top)^{-1}f_2]^\top+C_1C_1^\top\bigg\}ds,
\end{aligned}
\end{equation}
\begin{equation}\label{hat check X u2 and varphi}
\hspace{-20mm}\begin{aligned}
&\Sigma^i_{13}=\int_0^t\Big\{-\tilde{\Sigma}f_2^\top(K_2K_2^\top)^{-1}f_2\Sigma_{11}-B_{1i}R_{1i}^{-1}B_{1i}^\top\Sigma^i_{14}\\
&\qquad\quad -\big(\tilde{\Sigma}f^\top_2+C_2K_2^\top\big)(K_2K_2^\top)^{-1}f_2[\check{\Sigma}^{i}-(\check{\ga}^{1i})^\top]+\la^{1i,u_2}C_1^\top \Big\}ds,
\end{aligned}
\end{equation}
and
\begin{equation}\label{check varphi}
\begin{aligned}
&\Sigma^i_{14}=\int_t^T\Big\{\big[A-B_{1i}R_{1i}^{-1}\big(B_{1i}^\top P^{1i}+S_{1i}\big)\big]\Sigma^i_{14}+\Sigma^i_{14}\big[A-B_{1i}R_{1i}^{-1}\big(B_{1i}^\top P^{1i}+S_{1i}\big)\big]^\top\\
&\qquad\quad +2\check{\Sigma}^{i}f^\top_2(K_2K_2^\top)^{-1}f_2\check{\Sigma}^{i}-\la^{1i,u_2}(\la^{1i,u_2})^\top-\check{\ga}^{1i}(\check{\ga}^{1i})^\top\Big\}ds.
\end{aligned}
\end{equation}

Now, we can regard the partially observed problem of the leader in this section, as the completely observable case with the following state equation:
\begin{equation}\label{leader new state}
\left\{
\begin{aligned}
d\check{X}_t&=\bigg\{\Big[A-\sum_{i=1}^NB_{1i}R_{1i}^{-1}\big(B_{1i}^\top P^{1i}+S_{1i}\big)\Big]\check{X}-\sum_{i=1}^NB_{1i}R_{1i}^{-1}B_{1i}^\top\check{\vp}^{1i}+B_2u_2\\
&\qquad -\sum_{i=1}^NB_{1i}R_{1i}^{-1}r_{1i}+\al+\big(\tilde{\Sigma}f^\top_2+C_2K_2^\top\big)(K_2K_2^\top)^{-1}f_2\tilde{X}\bigg\}dt\\
&\quad +\Big(\tilde{\Sigma}f^\top_2(K_2^\top)^{-1}+C_2\Big)dW^2_t,\\
-d\check{\vp}^{1i}_t&=\Big\{\Big[A-B_{1i}R_{1i}^{-1}\big(B_{1i}^\top P^{1i}+S_{1i}\big)\Big]\check{\vp}^{1i}+P^{1i}\al+q_{1i}\\
&\qquad -R_{1i}^{-1}\big(B_{1i}^\top P^{1i}+S_{1i}\big)r_{1i}+P^{1i}B_2u_2\Big\}dt-\check{\ga}^{1i}_td\check{U}_t,\\
\check{X}_0&=x_0,\ \check{\vp}^{1i}_T=g_{1i},\quad i=1,2,\dots,N.
\end{aligned}
\right.
\end{equation}

In the following, we will consider the following reduced case of the leader's cost functional with $g_2=S_2=q_2=r_2=0$ (taking non-zero values has no essential influence besides the complex computation) and also by same orthogonal decomposition we have:
\begin{equation}\label{leader new cf}
\begin{aligned}
\tilde{J}_2(x_0;u_2)&\equiv J_2(x_0;\bar{u}_1,u_2)=\EE\bigg[\langle G_2X_T,X_T\rangle+\int_0^T\Big(\langle Q_2X,X\rangle+\langle R_2u_2,u_2\rangle\Big)dt\bigg]\\
&=\EE\bigg[\langle G_2\check{X}_T,\check{X}_T\rangle+\int_0^T\Big(\langle Q_2\check{X},\check{X}\rangle+\langle R_2u_2,u_2\rangle\Big)dt\bigg]\\
&\quad+\EE\bigg[\langle G_2\tilde{\tilde{X}}_T,\tilde{\tilde{X}}_T\rangle+\int_0^T\langle Q_2\tilde{\tilde{X}},\tilde{\tilde{X}}\rangle dt\bigg]\triangleq\check{J}_2(x_0;u_2)+\tilde{\tilde{J}}_2(x_0).
\end{aligned}
\end{equation}
\begin{remark}
We claim that $\tilde{\tilde{J}}_2(x_0)$ is independent of the control $u_2$ after orthogonal decomposition. Indeed, combining \eqref{leader_state} with \eqref{leader new state}, the filtering error processes $\tilde{\tilde{X}}$ and $\tilde{\tilde{\vp}}^{1i}\triangleq\vp^{1i}-\check{\vp}^{1i}$ satisfy
\begin{equation}\label{leader new error}
\left\{
\begin{aligned}
d\tilde{\tilde{X}}_t&=\bigg\{A\tilde{\tilde{X}}-\sum_{i=1}^NB_{1i}R_{1i}^{-1}\big(B_{1i}^\top P^{1i}+S_{1i}\big)(\hat{X}-\check{X})-\sum_{i=1}^NB_{1i}R_{1i}^{-1}B_{1i}^\top\tilde{\tilde{\vp}}^{1i}\\
&\qquad -\big(\tilde{\Sigma}f^\top_2+C_2K_2^\top\big)(K_2K_2^\top)^{-1}f_2\tilde{\tilde{X}}\bigg\}dt+C_1dW^1_t-\tilde{\Sigma}f^\top_2(K_2^\top)^{-1}dW^2_t,\\
-d\tilde{\tilde{\vp}}^{1i}_t&=\Big\{\Big[A-B_{1i}R_{1i}^{-1}\big(B_{1i}^\top P^{1i}+S_{1i}\big)\Big]\tilde{\tilde{\vp}}^{1i}+\check{\ga}^{1i}K^{-1}_2f_2\tilde{\tilde{X}}\Big\}dt-\la^{1i}_tdW^1_t+\check{\ga}^{1i}_tdW^2_t,\\
\tilde{\tilde{X}}_0&=0,\ \tilde{\tilde{\vp}}^{1i}_T=0,\quad i=1,2,\dots,N,
\end{aligned}
\right.
\end{equation}
where processes $\tilde{\tilde{X}}$ and $\tilde{\tilde{\vp}}^{1i}$ are independent of $u_2$. Indeed, the dependence of original state processes $X^{u_2},\vp^{1i,u_2}$ and $\la^{1i,u_2}$ on $u_2$ is due to the explicit appearance of these terms, $B_2u_2,P^{1i}B_2u_2$, in the drift and generator of FBSDE \eqref{leader_state}, and this is the only source. However, there is no such terms in the difference similar to the offset as we discussed for $\tilde{X}$ before, which implies that the filtering error processes $\tilde{\tilde{X}},\hat{X}-\check{X},\tilde{\tilde{\vp}}^{1i},\la^{1i}$ and $\check{\ga}^{1i}$ are independent of $u_2$ in \eqref{leader new error}. Therefore, $\tilde{\tilde{J}}_2(x_0)$ is independent of $u_2$.
\end{remark}
\begin{remark}
In the following discussion, the LQ optimal control theory of FBSDEs is necessary only if the control system is considered in the partially observed case on condition that the partially observed control is transformed into the completely observed case due to orthogonal decomposition of state process. The difference is that the fully coupled multi-dimensional system with inhomogeneous term should be considered because of the insert of optimal control of followers. This is always the case whether the followers own the full information or not when the leader can not observe the state system directly. However, another setting may have more interest that all the followers only know partial information but the leader knows all the information, in this case, the result in \cite{HJX23} is no longer applicable, now this question essentially boils down to a LQ stochastic optimal control of fully coupled FBSDEs with general conditional mean-field terms. The interesting extension has been studied in our another working paper while the decoupling technique of mean-filed LQ control are involved.
\end{remark}
Apply the stochastic LQ optimal control theory in the Appendix \ref{FBSLQ}, we now complete the leader's problem. For convenience, we rewrite \eqref{leader new state} as follows:
\begin{equation}\label{leader new state2}
\left\{
\begin{aligned}
d\check{X}_t&=\bigg\{\Big[A-\sum_{i=1}^NB_{1i}R_{1i}^{-1}\big(B_{1i}^\top P^{1i}+S_{1i}\big)\Big]\check{X}-\sum_{i=1}^NB_{1i}R_{1i}^{-1}B_{1i}^\top\check{\vp}^{1i}+B_2u_2\\
&\qquad-\sum_{i=1}^NB_{1i}R_{1i}^{-1}r_{1i}+\al\bigg\}dt+\sum_{j=1}^{l_2}\big[\tilde{\Sigma}f^\top_2(K_2^\top)^{-1}+C_2\big]^jd\check{U}^{j}_t,\\
-d\check{\vp}^{1i}_t&=\bigg\{\Big[A-B_{1i}R_{1i}^{-1}\big(B_{1i}^\top P^{1i}+S_{1i}\big)\Big]\check{\vp}^{1i}+P^{1i}\al+q_{1i}\\
&\qquad -R_{1i}^{-1}\big(B_{1i}^\top P^{1i}+S_{1i}\big)r_{1i}+P^{1i}B_2u_2\bigg\}dt-\sum_{j=1}^{l_2}\check{\ga}^{1ij}_td\check{U}^{j}_t,\\
\check{X}_0&=x_0,\ \check{\vp}^{1i}_T=g_{1i},\ i=1,2,\dots,N,
\end{aligned}
\right.
\end{equation}

Now, for $i=1,\cdots,N$ and $j=1,\cdots,l_2$, we have
\begin{equation*}
\begin{aligned}
&A_1=A-\sum_{i=1}^NB_{1i}R_{1i}^{-1}\big(B_{1i}P^{1i}+S_{1i}\big),\ \tilde{B}_{1i}=-B_{1i}R_{1i}^{-1}B_{1i}^\top,\ E_1=-\sum_{i=1}^NB_{1i}R_{1i}^{-1}r_{1i}+\al,\\
&D_1=B_2,\ E_2=\tilde{\Sigma}f_2^\top(K_2^\top)^{-1}+C_2,\ B^\top_{3i}=A-B_{1i}R_{1i}^{-1}\big(B^\top_{1i}P^{1i}+S_{1i}\big),\quad D_{3i}=P^{1i}B_2,\\
&E_{3i}=P^{1i}\al+q_{1i}-R_{1i}^{-1}\big(B_{1i}P^{1i}+S_{1i}\big)r_{1i},\ A_4=2Q_2,\quad D_4=2R_2,\ G=2G_2,\\
&\xi_i=g_{1i},\ A^j_2=A_{3i}=B^j_{2i}=B_{4i}=C^j_{1i}=C^{jj}_{2i}=C^j_{3i}=C^j_{4i}=D^j_2=H_i=F_i=0,\\
&L^j_1=I_{Nn\times Nn},\ L^j_2=I_{n\times n},\ L^j_3=0_{n\times n},\ L^j_4=0_{n\times Nn},\ L_5=2R_2,\ L_6=-(2R_2)^{-1}B_2P_1,\\
&L_7=-(2R_2)^{-1}B^\top_2((P^1)^\top+P_2),\ L^j_8=0_{n\times n},\ L^j_9=0_{n\times Nn},\ L^j_{10}=0_{Nn\times n},\\
&L^j_{11}=0_{Nn\times Nn},\ S^j_1=P_1,\ S^j_2=S^j_4=\check{V}^j_1+P_1\big[\tilde{\Sigma}f_2^\top(K_2^\top)^{-1}+C_2\big]^j,\\
&S_3=-(2R_2)^{-1}B^\top_2\check{\vp}_1,\ S^j_5=\check{V}^j_2+P_2^\top\big[\tilde{\Sigma}f_2^\top(K_2^\top)^{-1}+C_2\big]^j.
\end{aligned}
\end{equation*}
Thus, the optimal control of the leader is
\begin{equation}\label{optimal control of leader}
\begin{aligned}
\bar{u}_2(t)=&-(2R_2)^{-1}\big[B_2P_1+B_2^\top\big((P^1)^\top+P_2\big)P_3^{-1}P_2^\top\big]\bar{\check{X}}(t)\\
&+(2R_2)^{-1}B_2^\top\big((P^1)^\top+P_2\big)P_3^{-1}\bar{\check{\vp}}^1(t)\\
&-(2R_2)^{-1}B_2^\top\big((P^1)^\top+P_2\big)P_3^{-1}\check{\vp}_2(t)-(2R_2)^{-1}B_2^\top\check{\vp}_1(t),
\end{aligned}
\end{equation}
where $\bar{\check{X}}\equiv\check{X}^{\bar{u}_2}$ and $\bar{\check{\vp}}^1=\check{\vp}^{1,\bar{u}_2}$,
   $P^1\equiv(P^{11},P^{12},\cdots,P^{1N})^\top$
satisfies \eqref{P ^ 1i}, $P_1,P_2,P_3$ satisfy
\begin{equation*}
\left\{
\begin{aligned}
&\dot{P}_1+P_{1}\bigg[A-\sum_{i=1}^NB_{1i}R_{1i}^{-1}\big(B_{1i}P^{1i}+S_{1i}\big)+\tilde{B}_1^\top P_2^\top\bigg]\\
&\ +\bigg[A^\top-\sum_{i=1}^N\big(P^{1i}B_{1i}^\top+S_{1i}^\top\big)R_{1i}^{-1}B_{1i}^\top+P_2\tilde{B}_1\bigg]P_{1}-P_1B_2(2R_2)^{-1}B_2P_1+2Q_2=0,\\
&P_1(T)=2G_2,
\end{aligned}
\right.
\end{equation*}
\begin{equation*}
\hspace{-40mm}\left\{
\begin{aligned}
&\dot{P}_2-P_1\tilde{B}_1^\top P_3-P_1B_2(2R_2)^{-1}B_2^\top\big((P^1)^\top+P_2\big)+P_2B_3\\
&\ +\bigg[A^\top-\sum_{i=1}^N\big(P^{1i}B_{1i}^\top+S^\top_{1i}\big)R_{1i}^{-1}B_{1i}^\top+P_2\tilde{B}_1\bigg]P_2=0,\\
&P_2(T)=0,
\end{aligned}
\right.
\end{equation*}
\begin{equation*}
\left\{
\begin{aligned}
&\dot{P}_3+\big(P_2^\top \tilde{B}_1^\top+B_3^\top\big)P_{3}+P_3\big(\tilde{B}_1P_2+B_3\big)-\big(P^1+P_2^\top\big)B_2(2R_2)^{-1}B_2^\top\big((P^1)^\top+P_2\big)=0,\\
&P_3(T)=0,
\end{aligned}
\right.
\end{equation*}
respectively, and $(\check{\vp}_1,\check{V}_1)$, $(\check{\vp}_2,\check{V}_2)$ satisfy
\begin{equation*}
\left\{
\begin{aligned}
d\check{\vp}_1(t)&=-\bigg\{\bigg[P_2B_1+A^\top-\sum_{i=1}^N\big(P^{1i}B_{1i}^\top+S^\top_{1i}\big)R_{1i}^{-1}B_{1i}^\top-P_1B_2(2R_2)^{-1}B_2^\top\bigg]\check{\vp}_1\\
&\qquad +P_1B_1^\top\check{\vp}_{2}+P_1\bigg(-\sum_{i=1}^NB_{1i}R_{1i}^{-1}r_{1i}+\al\bigg)\bigg\}dt+\sum_{j=1}^{l_2}\check{V}^j_1d\check{U}^j_t,\\
d\check{\vp}_2(t)&=-\Big\{\big(P_2^\top B_1^\top+B_3^\top\big)\check{\vp}_2-\big[P_3B_1+\big(P^1+P_2^\top\big)B_2(2R_2)^{-1}B_2^\top\big]\check{\vp}_1\\
&\qquad +P_2^\top E_1+E_3\Big\}dt+\sum_{j=1}^{l_2}\check{V}^j_2d\check{U}^j_t,\\
\check{\vp}_1(T)&=0,\ \check{\vp}_2(T)=g_1.
\end{aligned}
\right.
\end{equation*}
The optimal cost of $\check{J}_2(x_0;u_2)$ is
\begin{equation*}
\check{J}_2(x_0;\bar{u}_2)=\frac{1}{2}\EE\bigg[\RR_2+\int_0^TM_5dt\bigg],
\end{equation*}
where $\RR_2=\big(X_{0}+P_{1}^{-1}(0) \check{\varphi}_{1}(0)\big)^\top P_{1}(0)\big(X_{0}+P_{1}^{-1}(0) \check{\varphi}_{1}(0)\big)$, and the explicit representation of $M_5$ is omitted, which indeed take values when insert the notations in this section into \eqref{M5}, with $\check{\tilde{\vp}}_1=-P_1^{-1}\check{\vp}_1, \check{\tilde{\vp}}_2=-P_2^\top P_1^{-1}\check{\vp}_1+\check{\vp}_2,
\check{\tilde{V}}^j_1=-P_1^{-1}\check{V}^j_1, \check{\tilde{V}}^j_2=-P_2^\top P_1^{-1}\check{V}^j_1+\check{V}^j_2$,
\begin{equation*}
\begin{aligned}
\tilde{P}&=\begin{pmatrix}P_1+P_2P_3^{-1}P_2^\top&-P_2P_3^{-1}\\-P_3^{-1}P_2^\top&P_3^{-1}\end{pmatrix}.
\end{aligned}
\end{equation*}
where the invertibility of $P_3$ can be discussed similarly as the counterpart in \cite[Section 4.1-4.3]{HJX23} by constructing a monotonous approximate sequence. So we omit the detailed derivations for avoiding repetition, and we still sketch the proof for the readers' convenience and completeness of the results in the Appendix \ref{FBSLQ}.

Next, we deal with $\tilde{\tilde{J}}_2(x_0)$. First, we need the equation of $\tilde{\tilde{X}}$:
\begin{equation*}
\begin{aligned}
d\tilde{\tilde{X}}_t&=\bigg\{A\tilde{\tilde{X}}-\sum_{i=1}^NB_{1i}R_{1i}^{-1}\big(B_{1i}^\top P^{1i}+S_{1i}\big)(\hat{X}-\check{X})-\sum_{i=1}^NB_{1i}R_{1i}^{-1}B_{1i}^\top(\vp^{1i}-\check{\vp}^{1i})\\
&\qquad -\big(\tilde{\Sigma}f_2^\top+C_2K_2^\top)(K_2K_2^\top)^{-1}f_2\tilde{\tilde{X}}\bigg\}dt +\sum_{j=1}^{l_1}C^j_1dW^{1j}_t-\sum_{j=1}^{l_2}\big[\tilde{\Sigma}f_2^\top(K_2^\top)^{-1}\big]^jdW^{2j}_t,
\end{aligned}
\end{equation*}
and introduce the ODE as follows
\begin{equation*}
\left\{
\begin{aligned}
&\dot{p}+p\big[A-\big(\tilde{\Sigma}f_2^\top+C_2K_2^\top)(K_2K_2^\top)^{-1}f_2\big]\\
&\quad +\big[A-\big(\tilde{\Sigma}f_2^\top+C_2K_2^\top)(K_2K_2^\top)^{-1}f_2\big]^\top p+Q_2=0,\\
&p_T=G_2.
\end{aligned}
\right.
\end{equation*}
Apply It\^o's formula to $\langle p\tilde{\tilde{X}},\tilde{\tilde{X}}\rangle$ and put it into $\tilde{\tilde{J}}_2(x_0)$, we obtain
\begin{equation*}
\begin{aligned}
\tilde{\tilde{J}}_2(x_0)=&\int_0^T\bigg\{-p\sum_{i=1}^NB_{1i}R_{1i}^{-1}\big(B_{1i}^\top P^{1i}+S_{1i}\big)\Sigma_{15}-\Sigma_{15}\sum_{i=1}^N\big(B_{1i}^\top P^{1i}+S_{1i}\big)^\top R_{1i}^{-1}B_{1i}^\top p\\
&\qquad -p\sum_{i=1}^NB_{1i}R_{1i}^{-1}B_{1i}^\top\check{\Sigma}^i-\sum_{i=1}^N(\check{\Sigma}^i)^\top B_{1i}R_{1i}^{-1}B_{1i}^\top p+\sum_{j=1}^{l_1}(C^j_1)^\top pC^j_1\\
&\qquad +\sum_{j=1}^{l_2}\big[(\tilde{\Sigma}f_2^\top(K_2^\top)^{-1})^j\big]^\top (\tilde{\Sigma}f_2^\top(K_2^\top)^{-1})^j\bigg\}dt.
\end{aligned}
\end{equation*}
Finally, we have the optimal cost of the leader as follows:
\begin{equation*}
\begin{aligned}
\tilde{J}_2(\bar{u}_2)&=\frac{1}{2}\RR_2+\frac{1}{2}\int_0^T\bigg\{\EE[M_5]-2p\sum_{i=1}^NB_{1i}R_{1i}^{-1}\big(B_{1i}^\top P^{1i}+S_{1i}\big)\Sigma_{15}\\
&\quad-2\Sigma_{15}\sum_{i=1}^N\big(B_{1i}^\top P^{1i}+S_{1i}\big)^\top R_{1i}^{-1}B_{1i}^\top p-2p\sum_{i=1}^NB_{1i}R_{1i}^{-1}B_{1i}^\top\check{\Sigma}^i\\
&\quad-2\sum_{i=1}^N(\check{\Sigma}^i)^\top B_{1i}R_{1i}^{-1}B_{1i}^\top p+2\sum_{j=1}^{l_1}(C^j_1)^\top pC^j_1\\
&\quad+2\sum_{j=1}^{l_2}\big[(\tilde{\Sigma}f_2^\top(K_2^\top)^{-1})^j\big]^\top (\tilde{\Sigma}f_2^\top(K_2^\top)^{-1})^j\bigg\}dt,
\end{aligned}
\end{equation*}
where $\Sigma_{15}$ satisfies the following ODE
\begin{equation}\label{Sigma15}
\begin{aligned}
&\Sigma_{15}=\int^t_0\Big\{\Sigma_{15}[A-(\tilde{\Sigma}f^\top_2(K^\top_2)^{-1}+C_2)K^{-1}_2f_2]^\top+[A-\sum^N_{i=1}B_{1i}R^{-1}_{1i}(B^\top_{1i}P^{1i}+S_{1i})]\Sigma_{15}\\
&-\Sigma_{12}[\sum^N_{i=1}B_{1i}R^{-1}_{1i}(B^\top_{1i}P^{1i}+S_{1i})]^\top-\sum^N_{i=1}\Sigma^i_{13}(B_{1i}R^{-1}_{1i}B^\top_{1i})^\top-\sum^N_{i=1}B_{1i}R^{-1}_{1i}B^\top_{1i}\check{\Sigma}^i+C_1C_1^\top\Big\}ds
\end{aligned}
\end{equation}
 and $\tilde{\Sigma},\check{\Sigma}^i,\Sigma_{12},\Sigma^i_{13}$ can be solved by equations \eqref{tilde Sigma}-\eqref{check varphi} similarly.

\section{Application to multi-agent formation control}\label{s5}

There are increasingly interest for the scholars using differential game methodology to model and analyze the optimization in multi-agent formation control problem. Gu \cite{Gu08} formulated the LQ leader-follower formation control where the Nash differential game theory are applied for the agents through the use of graph theory and an open-loop Nash equilibrium solution are investigated. Lin \cite{Lin14} considered a formation control problem for a multiple-UAV system, and the formation control problem is formulated and solved as a differential game problem due to the optimization of UAVs' differential objectives based on its local information. Mylvaganam and Astolfi \cite{MA15} studied a multi-agent formation system consisting of one leader and a group of followers, in which the leader is steered towards a target position while the followers seek to form and maintain a pre-defined formation about the leader, while avoiding collisions, is posed as a nonlinear differential game. In particular, the problem simplifies to a LQ differential game when collision avoidance is not taken into account. Moreover, some very foundational insights are provided in \cite{LZW20} which investigates the differential private algorithm for the average output consensus control of continuous-time heterogeneous multi-agent systems, then a sequence of further interesting new control and game results are first developed in \cite{WZH22}, which are of great inspiration for our future work. Especially, By making interesting connections to \cite{LZW20,WZH22}, it is worth noting that the consideration and introduction of key ingredients including differentially private, information exchange would motivate more practical applications in our paper and match our hierarchical feature in terms of asymmetrical roles, and try to apply the stochastic  Stackelberg differential game theory into the problem \cite{LZW20,WZH22}. Therefore, in this paper we aim to extend the existing formation control to fall within the stochastic Stackelberg game problem and develop the partially observed stochastic case with the introduction of hierarchical level and observation processes.
\subsection{Model Formulation}
A team has $N+1$ robots including one leader robot and $N$ follower robots, each of which is described by its dynamics in the following. For the single leader robot with $n$-dimensional coordinates $q^2\in\RR^n$, the state and control vectors are $X^2_t=((q^2_t)^\top,(\dot{q}^2_t)^\top)^\top\in\RR^{2n}$ and $u_2\in\RR^m$ and similarly for the follower robots with state and control vectors being $X^{1i}_t=((q^{1i}_t)^\top,(\dot{q}^{1i}_t)^\top)^\top\in\RR^{2n}$ and $u_{1i}\in\RR^m$ for $i=1,\dots,N$.

The dynamics of the leader robot and the follower robots are given as follows:
\begin{equation*}
\begin{aligned}
dX^2_t&=\big[aX^2_t+bu_2(t)\big]dt+C_1^2dW^1_t+C_2^2dW^2_t,\\
dX^{1i}_t&=\big[aX^{1i}_t+bu_{1i}(t)\big]dt+C_1^{1i}dW^1_t+C_2^{1i}dW^2_t,\ i=1,\cdots,N,
\end{aligned}
\end{equation*}
where
\begin{equation*}
\begin{aligned}
a=\begin{pmatrix}0&I_{n\times n}\\0&0\end{pmatrix}\in\RR^{2n\times2n},\quad b=\begin{pmatrix}0\\I_{n\times m}\end{pmatrix}\in\RR^{2n\times m},
\end{aligned}
\end{equation*}
and $C_1^2, C_1^{11}, \dots, C_1^{1N}\in\RR^{2n\times l_1}$, $C_2^2, C_2^{11}, \dots, C_2^{1N}\in\RR^{2n\times l_2}$.

Concatenating the states of all $N+1$ robots in a team into a vector:
\begin{equation*}
\begin{aligned}
X&=\begin{pmatrix}X^2\\X^{11}\\ \vdots\\X^{1N}\end{pmatrix}\in\RR^{2(N+1)n},\quad
C_1=\begin{pmatrix}C_1^2\\C_1^{11}\\ \vdots\\C_1^{1N}\end{pmatrix}\in\RR^{2(N+1)n\times l_1},\quad
C_2=\begin{pmatrix}C_2^2\\C_2^{11}\\ \vdots\\C_2^{1N}\end{pmatrix}\in\RR^{2(N+1)n\times l_2},\\
\tilde{X}&=\begin{pmatrix}X\\I_{2n}\end{pmatrix}\in\RR^{2(N+2)n},\quad
\tilde{C}_1=\begin{pmatrix}C_1\\0\end{pmatrix}\in\RR^{2(N+2)n\times l_1},\quad
\tilde{C}_2=\begin{pmatrix}C_2\\0\end{pmatrix}\in\RR^{2(N+2)n\times l_2},
\end{aligned}
\end{equation*}
where $I_{2n}\equiv[1,1,\dots,1]^\top$ is a $\RR^{2n}$-valued column vector and let
\begin{equation*}
\begin{aligned}
A&=\begin{pmatrix}I_{(N+1)\times(N+1)}& \\ &0\end{pmatrix}\otimes a\in\RR^{2(N+2)n\times 2(N+2)n},\\
B_2&=\begin{pmatrix}1\\0\\ \vdots\\0\\0\end{pmatrix}\otimes b\in\RR^{2(N+2)n\times m},\quad
B_{1i}=\begin{pmatrix}0\\1\\0\\ \vdots\\0\end{pmatrix}\otimes b\in\RR^{2(N+2)n\times m},\ i=1,\cdots, N,
\end{aligned}
\end{equation*}
where $\begin{pmatrix}0\ 1\ 0\ \cdots\ 0\end{pmatrix}^\top\in\RR^{N+2}$ denotes the vector that the remaining components are all 0 except for the $(i+1)$th component being 1. The operator $\otimes$ is the Kronecker product.

Then we have
\begin{equation*}
d\tilde{X}_t=\bigg[A\tilde{X}_t+B_2u_2(t)+\sum_{i=1}^{N}B_{1i}u_{1i}(t)\bigg]dt+\tilde{C}_1dW^1_t+\tilde{C}_2dW^2_t.
\end{equation*}
Let
\begin{equation*}
\begin{aligned}
X^{2,d}=\begin{pmatrix}q^{2,d}\\ \dot{q}^{2,d}\end{pmatrix},\quad X^{1i,d}=\begin{pmatrix}q^{1i,d}\\ \dot{q}^{1i,d}\end{pmatrix}
\end{aligned}
\end{equation*}
be the desired state vector for robots. Similarly,
$X^d=((X^{2,d})^\top,(X^{11,d})^\top, \cdots,(X^{1N,d})^\top)^\top$,
where
\begin{equation*}
\begin{aligned}
dX^{2,d}_t&=\big[aX^{2,d}_t+bu_2^d(t)\big]dt,\\
dX^{1i,d}_t&=\big[aX^{1i,d}_t+bu_{1i}^d(t)\big]dt,\ i=1,\cdots,N,
\end{aligned}
\end{equation*}
and the concatenating state equations is
\begin{equation*}
d\tilde{X}^d_t=\big[A\tilde{X}^d_t+B_2u_2^d(t)+\sum_{i=1}^{N}B_{1i}u_{1i}^d(t)\big]dt.
\end{equation*}

Now we use the graph to represent the formation control interconnection between robots based on the original work of Fax and Murray \cite{FM04} which provides a link between graph theory and the formation control for a given communication topology. A vertex of the graph corresponds to a robot and the edges of the graph capture the dependence of the interconnections.
Formally, a directed graph $\mathcal{G}=(\mathcal{V},\mathcal{E})$ consists of a set of vertices $\mathcal{V}=\{v_{11},\dots,v_{1N},v_2\}$, indexed by the robots in a team, and a set of edges $\mathcal{E}=\{(v_i,v_j)\in\mathcal{V}\times\mathcal{V}\}$, containing ordered pairs of distinct vertices. Assuming the graph has no loops, i.e., $(v_i,v_j)\in\mathcal{E}$ implies $v_i\neq v_j$. A graph is connected if for any vertices $v_i,v_j\in\mathcal{V}$, there exists a path of edges in $\mathcal{E}$ from $v_i$ to $v_j$. An edge-weighted graph is a graph in which each edge is assigned a weight. The edge $(v_i,v_j)$ is associated with weight $w_{ij}\geq 0$. To control a team to keep a formation, the graph should be assumed to be connected.

The incidence matrix $D$ of a directed graph $\mathcal{G}$ is the $\{0,\pm1\}$-matrix with rows and columns indexed by vertices of $\mathcal{V}$ and edges of $\mathcal{E}$, respectively, such that the $uv$th entry of $D$ is equal to 1 if the vertex $u$ is the head of the edge $v$, -1 if the vertex $u$ is the tail of the edge $v$, and 0, otherwise. If graph $\mathcal{G}$ has $m$ vertices and $|\mathcal{E}|$ edges, then incidence matrix $D$ of the graph $\mathcal{G}$ has order $(N+1)\times|\mathcal{E}|$.

The cohesion and separation of formation control is defined by the desired distance vector $d^d_{ij}=X^{i,d}-X^{j,d}$ between two neighbors $v_i$ and $v_j$ in $\{v_{11},v_{12},\cdots,v_{1N},v_2\}$. The formation error vector is defined as $X^i-X^j-d_{ij}^d$ for edge $(v_i,v_j)$.

Next, we consider the representation of the whole team formation error in a matrix form as follows:
\begin{equation*}
\begin{aligned}
&\sum_{(i,j)\in\mathcal{E}}w_{ij}\Vert X^i-X^j-d_{ij}^d\Vert^2=\sum_{(i,j)\in\mathcal{E}}w_{ij}\Vert X^i-X^j\Vert^2+w_{ij}\Vert d_{ij}^d\Vert^2-2\big\langle w_{ij}(X^i-X^j),d_{ij}^d\big\rangle\\
&=X^\top\hat{D}\hat{W}\hat{D}^\top X+(d^d)^\top\hat{W}d^d-X^\top\hat{D}\hat{W}d^d-(\hat{D}\hat{W}d^d)^\top X\triangleq\tilde{X}^\top\tilde{L}\tilde{X},
\end{aligned}
\end{equation*}
where
\begin{equation*}
\tilde{L}=\begin{pmatrix}\hat{L}&-\hat{D}\hat{W}d^dI_{2n}^\top\\-\big(\hat{D}\hat{W}d^dI_{2n}^\top\big)^\top&\frac{1}{4n^2}(d^dI_{2n}^\top)^\top \hat{W}d^dI_{2n}^\top\end{pmatrix}
\end{equation*}
and $\hat{W}=W\otimes I_{2n\times2n}$, with $W=diag[w_{ij}]$ being a diagonal weight matrix with dimension $|\mathcal{E}|$. Define the Laplacian of a graph $\mathcal{G}$ as $L=DWD^\top$ which is symmetric and positive semi-definite. Based on $(X\otimes Y)^\top=(X^\top\otimes Y^\top)$ and $(X\otimes Y)(U\otimes V)=(XU)\otimes(YV)$, we have $\hat{L}=\hat{D}\hat{W}\hat{D}^\top=L\otimes I_{2n\times2n}$ which is also symmetric and positive semi-definite.

In this case, we consider a more challenging formation control problem that all the robots have their own individual cost functional, in fact, they can choose their objectives based on relative displacement and velocity errors with the immediate neighbours in the graph topology.  Now we give the cost functional of the follower robots, for $i=1,\cdots,N$,
\begin{equation}\label{J_i}
J_{1i}(u_{1i},u_2)=\EE\bigg[\sum_{(i,j)\in\mathcal{E}}w_{ij}\Vert X_T^i-X_T^j-d_{ij}^d\Vert^2+\int_0^T\sum_{(i,j)\in\mathcal{E}}\mu_{ij}\Vert X_t^i-X_t^j-d_{ij}^d\Vert^2+u_{1i}^\top R_{1i}u_{1i}dt\bigg],
\end{equation}
where $\mu_{ij},R_{1i}>0$.

Let $K_{1i}=\hat{L}_{1i}=\hat{D}\hat{W}_{1i}\hat{D}^\top$, $\hat{W}_{1i}=W_{1i}\otimes I_{2n\times 2n}$, $W_{1i}=diag[w_{ij}]$, $Q_{1i}=\hat{L}_{i}=\hat{D}\hat{W}_{i}\hat{D}^\top$, $\hat{W}_{i}=W_{i}\otimes I_{2n\times2n}$, $W_{i}=diag[\mu_{ij}]$. $K_{1i}$ and $Q_i$ are symmetric and positive semi-definite. Then \eqref{J_i} becomes
\begin{equation*}
J_{1i}(u_{1i},u_2)=\EE\bigg[\tilde{X}_T^\top\tilde{K}_{1i}\tilde{X}_T+\int_0^T\tilde{X}^T\tilde{Q}_i\tilde{X}+u_{1i}^\top R_{1i}u_{1i}dt\bigg],
\end{equation*}
where
\begin{equation*}
\begin{aligned}
\tilde{K}_{1i}&=\begin{pmatrix}K_{1i}&-\hat{D}\hat{W}_{1i}d^dI_{2n}^\top\\-(\hat{D}\hat{W}_{1i}d^dI_{2n}^\top)^\top&\frac{1}{4n^2}(d^dI_{2n}^\top)^\top\hat{W}_{1i}d^dI_{2n}^\top\end{pmatrix},\\ \tilde{Q}_{1i}&=\begin{pmatrix}Q_{1i}&-\hat{D}\hat{W}_id^dI_{2n}^\top\\-(\hat{D}\hat{W}_id^dI_{2n}^\top)^\top&\frac{1}{4n^2}(d^dI_{2n}^\top)^\top\hat{W}_id^dI_{2n}^\top\end{pmatrix}.
\end{aligned}
\end{equation*}
To track a special trajectory $X^{2,d}$, the leader robot should track $X^{2,d}$, so its cost functional is
\begin{equation*}
\begin{aligned}
J_2(u_{1i},u_2)&=\EE\bigg[\sum_{(i,j)\in\mathcal{E}}\nu_{ij}\Vert X_T^i-X_T^j-d_{ij}^d\Vert+\Vert X_T^2-X_T^{2,d}\Vert_{k_2}\\
&\qquad +\int_0^T\sum_{(i,j)\in\mathcal{E}}\theta_{ij}\Vert X^i-X^j-d_{ij}^d\Vert+\Vert X_t^2-X_t^{2,d}\Vert_{\tilde{q}_2}+u_2^\top R_{22}u_2dt\bigg]\\
&=\EE\bigg[\tilde{X}^\top_T\big(\tilde{K}_2+\tilde{k}_2^1\big)\tilde{X}_T+\int_0^T\tilde{X}^\top\big(\tilde{Q}_2+\tilde{q}_2^1\big)\tilde{X}+u_2^\top R_{22}u_2dt\bigg]\\
&\triangleq\EE\bigg[\tilde{X}^\top_T\tilde{K}_2^1\tilde{X}_T+\int_0^T\tilde{X}^\top\tilde{Q}_2^1\tilde{X}+u_2^\top R_{22}u_2dt\bigg],
\end{aligned}
\end{equation*}
where
\begin{equation*}
\begin{aligned}
\tilde{K}_2&=\begin{pmatrix}\hat{K}_2&-\hat{D}\hat{W}_2d^dI_{2n}^\top\\-\big(\hat{D}\hat{W}_2d^dI_{2n}^\top\big)^\top&\frac{1}{4n^2}(d^dI_{2n}^\top)^\top\hat{W}_2d^dI_{2n}^\top\end{pmatrix},\\ \tilde{k}_2^1&=\begin{pmatrix}k_2^1&-k_2^2X_T^{2,d}I_{2n}^\top\\-\big(k_2^2X_T^{2,d}I_{2n}^\top\big)^\top&\frac{1}{4n^2}(X_T^{2,d}I_{2n}^\top)^\top k_2X_T^{2,d}I_{2n}^\top\end{pmatrix}, \\
\tilde{Q}_2&=\begin{pmatrix}Q_2&-\hat{D}\hat{U}_2d^dI_{2n}^\top\\-\big(\hat{D}\hat{U}_2d^dI_{2n}^\top\big)^\top&\frac{1}{4n^2}(d^dI_{2n}^\top)^\top\hat{U}_2d^dI_{2n}^\top\end{pmatrix},\\ \tilde{q}_2^1&=\begin{pmatrix}q_2^1&-q_2^2X^{2,d}I_{2n}^\top\\-\big(q_2^2X^{2,d}I_{2n}^\top\big)^\top&\frac{1}{4n^2}(X^{2,d}I_{2n}^\top)^\top \tilde{q}_2X^{2,d}I_{2n}^\top\end{pmatrix},
\end{aligned}
\end{equation*}
in which $\hat{K}_2=\hat{L}_2=\hat{D}\hat{W}_2\hat{D}^\top$, $\hat{W}_2=V_2\otimes I_{2n\times 2n}$, $V_2=diag[\nu_{ij}]$, $Q_2=\hat{P}_2=\hat{D}\hat{U}_2\hat{D}^\top$, $\hat{U}_2=H_2\otimes I_{2n\times2n}$, $H_2=diag[\theta_{ij}]$. $\hat{K}_2$ and $Q_2$ are symmetric and positive semi-definite. $k_2=diag[w_2]$, $\tilde{q}_2=diag[\mu_2]$, and $k_2^1=diag[k_2,0,\dots,0], k_2^2=[k_2,0,\dots,0]^\top$, $q_2^1=diag[\tilde{q}_2,0,\dots,0], \\q_2^2=[\tilde{q}_2,0,\dots,0]^\top$. $\tilde{K}_2^1,\tilde{Q}_2^1$ are also symmetric and positive semi-definite. $R_{22}$ is a negative-definite matrix.

We assume that $\tilde{K}_2,\tilde{Q}_2=0$, i.e., the leader do not take the formation error into consideration and only track the desired trajectory and following him with a fixed distance. Moreover, a robot addresses its cost functional that is only related to part of other robots based on the partial information when fixing a formation shown by connected graph topology. Therefore, we set that all of them can not observe the complete information where the followers know more information than the leader. Then the observable information available to the followers and the leader are described by the following equations:
\begin{equation*}
dY_1^{u_1,u_2}(t)=\big[f_1(t)\tilde{X}_t+g_1(t)\big]dt+K_1(t)dW^1_t,
\end{equation*}
specially, $f_1(t)=g_1(t)=0,K_1(t)=I$ implies that $Y_1(t)=W_t^1$; and
\begin{equation*}
dY_2^{u_1,u_2}(t)=\big[f_2(t)\tilde{X}_t+g_2(t)\big]dt+K_2(t)dW_t^2.
\end{equation*}

Then apply the theoretical results in the previous sections \ref{follower pro} and \ref{leader pro}, and let $A=A, B_2=B_2, B_{1i}=B_{1i}, C_1=\tilde{C}_1,C_2=\tilde{C}_2,Q_{1i}=\tilde{Q}_{1i}, R_{1i}=R_{1i},G_{1i}=\tilde{K}_{1i}, Q_2=\tilde{Q}_2^1, R_2=R_{22}, G_2=\tilde{K}_2^1, S_{1i}=S_2=q_{1i}=q_2=r_{1i}=r_2=g_{1i}=g_2=\al=0$, then we obtain the optimal control of the $i$th follower robot:
\begin{equation*}
\bar{u}_{1i}(t)=-R_{1i}^{-1}(t)B_{1i}^\top(t) P^{1i}_t\hat{\tilde{X}}_t-R_{1i}^{-1}(t)B_{1i}^\top(t)\vp^{1i}_t,
\end{equation*}
where $\hat{\tilde{X}},P^{1i},\vp^{1i}$ satisfy
\begin{equation*}
\begin{aligned}
d\hat{\tilde{X}}_t&=\bigg[\Big(A(t)-\sum_{i=1}^NB_{1i}(t)R_{1i}^{-1}(t)B_{1i}^\top(t) P^{1i}_t\Big)\hat{\tilde{X}}_t-\sum_{i=1}^NB_{1i}(t)R_{1i}^{-1}(t)B_{1i}^\top(t)\vp^{1i}_t\\
&\quad+B_2(t)u_2(t)\bigg]dt+\tilde{C}_1(t)dW_t^1,
\end{aligned}
\end{equation*}
\begin{equation}\label{P1i}
\left\{
\begin{aligned}
&\dot{P}^{1i}_t+A^\top(t)P^{1i}_t+P^{1i}_tA(t)-P^{1i}_tB_{1i}(t)R_{1i}^{-1}(t)B^\top_{1i}(t)P^{1i}_t+\tilde{Q}_{i}(t)=0,\\
&P^{1i}_T=\tilde{K}_{1i},
\end{aligned}
\right.
\end{equation}
and
\begin{equation*}\left\{
\begin{aligned}
-d\vp^{1i}_t&=\big[\big(A(t)-B_{1i}(t)R_{1i}^{-1}(t)B_{1i}^\top(t) P^{1i}_t\big)^\top\vp^{1i}_t+P^{1i}_tB_2(t)u_2(t)\big]dt-\la^{1i}_tdW_t^1,\\
\vp^{1i}_T&=0,
\end{aligned}
\right.\end{equation*}
respectively.
Now, for any fixed $u_2$, the optimal cost of the follower is given by
\begin{equation*}
\begin{aligned}
J_{1i}(u_2)&=\EE\big[\langle P^{1i}_0x_0+2\vp^{1i}_0,x_0\rangle\big]+\EE\int_0^T\bigg[\sum_{j=1}^{l_2}\langle\Pi^{1i}\tilde{C}_{2j},\tilde{C}_{2j}\rangle+\sum_{j=1}^{l_1}\langle P^{1i}\tilde{C}_{1j},\tilde{C}_{1j}\rangle\\
&-\langle R_{1i}^{-1}B_{1i}^\top\vp^{1i},B_{1i}^\top\vp^{1i}\rangle+2\langle B_2u_2,\vp^{1i}\rangle+2\langle\tilde{C}_1,\la^{1i}\rangle\bigg]dt,
\end{aligned}
\end{equation*}
where
\begin{equation*}
\left\{
\begin{aligned}
&\dot{\Pi}^{1i}_t+A^\top\Pi^{1i}_t+\Pi^{1i}_tA+\tilde{Q}_i=0,\\
&\Pi^{1i}_T=\tilde{K}_{1i}.
\end{aligned}
\right.
\end{equation*}

Next, the optimal control of the leader robot has the feedback representation as follows:
\begin{equation*}
\begin{aligned}
\bar{u}_2(t)&=-(2R_{22})^{-1}\big[B_2^\top P_1+B_2^\top\big((P^1)^\top+P_2\big)P_3^{-1}P_2^\top\big]\check{\tilde{X}}_t\\
&\quad +(2R_{22})^{-1}B_2^\top\big((P^1)^\top+P_2\big)P_3^{-1}\check{\vp}^1_t,
\end{aligned}
\end{equation*}
where $P^{1i}$ satisfies \eqref{P1i} and $P_1,P_2,P_3$ satisfy
\begin{equation*}
\left\{
\begin{aligned}
&\dot{P}_1+P_1\bigg(A-\sum_{i=1}^NB_{1i}R_{1i}^{-1}B^\top_{1i}P^{1i}+\tilde{B}_1^\top P_2^\top\bigg)+\bigg(A^\top-\sum_{i=1}^NP^{1i}B_{1i} R_{1i}^{-1}B_{1i}^\top+P_2\tilde{B}_1\bigg)P_1\\
&\quad -P_1B_2(2R_{22})^{-1}B^\top_2P_1+2\tilde{Q}_2^1=0,\\
&P_1(T)=2\tilde{K}_2^1,
\end{aligned}
\right.
\end{equation*}
\begin{equation*}
\left\{
\begin{aligned}
&\dot{P}_2-P_1\tilde{B}_1^\top P_3-P_1B_2(2R_{22})^{-1}B_2^\top((P^1)^\top+P_2)+P_{2}B_3\\
&\quad +\bigg(A^\top-\sum_{i=1}^NP^{1i}B_{1i} R_{1i}^{-1}B_{1i}^\top+P_2\tilde{B}_1\bigg)P_2=0,\\
&P_2(T)=0,
\end{aligned}
\right.
\end{equation*}
\begin{equation*}
\left\{
\begin{aligned}
&\dot{P}_3+\big(P_2^\top \tilde{B}_1^\top+B_3^\top\big)P_3+P_3\big(\tilde{B}_1P_2+B_3\big)-(P^1+P_2^\top)B_2(2R_{22})^{-1}B_2^\top((P^1)^\top+P_2)=0,\\
&P_3(T)=0,
\end{aligned}
\right.
\end{equation*}
$\check{\tilde{X}}$ satisfies
\begin{equation*}\left\{
\begin{aligned}
d\check{\tilde{X}}_t&=\bigg\{\bigg[A-\sum_{i=1}^NB_{1i}R_{1i}^{-1}B_{1i}^\top P^{1i}-B_2(2R_{22})^{-1}\big(B^\top_2P_1+B_2^\top((P^1)^\top+P_2)P_3^{-1}P_2^\top\big)\bigg]\check{\tilde{X}}\\
&\qquad -\sum_{i=1}^NB_{1i}R_{1i}^{-1}B_{1i}^\top\check{\vp}^{1i}+B_2(2R_{22})^{-1}B_2^\top((P^1)^\top+P_2)P_3^{-1}\check{\vp}^1\bigg\}dt\\
&\quad +\big[\tilde{\Sigma}f^\top_2(K_2^\top)^{-1}+\tilde{C}_2\big]d\check{U}_t,\\
-d\check{\vp}^{1i}_t&=\bigg[\big(A-B_{1i}R_{1i}^{-1}B_{1i}^\top P^{1i}\big)\check{\vp}^{1i}+P^{1i}B_2(2R_{22})^{-1}B_2^\top((P^1)^\top+P_2)P_3^{-1}\check{\vp}^1\\
&\qquad -P^{1i}B_2(2R_{22})^{-1}\big(B^\top_2P_1+B_2^\top((P^1)^\top+P_2)P_3^{-1}P_2^\top\big)\check{\tilde{X}}\bigg]dt\\
&\quad -\check{\ga}^{1i}_td\check{U}_t,\ i=1,2,\dots,N,\\
d \check{\vp}_t^1&=-\bigg\{\big[-D_3D_4^{-1} D_3^\top P_3^{-1}P_2^\top-D_3D_4^{-1} D_1^\top(P_1+P_2P_3^{-1}P_2^\top)\big]\check{\tilde{X}}\\
&\qquad +\big[B_3^\top +D_3D_4^{-1} D_3^\top P_3^{-1}+D_3D_4^{-1} D_1^\top P_2P_3^{-1}\big]\check{\vp}^1\bigg\} dt+\check{\ga}_t^{1}d\check{U}_t, \\
\check{X}_0&=x_0,\ \check{\vp}^{1i}_T=0,
\end{aligned}
\right.\end{equation*}
and $(\check{\vp}_1,\check{V}_1)$, $(\check{\vp}_2,\check{V}_2)$ satisfy
\begin{equation*}\left\{
\begin{aligned}
d\check{\vp}_1&=-\bigg[\bigg(P_2B_1+A^\top-\sum_{i=1}^NP^{1i}B_{1i}^\top R_{ii}^{-1}B_{1i}^\top-P_1B_2(2R_{22})^{-1}B_2^\top\bigg)\check{\vp}_1\\
&\qquad +P_1B_1^\top\check{\vp}_2\bigg]dt+\sum_{j=1}^{l_2}\check{V}^j_1d\check{U}^j_t,\\
d\check{\vp}_2&=-\big[\big(P_2^\top B_1^\top+B_3^\top\big)\check{\vp}_2-\big(P_3B_1+(P^1+P_2^\top)B_2(2R_{22})^{-1}B_2^\top\big)\check{\vp}_1\big]dt+\sum_{j=1}^{l_2}\check{V}^j_2d\check{U}^j_t,\\
\check{\vp}_1(T)&=0,\ \check{\vp}_2(T)=0,
\end{aligned}
\right.\end{equation*}
which admit unique zero solutions.

Then the optimal cost of the leader is as follows:
\begin{equation*}
\begin{aligned}
J_2(\bar{u}_2)&=\frac{1}{2}\RR_2+\frac{1}{2}\int_0^T\bigg\{\EE[M_5]-2p\sum_{i=1}^NB_{1i}R_{1i}^{-1}B_{1i}^\top P^{1i}\Sigma_{15}-2\Sigma_{15}\sum_{i=1}^N(B_{1i}^\top P^{1i})^\top R_{1i}^{-1}B_{1i}^\top p\\
&\quad -2p\sum_{i=1}^NB_{1i}R_{1i}^{-1}B_{1i}^\top\check{\Sigma}^i-2\sum_{i=1}^N(\check{\Sigma}^i)^\top B_{1i}R_{1i}^{-1}B_{1i}^\top p+2\sum_{j=1}^{l_1}(\tilde{C}^j_1)^\top p\tilde{C}^j_1\\
&\quad +2\sum_{j=1}^{l_2}\big[(\tilde{\Sigma}f_2^\top(K_2^\top)^{-1})^j\big]^\top (\tilde{\Sigma}f_2^\top(K_2^\top)^{-1})^j\bigg\}dt,
\end{aligned}
\end{equation*}
where $\tilde{\tilde{\tilde{X}}}\triangleq\tilde{X}-\check{\tilde{X}}$ satisfies
\begin{equation*}
\begin{aligned}
d\tilde{\tilde{\tilde{X}}}_t&=\bigg[A\tilde{\tilde{X}}_t-\sum_{i=1}^NB_{1i}R_{1i}^{-1}B_{1i}^\top P^{1i}(\hat{\tilde{X}}_t-\check{\tilde{X}}_t)-\sum_{i=1}^NB_{1i}R_{1i}^{-1}B_{1i}^\top(\vp^{1i}_t-\check{\vp}^{1i}_t)\\
&\qquad -(\tilde{\Sigma}f_2^\top+\tilde{C}_2K_2^\top)(K_2K_2^\top)^{-1}f_2\tilde{\tilde{\tilde{X}}}_t\bigg]dt+\sum_{j=1}^{l_1}\tilde{C}^j_1dW^{1j}-\sum_{j=1}^{l_2}\big[\tilde{\Sigma}f_2^\top(K_2^\top)^{-1}\big]^jdW^{2j},
\end{aligned}
\end{equation*}
\begin{equation*}
\left\{
\begin{aligned}
&\dot{p}+p\big[A-(\tilde{\Sigma}f_2^\top+\tilde{C}_2K_2^\top)(K_2K_2^\top)^{-1}f_2\big]\\
&\quad +\big[A-(\tilde{\Sigma}f_2^\top+\tilde{C}_2K_2^\top)(K_2K_2^\top)^{-1}f_2\big]^\top p+\tilde{Q}_2^1=0,\\
&p_T=\tilde{K}_2^1,
\end{aligned}
\right.
\end{equation*}
with $\RR_2=X_0^\top P_1(0)X_0$, and $\tilde{\Sigma},\check{\Sigma}^i,\Sigma_{12},\Sigma^i_{13}$ and $\Sigma_{15}$ can be solved by equations \eqref{tilde Sigma}-\eqref{check varphi} and \eqref{Sigma15} with the corresponding values when replacing $X$ with $\tilde{X}$.

\subsection{Simulations}
In this section, we present a numerical example to illustrate the partially observed Stackelberg formation control strategy. We consider a planar multi-robot system with one leader and three followers. Hence, we consider 
    $n=2, m=2, N=3$.
For each robot, the position and velocity are denoted by
\begin{equation}
\begin{aligned}
    q^{1i}_t&=(q^{1i}_{x,t},q^{1i}_{y,t})^\top,\quad q^2_t=(q^2_{x,t},q^2_{y,t})^\top\\
    \dot q^{1i}_t&=(\dot q^{1i}_{x,t},\dot q^{1i}_{y,t})^\top,\quad \dot q^2_t=(\dot q^2_{x,t},\dot q^2_{y,t})^\top,i=1,2,3,
\end{aligned}
\end{equation}
and the state is
\begin{equation}
\begin{aligned}
    X^{1i}_t&=\big((q^{1i}_t)^\top,(\dot q^{1i}_t)^\top\big)^\top\in\mathbb R^{4},\\
    X^2_t&=\big((q^2_t)^\top,(\dot q^2_t)^\top\big)^\top\in\mathbb R^{4},i=1,2,3.
\end{aligned}
\end{equation}
The control input is two-dimensional,
\begin{equation}
\begin{aligned}
    u_{1i}(t)&=\big(u_{1i,x}(t),u_{1i,y}(t)\big)^\top\in\mathbb R^{2},\\
    u_2(t)&=\big(u_{2,x}(t),u_{2,y}(t)\big)^\top\in\mathbb R^{2},i=1,2,3.
\end{aligned}
\end{equation}
The whole team contains \(N+1=4\) robots, and the augmented state is consistent with the theoretical formulation
\[
    \tilde{X}
    =
    \big((X^2)^\top,(X^{11})^\top,(X^{12})^\top,(X^{13})^\top,
    I_{2n}^\top\big)^\top .
\]
The desired trajectories are generated from the augmented desired state
\[
    \tilde{X}^d
    =
    \big((X^{2,d})^\top,(X^{11,d})^\top,(X^{12,d})^\top,
    (X^{13,d})^\top,I_{2n}^\top\big)^\top ,
\]
where the physical reference components correspond to the desired planar paths of the leader and followers. The desired formation is a leader-centered V-shape, meaning that the relative positions
\[
    q^{1i}_t-q^2_t,\qquad i=1,2,3,
\]
are expected to approach prescribed V-shaped offsets.

We compare three cases: a deterministic case, a small-noise case, and an original-noise case. The deterministic case removes the Brownian perturbations and serves as a baseline. The small-noise and original-noise cases illustrate the robustness of the state feedback strategy under stochastic disturbances. The Monte Carlo sample paths are used to compute empirical means, standard deviations, and covariance quantities.
\begin{figure}[t]
    \centering
    \includegraphics[width=\columnwidth]{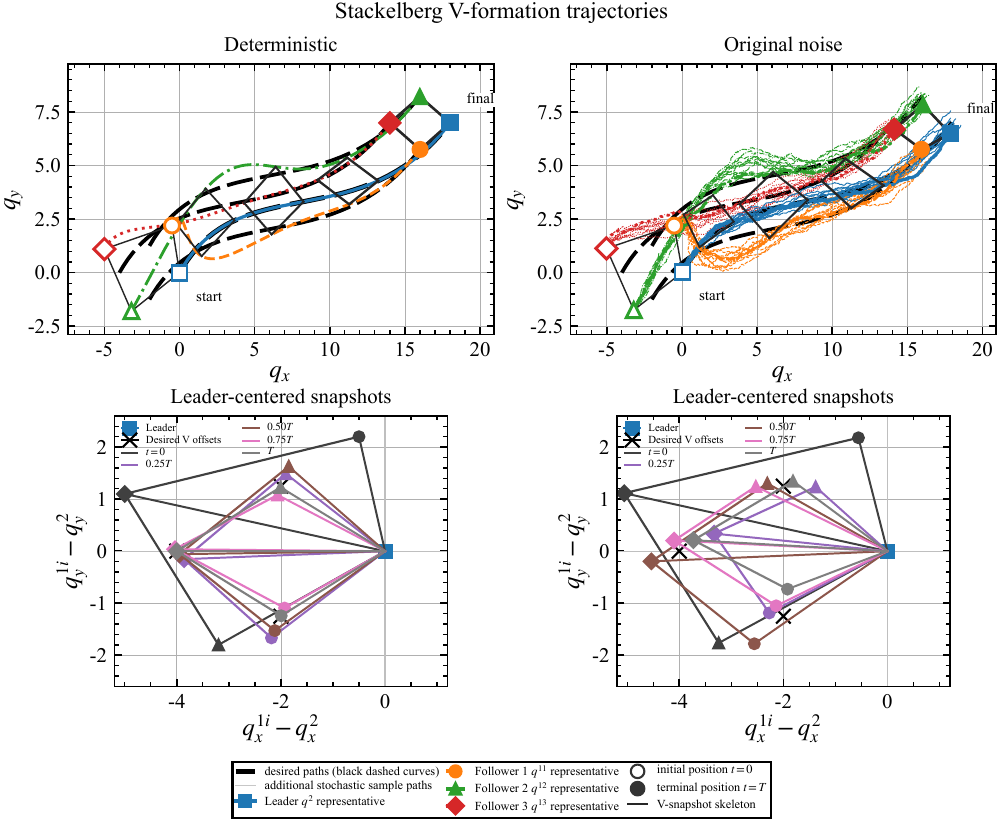}
    \caption{Stackelberg V-formation trajectories under deterministic and original-noise cases.
    The upper panels show physical trajectories,
    while the lower panels show leader-centered relative snapshots.}
    \label{fig:trajectories}
\end{figure}

Figure~\ref{fig:trajectories} shows the physical trajectories and the corresponding leader-centered snapshots. The robots are initialized away from the desired V-shape. In the deterministic case, the trajectories are smooth and the followers gradually approach the desired relative offsets. In the original-noise case, additional sample paths show visible stochastic fluctuations caused by Brownian perturbations. Nevertheless, the representative trajectory remains close to the desired paths, and the leader-centered snapshots show that the relative formation remains near the prescribed V-shape. This indicates that the Stackelberg state feedback strategy stabilizes the formation even under stochastic disturbances.

\begin{figure}[t]
    \centering
    \includegraphics[width=\columnwidth]{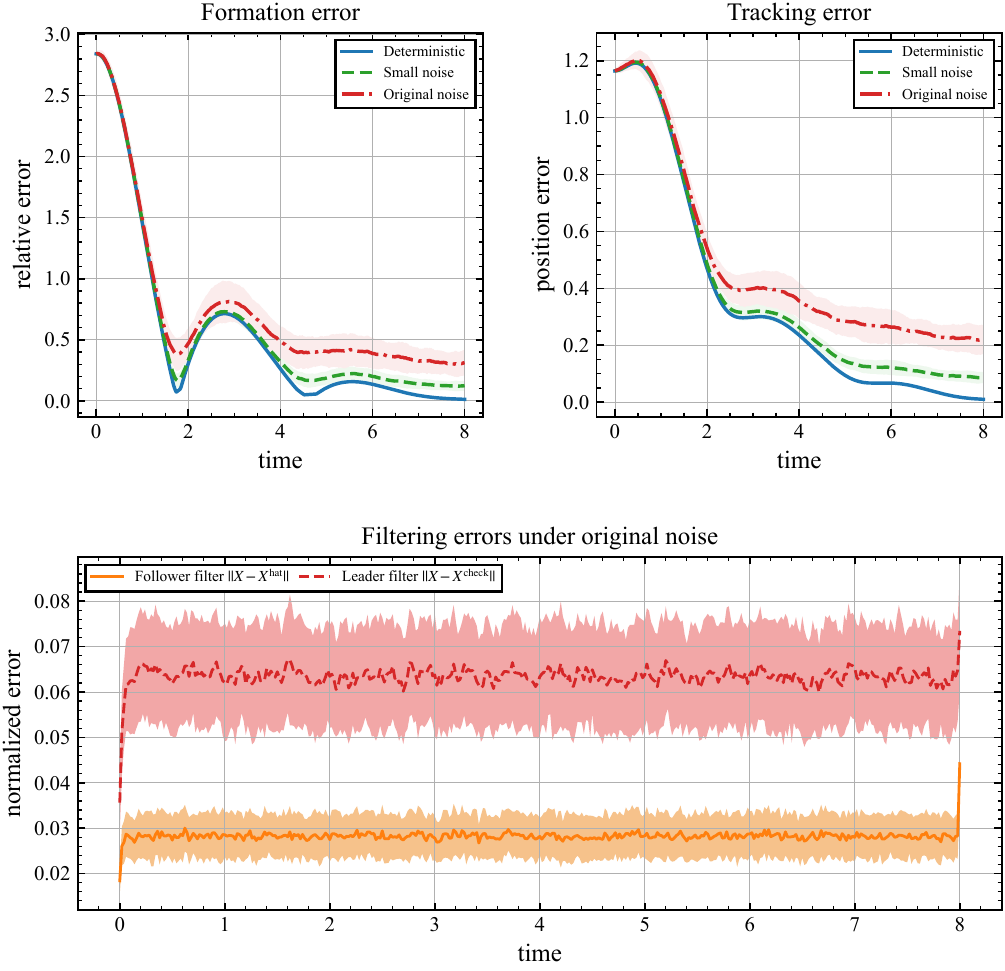}
    \caption{Formation, tracking, and filtering errors.
    The deterministic case achieves the smallest residual errors,
    whereas the stochastic cases exhibit bounded fluctuations.}
    \label{fig:errors}
\end{figure}

Figure~\ref{fig:errors} reports the formation, tracking, and filtering errors. The formation error measures the deviation of the relative configuration from the desired V-shape, while the tracking error measures the deviation from the desired augmented trajectory. Both errors decrease after an initial transient. The deterministic case gives the smallest residual error, whereas the original-noise case yields a larger but bounded residual error. This is consistent with stochastic stabilization, where the Stackelberg system with state feedback equilibrium remains in a neighborhood of the desired formation. The filtering errors remain bounded under the original-noise case. The leader-side filtering error (based on less information) is larger than the follower-side filtering error (based on more information), which reflects the information asymmetry in the partially observed Stackelberg setting.

\begin{figure}[t]
    \centering
    \includegraphics[width=\columnwidth]{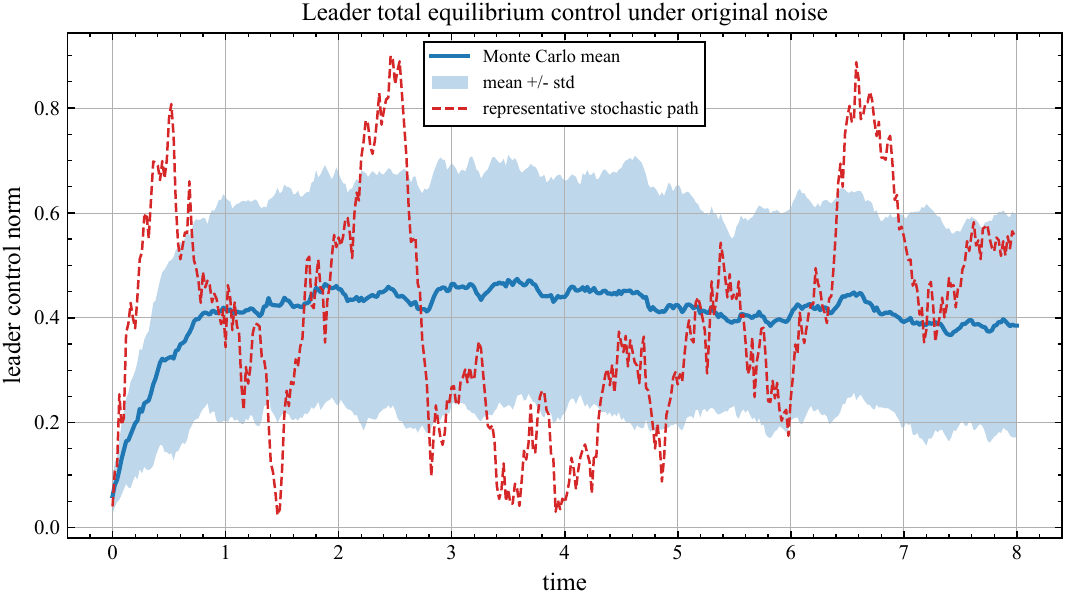}
    \caption{Leader total equilibrium control under the original-noise case.}
    \label{fig:leader-control}
\end{figure}

Figure~\ref{fig:leader-control} plots the norm of the leader equilibrium control. Since the control is vector-valued, its Euclidean norm is used to represent the instantaneous control effort, which is consistent with the quadratic control cost. The blue curve denotes the Monte Carlo mean, the shaded region represents one standard deviation, and the dashed curve is a representative stochastic realization. The representative path fluctuates around the mean due to Brownian perturbations and filtering-induced corrections. The bounded behavior of the control norm indicates that the leader does not require unbounded effort to maintain tracking and formation performance under the original-noise case.

\begin{figure}[t]
    \centering
    \includegraphics[width=\columnwidth]{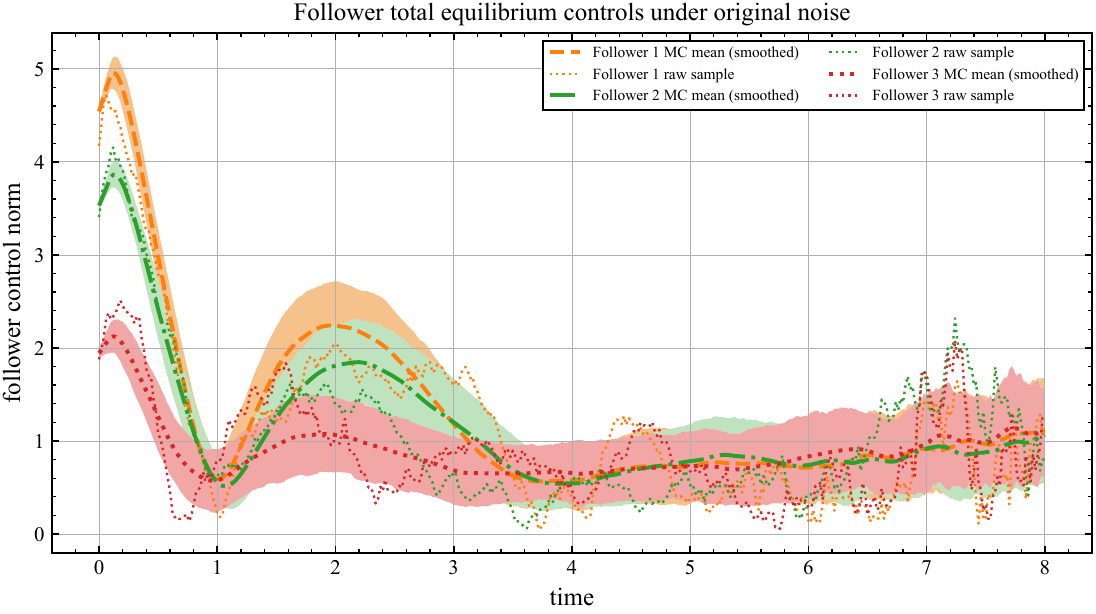}
    \caption{Follower total equilibrium control norms under the original-noise case.}
    \label{fig:follower-controls}
\end{figure}

Figure~\ref{fig:follower-controls} shows the total control norms of the three followers. The follower controls are initially large because the initial configuration is far from the desired V-shape. As the formation is established, the control effort decreases and remains bounded. The differences among the three follower controls are caused by the asymmetric control weights and the different relative positions of the followers in the formation. The raw representative paths display stochastic fluctuations, while the smoothed Monte Carlo means show the average control trend.

\begin{figure}[t]
    \centering
    \includegraphics[width=\columnwidth]{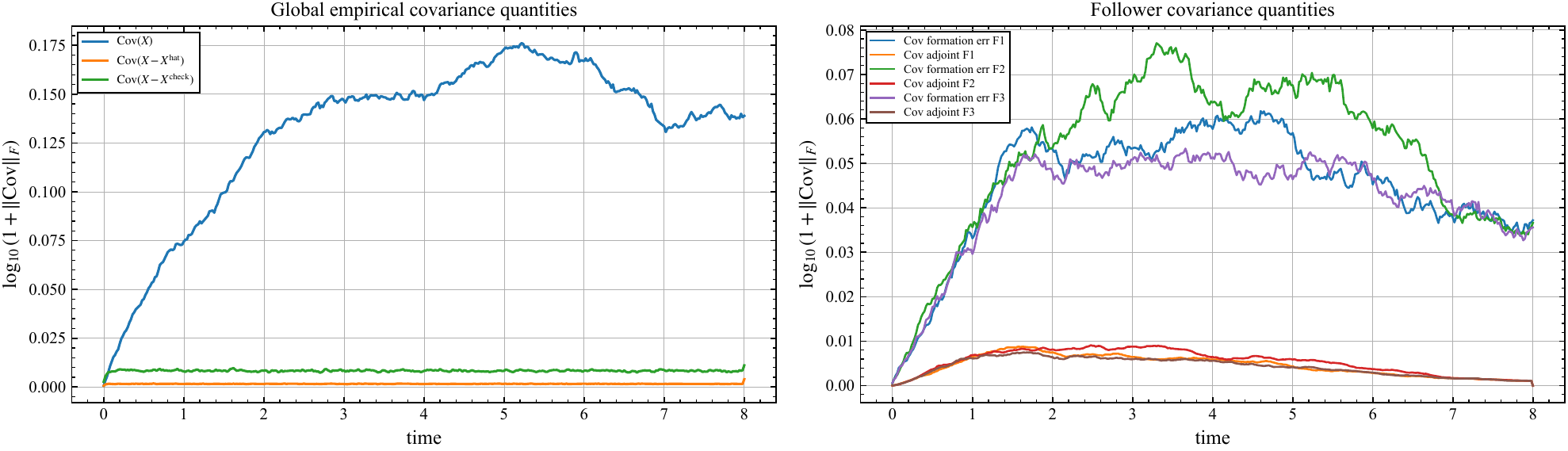}
    \caption{Empirical covariance quantities computed from Monte Carlo sample paths.}
    \label{fig:covariances}
\end{figure}

Figure~\ref{fig:covariances} presents empirical covariance quantities. The global covariance of the state reflects the accumulated effect of stochastic forcing. The covariance of the filtering errors remains much smaller (in which the covariance of the followers' filtering errors is smaller than that of leader because of more information known by the followers) than the state covariance, showing the stabilizing role of the filtering mechanism. The follower covariance quantities describe the stochastic variability of the formation errors and the adjoint-related correction terms. The formation covariance is larger than the adjoint covariance, which is consistent with the fact that the state and formation errors are the dominant stochastic quantities, whereas the adjoint terms act as correction components in the state feedback equilibrium.

\begin{figure}[t]
    \centering
   \includegraphics[width=\columnwidth]{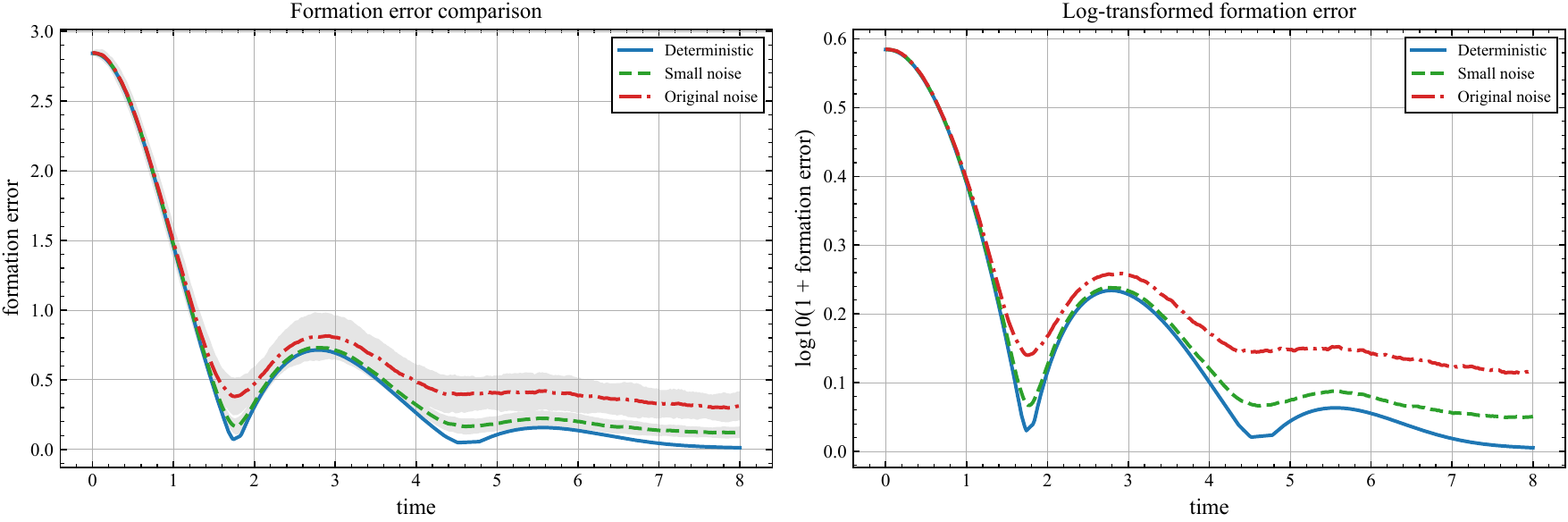}
    \caption{Formation error comparison and log-transformed formation error under different noise levels.}
    \label{fig:formation-decay}
\end{figure}

Figure~\ref{fig:formation-decay} further compares the formation-error decay under deterministic, small-noise, and original-noise cases. The left panel shows the original formation error, and the right panel shows the log-transformed error. The deterministic case converges to the smallest residual error. The small-noise case remains close to the deterministic case, while the original-noise case has a larger residual level due to persistent stochastic excitation. The transient overshoot is expected for second-order robot dynamics, since the system has inertial effects. Overall, the error decay confirms that the proposed state feedback strategy stabilizes the formation in a stochastic neighborhood of the desired V-shape.

The above simulations demonstrate that the proposed partially observed Stackelberg strategy achieves three objectives: formation stabilization, reference tracking, and bounded filtering/control performance under stochastic perturbations. The numerical results are therefore consistent with the theoretical state feedback structure derived in the preceding sections.

\section{Conclusion}\label{s6}
This paper investigates a class of LQ partially observed stochastic Stackelberg differential game, where a single leader and multiple followers are considered. In our setting, the information available to the leader is less than that known by the followers and assume that all the followers acquire the same partially observable information. Motivated by \cite{SX23}, introducing a novel state decomposition and using the orthogonal decomposition technique, we derive the state filtering feedback form of optimal control of the $i$'s follower by the method of completion of square. Next, motivated by the new decoupling technique in \cite{HJX23}, we first generalize and extend the LQ stochastic control of multi-dimensional fully coupled FBSDEs with inhomogeneous terms and then is applied to study the LQ partially observed optimal control problem of FBSDEs, where we relax the restrict condition imposed on control process in \cite{WWX15, ZS22} by the distinguished state decomposition and orthogonal decomposition of forward and backward state processes. A multi-agent formation control is extended to the stochastic framework and studied by the partially observed Stackelberg differential game approach which differs from the deterministic case.

It is worth pointing out that we extend the partially observed optimal control problem in \cite{SX23} into the Stackelberg game with hierarchical framework. Therefore, we aim to obtain the open-loop optimal control of the leader and followers represented by adjoint processes first and then give the closed-loop feedback representations \eqref{bar u1} and \eqref{optimal control of leader} of open-loop optimal control, characterized by filtering processes via decoupling Hamiltonian system. Similar to \cite{SX23}, the discussion in our paper is not considered as the closed-loop controls which should be assumed to be of close-loop form first. Moreover, compared with the discussion in \cite{BCS15}, we mainly solve the open-loop optimal control problems of followers and leader sequentially instead of considering the adapted closed-loop memoryless information structure in \cite{BCS15}. In our future publications, we will focus on the existence of memory effects due to the structure of optimization problems, and hope to report some results relevant to these problems under closed-loop controls and adapted closed-loop memoryless information.

Another consideration is that the leader knows the complete information and the information available to the followers is comparatively less. In this case, the leader need to solve the LQ optimal control problem of fully coupled conditional mean-field FBSDEs, which is not a trivial extension, and it is under investigation in our another working paper. Another interesting extension is to construct the Stackelberg game theory for providing the first continuous-time partially observed leader-follower analysis of practical problem in \cite{LZW20, WZH22} in which the differentially private algorithm is an useful tool to our numerical simulation, as our future research direction.

\appendix

\setcounter{equation}{0}
\numberwithin{equation}{section}
\section{Independence Between Filtering Error process and Filtration}\label{Appendix_A}
In this section, we give the mathematical explanation for the independence between $\tilde{X}$ and $\cF^{Y^{u_1,u_2}_1}_t$, which is the key result for obtaining the filtering equation.

First, we recall that $\tilde{X}_t= X^{u_2}_t-\hat{X}^{u_2}_t$, which is independent of control processes $(u_1,u_2)$. Indeed, by the filtration inclusion \eqref{assumption_filtration}, $\cF_t^{Y_2^{u_1,u_2}}\subset\cF_t^{Y_1^{0,u_1}}$, we have $u_1$ and $u_2$ are both adapted to $\cF_t^{Y^{0,u_1}_1}$, which implies that these forms for both $u_1,u_2$ and their corresponding filtering keep the same. Therefore, the difference $\tilde{X}$ do not depend on control processes. Then for a fixed control $u_1$, we prove that $\tilde{X}_t$ is independent of $\cF^{Y^{u_1,u_2}_1}_t$. For any $u_2$, the followers' control problem can be solved as the partially observed optimal control problem in \cite{SX23}. Therefore, we aim to transfer the original system into the equivalent one without $u_2$. Indeed, we recall the decomposition \eqref{decomposition} and systems composed of state process $X^{1,u_2}$ with control $u_2$
\begin{equation}
\left\{
\begin{aligned}
dX^{1,u_2}_t&=\big[A(t)X^{1,u_2}_t+B_2(t)u_2(t)+\alpha_t\big]dt,\quad X^{1,u_2}_0=0,\\
dY^{1,u_2}_1(t)&=\big[f_1(t)X^{1,u_2}_t+g_1(t)\big]dt,\quad Y^{1,u_2}_1(0)=0,
\end{aligned}
\right.
\end{equation}
and new state process $X^{0,u_1}$ without $u_2$
\begin{equation}\label{X0Y0_u1}
\left\{
\begin{aligned}
dX^{0,u_1}_t&=\big[A(t)X^{0,u_1}_t+B_1^\top(t) u_1(t) \big]dt+C_1(t)dW^1_t+C_2(t)dW^2_t,\quad X^{0,u_1}_0=x_0,\\
dY^{0,u_1}_1(t)&=f_1(t)X^{0,u_1}_tdt+K_1(t)dW^1_t,\quad Y^{0,u_1}_1(0)=0,
\end{aligned}
\right.
\end{equation}
where it should be noted that the controlled (only by $u_1$) signal-observation system is consistent with that studied in \cite{SX23}. Meanwhile, it has been proved that the equivalence of filtrations $\cF^{Y^{u_1,u_2}_1}_t$ and $\cF^{Y^{0,u_1}_1}_t$ in Lemma \ref{equivalent_filtration}. Then by $X^{u_1,u_2}_t=X^{1,u_2}_t+X^{0,u_1}_t$ and $X^{1,u_2}\in\cF^{Y^{0,u_1}_1}_t$, we obtain
\begin{equation}
\begin{aligned}
\hat{X}^{u_1,u_2}_t&=\EE[X^{u_1,u_2}_t|\cF^{Y^{u_1,u_2}_1}_t]=\EE[X^{u_1,u_2}_t|\cF^{Y^{0,u_1}_1}_t]=\EE[X^{1,u_2}_t+X^{0,u_1}_t|\cF^{Y^{0,u_1}_1}_t]\\
&=X^{1,u_2}_t+\EE[X^{0,u_1}_t|\cF^{Y^{0,u_1}_1}_t]=X^{1,u_2}_t+\hat{X}^{0,u_1}_t.
\end{aligned}
\end{equation}
So the filtering error $\tilde{X}_t$ satisfies
\begin{equation}\label{tilde_X_equivalent}
\tilde{X}_t=X^{u_1,u_2}_t-\hat{X}^{u_1,u_2}_t=(X^{1,u_2}_t+X^{0,u_1}_t)-(X^{1,u_2}_t+\hat{X}^{0,u_1}_t)=X^{0,u_1}_t-\hat{X}^{0,u_1}_t,
\end{equation}
which means that the error of original system is equivalent to that of $u_2$-independent system \eqref{X0Y0_u1}. Now we focus on sub-system \eqref{X0Y0_u1}. For a fixed but arbitrary admissible control $u_1$, on one hand, $\hat{X}^{0,u_1}_t=\EE[X^{0,u_1}_t|\cF^{Y^{0,u_1}_1}_t]$ is the orthogonal projection of $X^{0,u_1}_t$ on $L^2\big(\Om,\cF^{Y^{0,u_1}_1}_t\big)$, so
\begin{equation}
\begin{aligned}
\EE\Big[(X^{0,u_1}_t-\hat{X}^{0,u_1}_t)\Big|\cF^{Y^{0,u_1}_1}_t\Big]=0.
\end{aligned}
\end{equation}
Therefore, for any square-integrable $\cF^{Y^{0,u_1}_1}_t$-adapted stochastic process $\kappa\in L^2\big(\cF^{Y^{0,u_1}_1}_t\big)$, we have
\begin{equation}
\begin{aligned}
&\EE\big[(X^{0,u_1}_t-\hat{X}^{0,u_1}_t)\kappa\big]=\EE\Big[\EE\big[(X^{0,u_1}_t-\hat{X}^{0,u_1}_t)\kappa\big|\cF^{Y^{0,u_1}_1}_t\big]\Big]\\
&=\EE\Big[\kappa\EE\big[(X^{0,u_1}_t-\hat{X}^{0,u_1}_t)\big|\cF^{Y^{0,u_1}_1}_t\big]\Big]=0,
\end{aligned}\end{equation}
which implies that
\begin{equation}
(X^{0,u_1}_t-\hat{X}^{0,u_1}_t)\perp L^2\big(\cF^{Y^{0,u_1}_1}_t\big).
\end{equation}
In particular, for any $s\leq t$, $\EE\big[(X^{0,u_1}_t-\hat{X}^{0,u_1}_t)Y^{0,u_1}_1(s)\big]=0$. On the other hand, it is obvious that $X^{0,u_1}$ and $Y^{0,u_1}_1$ are both linear Gaussian functions of $(W^1,W^2)$. Then for any $0\leq s_1<\cdots<s_k\leq t$, $(X^{0,u_1},Y^{0,u_1}_{1}(s_1),\dots,Y^{0,u_1}_{1}(s_k))$ is joint Gaussian vector. So is $(X^{0,u_1}_t-\hat{X}^{0,u_1}_t,Y^{0,u_1}_{1}(s_1),\dots,Y^{0,u_1}_{1}(s_k))$. Moreover, for any finite time $0\leq s_1<\cdots<s_k\leq t$, it implies that $\EE\big[(X^{0,u_1}_t-\hat{X}^{0,u_1}_t)(Y^{0,u_1}_{1}(s_1),\dots,Y^{0,u_1}_{1}(s_k))\big]=0$.

Next, for fixed $t$, take any finite time $0\leq s_1<\cdots<s_k\leq t$, we have $Cov\big((X^{0,u_1}_t-\hat{X}^{0,u_1}_t),(Y^{0,u_1}_{1}(s_1),\dots,Y^{0,u_1}_{1}(s_k))\big)=0$, so we obtain
\begin{equation}
(X^{0,u_1}_t-\hat{X}^{0,u_1}_t)\ \text{is independent of}\ (Y^{0,u_1}_{1}(s_1),\dots,Y^{0,u_1}_{1}(s_k)).
\end{equation}
Continuously, set
\begin{equation}
\cC_t=\Big\{\{(Y^{0,u_1}_{1}(s_1),\dots,Y^{0,u_1}_{1}(s_k))\in\CC\}:\ k\geq1,\ 0\leq s_1<\cdots<s_k\leq t,\ \CC\in\cB(\RR^{l_1k})\Big\}.
\end{equation}
which is a $\pi$-system and generate $\sigma$-field $\sigma(Y^{0,u_1}_1(s):0\leq s\leq t)$ and its completed form $\cF^{Y^{0,u_1}_1}_t$. For any $\mathbb{A}\in\cC_t$ and any Borel set $\BB\subset\RR^n$, we have
\begin{equation}
\PP\big((X^{0,u_1}_t-\hat{X}^{0,u_1}_t)\in\BB,\mathbb{A}\big)=\PP\big((X^{0,u_1}_t-\hat{X}^{0,u_1}_t)\in\BB\big)\PP(\mathbb{A}).
\end{equation}
Then define
\begin{equation}
\La:=\Big\{\mathbb{A}\in\cF^{Y^{0,u_1}_1}_t:\ \PP\big((X^{0,u_1}_t-\hat{X}^{0,u_1}_t)\in\BB,\mathbb{A}\big)=\PP\big((X^{0,u_1}_t-\hat{X}^{0,u_1}_t)\in\BB\big)\PP(\mathbb{A}),\ \forall\ \BB\in\cB(\RR^n)\Big\}.
\end{equation}
which can be verified as a $\la$-system, and $\cC_t\subset\La$. By $\pi-\la$ theorem, $\cF^{Y^{0,u_1}_1}_t=\sigma(\cC_t)\subset\La$. Therefore, for any $\mathbb{A}\in\cF^{Y^{0,u_1}_1}_t$ and any Borel set $\BB\subset\RR^n$, we have
\begin{equation}
\PP\big((X^{0,u_1}_t-\hat{X}^{0,u_1}_t)\in\BB,\mathbb{A}\big)=\PP\big((X^{0,u_1}_t-\hat{X}^{0,u_1}_t)\in\BB\big)\PP(\mathbb{A}),
\end{equation}
which implies that $(X^{0,u_1}_t-\hat{X}^{0,u_1}_t)$ is independent of $\cF^{Y^{0,u_1}_1}_t$. Therefore, by \eqref{tilde_X_equivalent} and Lemma \ref{equivalent_filtration}, it implies $\tilde{X}_t$ is independent of $\cF^{Y^{u_1,u_2}_1}_t$.
\begin{remark}
Note that the control $u_1$ has no influence on the procedure of obtaining the independence since we have fixed a control $u_1$ in advance. Strictly speaking, we claim that for each fixed admissible control $u_1$, $\tilde{X}_t$ is independent of $\cF^{Y^{0,u_1}_1}_t$. Meanwhile, we take a fixed one but arbitrary, so it also hold for the fact that $\tilde{X}_t$ is independent of $\cF^{Y^{0,\bar{u}_1}_1}_t$ after the optimal control $\bar{u}_1$ is obtained finally, which also implies that $\tilde{X}_t$ is independent of $\cF^{Y^{\bar{u}_1,\bar{u}_2}_1}_t$.
\end{remark}

\section{Results of Fully Coupled Forward-Backward Stochastic LQ Optimal Control Problem}\label{FBSLQ}

We consider an LQ optimal control problem of fully coupled forward-backward stochastic differential controlled system, for the multi-dimensional case with inhomogeneous terms, which is an extension of \cite{HJX23}.

We consider the following fully coupled forward-backward stochastic differential controlled system, for the state triple $(X,Y,Z)\equiv\big(X,(Y^1,\cdots,Y^N)^\top,Z^1\equiv(Z^{11},\cdots,Z^{N1})^\top,\\ \cdots, Z^{l_2}\equiv(Z^{1l_2},\cdots,Z^{Nl_2})^\top\big)$:
\begin{equation}\label{A1}
\left\{
\begin{aligned}
dX_t&=\bigg[A_1X+\tilde{B}_1^\top Y+\sum_{k=1}^{l_2}(C^k_1)^\top Z^k+D_1u+E_1\bigg] d t\\
&\quad +\sum_{j=1}^{l_2}\bigg[ A^j_2X+(B_2^j)^\top Y+\sum_{k=1}^{l_2}(C_2^{jk})^\top Z^k+D^j_2u+E^j_2\bigg] d B^j_t,\\
dY^i_t&=-\bigg[A_{3i}X_t+B_{3i}^\top Y^i+\sum_{j=1}^{l_2}(C^j_{3i})^\top Z^{ij}+D_{3i}u+E_{3i}\bigg] d t +\sum_{j=1}^{l_2}Z^{ij}_{t} d B^j_t,\\
X_0&=x_0,\ Y_T^i=F_i X_T+\xi_i,
\end{aligned}
\right.
\end{equation}
where $u\in\mathbb{R}$, $X,Y^i,Z^{ij}\in\mathbb{R}^n$ and $B^j\in\RR$ is classical Brownian motion, $i=1,\cdots,N,j=1,\cdots,l_2$, and cost functional:
\begin{equation}\label{A2}
\begin{aligned}
J(u(\cdot))&=\frac{1}{2} \mathbb{E}\bigg[\int_0^T \bigg(\langle A_4 X,X\rangle+Y^\top B_4 Y+\sum_{j=1}^{l_2}(Z^j)^\top C^j_4 Z^j+\langle D_4 u,u\rangle \bigg)d t\\
&\qquad\qquad +\langle G X_T,X_T\rangle+Y_0^\top H Y_0\bigg].
\end{aligned}
\end{equation}
In the above, the coefficients are matrix-valued functions of proper dimensions. Thus, by the maximum principle, the optimal control $\bar{u}$ satisfies
\begin{equation*}
 \sum_{i=1}^ND^\top_{3i} h^i_t+D_1^\top m_t+\sum_{j=1}^{l_2}(D^j_2)^\top n^j_t+D_4 \bar{u}_t=0,
\end{equation*}
where adjoint process triple $(h,m,n)\equiv\big((h^1,\cdots,h^N)^\top,m,(n^1,\cdots,n^{l_2})^\top\big)$ for $h^i,m,n^j\in\mathbb{R}^n$, $i=1,\cdots,N,j=1,\cdots,l_2$, satisfies:
\begin{equation*}
\left\{
\begin{aligned}
dh^i_t&=\bigg[B_{3i} h^i+\tilde{B}_{1i} m+\sum_{j=1}^{l_2}B^j_{2i} n^j+B_{4i} \bar{Y}^i\bigg] d t\\
&\quad +\sum_{j=1}^{l_2}\bigg[C^j_{3i} h^i+C^j_{1i} m+\sum_{k=1}^{l_2}C^{jk}_{2i} n^k+\sum_{k=1}^{l_2}C^{jk}_{4i} \bar{Z}^{ik}\bigg] d B^j_{t}, \\
dm_t&=-\bigg[\sum_{i=1}^NA_{3i} h^i+A_1^\top m+\sum_{j=1}^{l_2}(A^j_2)^\top n^j+A_4 \bar{X}\bigg] d t +\sum_{j=1}^{l_2}n^j_td B^j_t, \\
h^i_0&=H_i \bar{Y}^i_0,\ m_{T}=G\bar{X}_T+\sum_{i=1}^NF^\top_i h^i_T,\quad i=1,\cdots,N,
\end{aligned}
\right.
\end{equation*}
with the optimal state triple $(\bar{X},\bar{Y}^i,\bar{Z}^{ik})$.

\subsection{Definite case}

We suppose that $D_4>0$ in this subsection. Thus the optimal control is:
\begin{equation*}
 \bar{u}_t=-D_4^{-1}D^\top_3h_t-D_4^{-1}D^\top_1m_t-D_4^{-1}\sum_{j=1}^{l_2}(D^j_2)^\top n^j_t,
\end{equation*}
where $D_3^\top h\equiv\sum\limits_{i=1}^ND^\top_{3i}h^i$. Insert $\bar{u}$ back then we have the following matrix-valued FBSDEs:
\begin{equation*}
\left\{
\begin{aligned}
d \bar{X}_t&=\bigg[A_{1} \bar{X}+\tilde{B}_1^\top \bar{Y}+\sum_{k=1}^{l_2}(C^k_1)^\top \bar{Z}^k- D_1D_4^{-1}D_3^\top h-D_1D_4^{-1} D^\top_1 m\\
&\qquad -D_1D_4^{-1} \sum_{i=1}^{l_2}(D^i_2)^\top n^i_{t}+E_1\bigg] d t +\sum_{j=1}^{l_2}\bigg[A^j_2 \bar{X}+(B_2^j)^\top\bar{Y}+\sum_{k=1}^{l_2}(C_2^{jk})^\top \bar{Z}^k\\
&\qquad -D^j_2D_4^{-1} D_3^\top h-D^j_2D_4^{-1} D_1^\top m-D^j_2D_4^{-1} \sum_{j=1}^{l_2}(D^j_2)^\top n^j+E^j_2\bigg] dB^j_t, \\
d \bar{Y}_{t}&=-\bigg[A_3\bar{X}+B_3^\top \bar{Y}+\sum_{j=1}^{l_2}(C^j_3)^\top \bar{Z}^j-D_3D_4^{-1} D_3^\top h-D_3D_4^{-1} D_1^\top m\\
&\qquad -\sum_{j=1}^{l_2}D_3D_4^{-1} (D^j_2)^\top n^j+E_3\bigg] d t +\sum_{j=1}^{l_2}\bar{Z}^jdB^j_t, \\
d h_{t}&=\bigg[B_3 h+\tilde{B}_1 m+\sum_{j=1}^{l_2}B^j_2 n^j+B_4 \bar{Y}\bigg] dt\\
&\quad +\sum_{j=1}^{l_2}\bigg[C^j_3 h+C^j_1 m+\sum_{k=1}^{l_2}C^{jk}_2 n^k+\sum_{k=1}^{l_2}C^{jk}_4 \bar{Z}^k\bigg] d B^j_t,\\
dm_{t}&=-\bigg[A_3^\top h+A_1^\top m+\sum_{j=1}^{l_2}(A^j_2)^\top n^j+A_4 \bar{X}\bigg]dt +\sum_{j=1}^{l_2}n^j_td B^j_{t}, \\
\bar{X}_0&=x_0,\ \bar{Y}_T=F \bar{X}_{T}+\xi,\ h_0=H \bar{Y}_0,\ m_T=G\bar{X}_T+F^\top h_T,
\end{aligned}
\right.
\end{equation*}
where
\begin{equation*}
\begin{aligned}
C^j_{4i}&\triangleq(C_{4i}^{j1},C_{4i}^{j2},\cdots,C_{4i}^{jl_2}),\ i=1,\cdots,N,\ j=1,\cdots,l_2,\\
A_3&\triangleq\left(\begin{array}{c}
 A_{31}\\
 A_{32}\\
 \vdots\\
 A_{3N}
 \end{array}\right),\quad \tilde{B}_1\triangleq\left(\begin{array}{c}
\tilde{B}_{11}\\
\tilde{B}_{12}\\
\vdots\\
\tilde{B}_{1N}
\end{array}\right),\quad
 B_3\triangleq\begin{pmatrix}B_{31}& & & \\ &B_{32}& & \\ & &\ddots& \\ & & &B_{3N}\end{pmatrix},\\
C^j_3&\triangleq\begin{pmatrix}C^j_{31}& & & \\ &C^j_{32}& & \\ & &\ddots& \\ & & &C^j_{3N}\end{pmatrix},\ j=1,\dots,l_2,\quad
 B^j_2\triangleq\left(\begin{array}{c}
B^j_{21}\\
B^j_{22}\\
\vdots\\
B^j_{2N}
\end{array}\right),\\
D_3&(D_4)^{-1}D_3^\top\triangleq
 \begin{pmatrix}D_{31}D_4^{-1}D_{31}^\top&D_{31}D_4^{-1}D_{32}^\top &\cdots&D_{31}D_4^{-1}D^\top_{3N} \\D_{32}D_4^{-1}D^\top_{31} &D_{32}D_4^{-1}D^\top_{32}&\cdots &D_{32}D_4^{-1}D^\top_{3N} \\
 \vdots &\vdots &\ddots&\vdots \\D_{3N}D_4^{-1}D_{31}^\top &D_{3N}D_4^{-1}D^\top_{32} &\cdots &D_{3N}D_4^{-1}D^\top_{3N}\end{pmatrix},\\
\end{aligned}
\end{equation*}
and
\begin{equation*}
\begin{aligned}
D_3&(D_4)^{-1}D_2^\top\triangleq\begin{pmatrix}D_{31}D_4^{-1}(D^1_2)^\top&D_{31}D_4^{-1}(D^2_2)^\top &\cdots&D_{31}D_4^{-1}(D^{l_2}_2)^\top \\D_{32}D_4^{-1}(D^1_2)^\top &D_{32}D_4^{-1}(D^2_2)^\top&\cdots &D_{32}D_4^{-1}(D^{l_2}_2)^\top \\\vdots &\vdots &\ddots&\vdots \\D_{3N}D_4^{-1}(D^{1}_2)^\top &D_{3N}D_4^{-1}(D^{2}_2)^\top &\cdots &D_{3N}D_4^{-1}(D^{l_2}_2)^\top\end{pmatrix},\\
B_2&\triangleq\begin{pmatrix}B_{21}^1&B_{21}^2 &\cdots&B_{21}^{l_2} \\B_{22}^1&B_{22}^2 &\cdots&B_{22}^{l_2} \\\vdots &\vdots &\ddots&\vdots \\B_{2N}^1&B_{2N}^2 &\cdots&B_{2N}^{l_2}\end{pmatrix},\quad
B_4\triangleq\begin{pmatrix}B_{41}& & & \\ &B_{42}& & \\ & &\ddots& \\ & & &B_{4N}\end{pmatrix},\\
C^j_4&\triangleq\begin{pmatrix}C^j_{41}& & & \\ &C^j_{42}& & \\ & &\ddots& \\ & & &C^j_{4N}\end{pmatrix},\quad
C^{jk}_2=\left(\begin{array}{l}
C^{jk}_{21}\\
C^{jk}_{22}\\
\vdots\\
C^{jk}_{2N}
\end{array}\right),\quad
C^j_2=\begin{pmatrix}C_{21}^{j1}&C_{21}^{j2} &\cdots&C_{21}^{jl_2} \\C_{22}^{j1}&C_{22}^{j2} &\cdots&C_{22}^{jl_2} \\\vdots &\vdots &\ddots&\vdots \\C_{2N}^{j1}&C_{2N}^{j2} &\cdots&C_{2N}^{jl_2}\end{pmatrix},\\
C^j_1&\triangleq\left(\begin{array}{c}
C^j_{11}\\
C^j_{12}\\
\vdots\\
C^j_{1N}
\end{array}\right),\ j=1,\dots,l_2,\quad
D_3=\left(\begin{array}{c}
D_{31}\\
D_{32}\\
\vdots\\
D_{3N}
\end{array}\right),\quad E_3=\left(\begin{array}{c}
E_{31}\\
E_{32}\\
\vdots\\
E_{3N}
\end{array}\right),\\
H&\triangleq\begin{pmatrix}H_1& & & \\ &H_2& & \\ & &\ddots& \\ & & &H_N\end{pmatrix},\quad
F=(F_1,F_2,\cdots,F_N)^\top,\quad \xi=(\xi_1,\xi_2,\cdots,\xi_N)^\top.
\end{aligned}
\end{equation*}

Set
\begin{equation*}
\begin{aligned}
&\tilde{X}\equiv\left(\begin{array}{l}
\bar{X} \\
h
\end{array}\right),\quad \tilde{Y}\equiv\left(\begin{array}{l}
m \\
\bar{Y}
\end{array}\right),\quad \tilde{Z}^j\equiv\left(\begin{array}{l}
n^j \\
\bar{Z}^j
\end{array}\right)\in\RR^{(N+1)n},\ j=1,\dots,l_2, \\
&\tilde{A}_1\triangleq\left(\begin{array}{cc}
A_1 & -D_1D_4^{-1} D_3^\top \\
0 & B_3
\end{array}\right), \quad \tilde{\tilde{B}}_1\triangleq\left(\begin{array}{cc}
-D_1D_4^{-1} D_{1}^\top & \tilde{B}_1^\top \\
\tilde{B}_1 & B_4
\end{array}\right), \\
&\widetilde{C}^j_1\triangleq\left(\begin{array}{cc}
-D_1D_4^{-1} (D^j_2)^\top & (C^j_1)^\top \\
B^j_2 & 0
\end{array}\right),\quad
\tilde{A}^j_2\triangleq\left(\begin{array}{cc}
A^j_2 & -D_2^jD_4^{-1}  D_3^\top \\
0 & C^j_3
\end{array}\right),\\
&\tilde{B}^j_2\triangleq\left(\begin{array}{cc}
-D_2^jD_4^{-1} D_1^\top & (B^j_2)^\top \\
C^j_1 & 0
\end{array}\right),\quad \tilde{C}^{jk}_2\triangleq\left(\begin{array}{cc}
-D_2^jD_4^{-1} (D^k_2)^\top & (C^{jk}_2)^\top \\
C^{jk}_2 & C^{jk}_4
\end{array}\right),\\
&\tilde{A}_3\triangleq\left(\begin{array}{cc}
A_4 & A_3^\top \\
A_3 & -D_3D_4^{-1} D_3^\top
\end{array}\right),\quad \tilde{B}_3\triangleq\left(\begin{array}{cc}
A_1^\top & 0\\
-D_3D_4^{-1}D_1^\top & B_3^\top
\end{array}\right),\\
&\tilde{C}^j_3\triangleq\left(\begin{array}{cc}
(A^j_2)^\top & 0\\
-D_3D_4^{-1}(D^j_2)^\top & (C^j_3)^\top
\end{array}\right), \quad
\tilde{X}_0\triangleq\left(\begin{array}{c}
x_0 \\
H\bar{Y}_0
\end{array}\right),\quad \tilde{E}_1\triangleq\left(\begin{array}{c}
E_1 \\
0
\end{array}\right),\\
& \tilde{E}^j_2\triangleq\left(\begin{array}{c}
E^j_2 \\
0
\end{array}\right),\quad \tilde{E}_3\triangleq\left(\begin{array}{c}
0\\
E_3
\end{array}\right),\quad \tilde{F}\triangleq\left(\begin{array}{cc}
G & F^\top\\
F & 0
\end{array}\right),\quad \tilde{\xi}\triangleq\left(\begin{array}{c}
0\\
\xi
\end{array}\right),
\end{aligned}
\end{equation*}
we have
\begin{equation}\label{FBSDE1}
\left\{
\begin{aligned}
d \tilde{X}_t&=\bigg[\tilde{A}_1\tilde{X}+\tilde{\tilde{B}}_1\tilde{Y}+\sum_{j=1}^{l_2}\tilde{C}^j_1\tilde{Z}^j+\tilde{E}_1\bigg]dt
 +\sum_{j=1}^{l_2}\bigg[\tilde{A}^j_2\tilde{X}+\tilde{B}^j_2\tilde{Y}+\sum_{k=1}^{l_2}\tilde{C}^{jk}_2\tilde{Z}^k+\tilde{E}^j_2\bigg] d B^j_{t},\\
d \tilde{Y}_t&=-\bigg[\tilde{A}_3\tilde{X}+\tilde{B}_3\tilde{Y}+\sum_{j=1}^{l_2}\tilde{C}^j_3\tilde{Z}^j+\tilde{E}_3\bigg] dt+\sum_{j=1}^{l_2}\tilde{Z}^j_t d B^j_t,\\
\tilde{X}_0&=\left(\begin{array}{c}
x_0 \\
H\tilde{Y}_0
\end{array}\right),\ \tilde{Y}_T=\tilde{F}\tilde{X}_T+\tilde{\xi}.
\end{aligned}
\right.
\end{equation}

In the following, we only study the case of
\begin{equation*}
\begin{aligned}
C_2^{jk}=
\begin{cases}
0,&{1\leq j\neq k\leq l_2},\\
C_2^{jk},&{j=k},
\end{cases}
\end{aligned}
\end{equation*}
\begin{equation*}
\begin{aligned}
C_4^{jk}=
\begin{cases}
0,&{1\leq j\neq k\leq l_2},\\
C_4^{jk},&{j=k},
\end{cases}
\end{aligned}
\end{equation*}
and
\begin{equation*}
\begin{aligned}
D_4^{jk}=
\begin{cases}
0,&{1\leq j\neq k\leq m},\\
D_4^{jk},&{j=k},
\end{cases}
\end{aligned}
\end{equation*}
for simplification, and the general case has no intrinsic difficulty but complexity of computation and notation. Now \eqref{FBSDE1} becomes
\begin{equation}\label{FBSDE2}
\left\{
\begin{aligned}
d \tilde{X}_t&=\bigg[\tilde{A}_1\tilde{X}+\tilde{\tilde{B}}_1\tilde{Y}+\sum_{j=1}^{l_2}\tilde{C}^j_1\tilde{Z}^j+\tilde{E}_1\bigg]dt
+\sum_{j=1}^{l_2}\bigg[\tilde{A}^j_2\tilde{X}+\tilde{B}^j_2\tilde{Y}+\tilde{C}^{jj}_2\tilde{Z}^j+\tilde{E}^j_2\bigg] d B^j_t,\\
d \tilde{Y}_t&=-\bigg[\tilde{A}_3\tilde{X}+\tilde{B}_3\tilde{Y}+\sum_{j=1}^{l_2}\tilde{C}^j_3\tilde{Z}^j+\tilde{E}_3\bigg] dt+\sum_{j=1}^{l_2}\tilde{Z}^j_t d B^j_t,\\
\tilde{X}_0&=\left(\begin{array}{c}
x_0 \\
H\tilde{Y}_0
\end{array}\right),\ \tilde{Y}_T=\tilde{F}\tilde{X}_T+\tilde{\xi}.
\end{aligned}
\right.
\end{equation}
Similar to \cite{HJX23}, we can build the relation $\tilde{Y}_t=Q_{t} \tilde{X}_t+\vp_t$ with
\begin{equation*}
\left\{
\begin{aligned}
&\dot{Q}+Q \tilde{A}_1+\tilde{A}_1^\top Q+\sum_{j=1}^{l_2}\big(\tilde{B}^j_2Q+\tilde{A}^j_2\big)^\top\big(I-Q\tilde{C}^{jj}_2\big)^{-1}Q\big(\tilde{B}^j_2Q+\tilde{A}^j_2\big)+Q\tilde{\tilde{B}}_1Q+\tilde{A}_3=0,\\
&Q_T=\tilde{F},
\end{aligned}
\right.
\end{equation*}
and
\begin{equation*}
\left\{\begin{aligned}
d \varphi_t&=-\bigg[\bigg(Q \tilde{\tilde{B}}_1+\tilde{B}_3+\sum_{j=1}^{l_2}\big(Q \tilde{C}^j_1+\tilde{C}^j_3\big)\big(I-Q \tilde{C}^{jj}_2\big)^{-1} Q \tilde{B}^j_2\bigg) \varphi\\
&\qquad +\sum_{j=1}^{l_2}\big(Q \tilde{C}^j_1+\tilde{C}^j_3\big)\big(I-Q \tilde{C}^{jj}_2\big)^{-1} V^j\\
&\qquad +\sum_{j=1}^{l_2}\big(Q\tilde{C}^j_1+\tilde{C}^j_3\big)\big(I-Q \tilde{C}^{jj}_2\big)^{-1} Q \tilde{E}^j_2+Q \tilde{E}_1+\tilde{E}_3\bigg]dt+\sum_{j=1}^{l_2}V^j_t d B^j_t, \\
\vp_T&=\tilde{\xi},
\end{aligned}\right.
\end{equation*}
where, noting that $\tilde{B}_3\equiv\tilde{A}^\top_1,\tilde{C}^j_1\equiv(\tilde{B}^j_2)^\top,\tilde{C}^j_3\equiv(\tilde{A}^j_2)^\top$. Set
\begin{equation*}
\begin{aligned}
& Q\triangleq\left(\begin{array}{cc}
Q_1 & Q_2 \\
Q_3 & -Q_4
\end{array}\right),\quad k^j\triangleq\big(I-Q\tilde{C}^{jj}_2\big)^{-1}Q\big(\tilde{B}^j_2Q+\tilde{A}^j_2\big)\equiv\left(\begin{array}{ll}
k^j_1 & k^j_2 \\
k^j_3 & k^j_4
\end{array}\right),\\
&\varphi=\left(\begin{array}{l}
\varphi_1 \\
\varphi_2
\end{array}\right),\quad
V^j\triangleq\left(\begin{array}{l}
V^j_1 \\
V^j_2
\end{array}\right),
\quad\left(\begin{array}{ll}
J^j_1 & J^j_2 \\
J^j_3 & J^j_4
\end{array}\right)\triangleq\big(I-Q \tilde{C}^{jj}_2\big)^{-1} Q \tilde{B}^j_2, \\& \left(\begin{array}{ll}
I^j_1 & I^j_2 \\
I^j_3 & I^j_4
\end{array}\right)\triangleq\big(I-Q \tilde{C}^{jj}_{2}\big)^{-1},
\end{aligned}
\end{equation*}
where $Q_1,k^j_1,J^j_1,I^j_1\in\RR^{n\times n}$, $Q_2,k^j_2,J^j_2,I^j_2\in\RR^{n\times Nn}$, $Q_3,k^j_3,J^j_3,I^j_3\in\RR^{Nn\times n}$, $Q_4,k^j_4,J^j_4,I^j_4\in\RR^{Nn\times Nn}$, $\vp_1,V^j_1\in\RR^n$ and $\vp_2,V^j_2\in\RR^{Nn}$.

Moreover, the following relationships hold:
\begin{equation}\label{defrela}
\begin{aligned}
& \bar{X}=X^*,\quad h=h^*,\quad \bar{Y}=Q_3 X^*-Q_4 h^*+\varphi_2,\quad m=Q_1 X^*+Q_2 h^*+\varphi_1, \\
& \bar{Z}^j=k^j_3 X^*+k^j_4 h^*+J^j_3 \varphi_1+J^j_4 \varphi_2+(I^j_3 Q_1+I^j_4 Q_3) E^j_2+I^j_3 V^j_1+I^j_4 V^j_2, \\
& n^j=k^j_1 X^*+k^j_2 h^*+J^j_1 \varphi_1+J^j_2 \varphi_2+(I^j_1 Q_1+I^j_2 Q_3) E^j_2+I^j_1 V^j_1+I^j_2 V^j_2,
\end{aligned}
\end{equation}
and the optimal control is
\begin{equation*}
\begin{aligned}
\bar{u}&=D_4^{-1}\bigg(-D^\top_1Q_1-\sum_{j=1}^{l_2}(D^j_2)^\top k^j_1\bigg) X^*+D_4^{-1}\bigg(-D_3^\top-D_1^\top Q_2-\sum_{j=1}^{l_2}(D^j_2)^\top k^j_2\bigg) h^* \\
&\quad -D_4^{-1}\bigg(D^\top_1 \varphi_1+\sum_{j=1}^{l_2}(D^j_2)^\top J^j_1 \varphi_1+\sum_{j=1}^{l_2}(D^j_2)^\top J^j_2\varphi_2\bigg)\\
&\quad -D_4^{-1}\bigg(\sum_{j=1}^{l_2}(D^j_2)^\top I^j_1 V^j_1+\sum_{j=1}^{l_2}(D^j_2)^\top I^j_2 V^j_2\bigg)-D_4^{-1} \sum_{j=1}^{l_2}(D^j_2)^\top\big(I^j_1 Q_1+I^j_2 Q_3\big)E^j_2,
\end{aligned}
\end{equation*}
where
\begin{equation*}
\begin{aligned}
d\tilde{X}^*_t&=\bigg[\bigg(\tilde{A}_1+\tilde{\tilde{B}}_1 Q+\sum_{j=1}^{l_2}\tilde{C}^j_1 k^j\bigg) \tilde{X}^* +\bigg(\tilde{\tilde{B}}_1+\sum_{j=1}^{l_2}\tilde{C}^j_1\big(I-Q\tilde{C}^{jj}_2\big)^{-1} Q \widetilde{B}^j_2\bigg) \varphi\\
&\qquad +\sum_{j=1}^{l_2}\tilde{C}^j_1\big(I-Q\tilde{C}^{jj}_2\big)^{-1}V^j+\tilde{E}_1+\sum_{j=1}^{l_2}\tilde{C}^j_1\big(I-Q \tilde{C}^{jj}_2\big)^{-1} Q \tilde{E}^j_2\bigg]dt \\
&\quad +\sum_{j=1}^{l_2}\bigg[\Big(\tilde{A}^j_2+\tilde{B}^j_2 Q+\tilde{C}^{jj}_2 k^j\Big) \tilde{X}^* +\Big(\tilde{B}^j_2+\tilde{C}^{jj}_2\big(I-Q \tilde{C}^{jj}_2\big)^{-1} Q \tilde{B}^j_2\Big)\vp\\
&\qquad +\tilde{C}^{jj}_2\big(I-Q\tilde{C}^{jj}_2\big)^{-1}V^j+\tilde{E}^j_2 +\tilde{C}^{jj}_2\big(I-Q \tilde{C}^{jj}_2\big)^{-1} Q \tilde{E}^j_2\bigg]d B^j_t.
\end{aligned}
\end{equation*}

\subsection{Indefinite case}

In this subsection, we do not need the condition $D_{4}>0$. Therefore, in order to obtain the optimal control, motivated by relations \eqref{defrela}, we set
\begin{equation}\label{indefrela}
\begin{aligned}
m_t&=P_1(t) \bar{X}_t+P_2(t) h_t+\vp_1(t), \\
\bar{Y}_t&=P_2^\top(t) \bar{X}_t-P_3(t) h_t+\vp_2(t),
\end{aligned}
\end{equation}
and $\left(\varphi_1, V_1\right),\left(\varphi_2, V_2\right)$ satisfy
\begin{equation*}
d \varphi_1=-\gamma_1 dt+\sum_{j=1}^{l_2}V^j_1 d B^j_t, \quad d\varphi_2=-\gamma_2 dt+\sum_{j=1}^{l_2}V^j_2 d B^j_t.
\end{equation*}
By applying It\^o's formula, we have
\begin{equation*}
\begin{aligned}
& d \bar{Y}_t=P_2^\top\bigg[A_1 \bar{X}+\tilde{B}_1^\top \bar{Y}+\sum_{k=1}^{l_2}(C^k_1)^\top \bar{Z}^k+D_1 \bar{u}+E_1\bigg] d t\\
&\qquad +P_2^\top\sum_{j=1}^{l_2}\big[A^j_2 \bar{X}+(B^j_2)^\top \bar{Y}+(C^{jj}_2)^\top \bar{Z}^j+D^j_2 \bar{u}+E^j_2\big] d B^j_t\\
&\qquad +\dot{P}_2^\top \bar{X} d t-P_3\bigg[B_3 h+\tilde{B}_1 m+\sum_{j=1}^{l_2}B^j_2 n^j+B_4 \bar{Y}\bigg] d t\\
&\qquad -P_3\sum_{j=1}^{l_2}\big[C^j_3 h+C^j_1 m+C^{jj}_2 n^j+C^{jj}_4 \bar{Z}^j\big] d B^j_t -\dot{P}_3 h d t-\gamma_2 d t+\sum_{j=1}^{l_2}V^j_2 d B^j_t, \\
\end{aligned}
\end{equation*}\begin{equation*}
\begin{aligned}
& d m_t=P_1\bigg[A_1 \bar{X}+\tilde{B}_1^\top \bar{Y}+\sum_{k=1}^{l_2}(C^k_1)^\top \bar{Z}^k+D_1 \bar{u}+E_1\bigg] d t\\
&\qquad +P_1\sum_{j=1}^{l_2}\big[A^j_2 \bar{X}+(B^j_2)^\top \bar{Y}+(C^{jj}_2)^\top \bar{Z}^j+D^j_2 \bar{u}+E^j_2\big] d B^j_t \\
&\qquad +\dot{P}_1 \bar{X} d t+P_2\bigg[B_3 h+\tilde{B}_1 m+\sum_{j=1}^{l_2}B^j_2 n^j+B_4 \bar{Y}\bigg] d t\\
&\qquad +P_2\sum_{j=1}^{l_2}\big[C^j_3 h+C^j_1 m+C^{jj}_2 n^j+C^{jj}_4 \bar{Z}^j\big] d B^j_t +\dot{P}_2 h d t-\ga_1 d t+\sum_{j=1}^{l_2}V^j_1 d B^j_t.
\end{aligned}
\end{equation*}
Comparing the coefficients, we have
\begin{equation*}
\begin{aligned}
& P_2^\top A_1 \bar{X}+P_2^\top \tilde{B}_1^\top \bar{Y}+P_2^\top\sum_{k=1}^{l_2}(C^k_1)^\top\bar{Z}^k+P_2^\top D_1 \bar{u}+P_2^\top E_1+\dot{P}_2^\top \bar{X}-P_3 B_3 h-P_3 \tilde{B}_1 m\\
&\quad -P_3 \sum_{j=1}^{l_2}B^j_2 n^j -P_3 B_4 \bar{Y} -\dot{P}_3 h-\gamma_2= -\bigg[A_3 \bar{X}+B_3^\top \bar{Y}+\sum_{j=1}^{l_2}(C^j_3)^\top \bar{Z}^j+D_3 \bar{u}+E_3\bigg],\\
& P_1 A_1 \bar{X}+P_1 \tilde{B}_1^\top \bar{Y}+P_1\sum_{k=1}^{l_2} (C^k_1)^\top \bar{Z}^k+P_1 D_1 \bar{u}+P_1 E_1+\dot{P}_1 \bar{X}+P_2 B_3 h+P_2 \tilde{B}_1 m\\
&\quad +P_2 \sum_{j=1}^{l_2}B^j_2 n^j+P_2 B_4 \bar{Y} +\dot{P}_2 h-\gamma_1= -\bigg[A_3^\top h+A^\top_1 m+\sum_{j=1}^{l_2}(A^j_2)^\top n^j+A_4 \bar{X}\bigg],
\end{aligned}
\end{equation*}
and
\begin{equation*}
\begin{aligned}
& \bar{Z}^j=P_2^\top A^j_2 \bar{X}+P_2^\top (B^j_2)^\top \bar{Y}+P_2^\top (C^{jj}_2)^\top \bar{Z}^j+P_2^\top D^j_2 \bar{u}+P_2^\top E^j_2-P_3 C^j_3 h\\
&\qquad -P_3 C^j_1 m-P_3 C^{jj}_2 n^j-P_3 C^{jj}_4 \bar{Z}^j+V^j_2, \\
& n^j_t=P_1 A^j_2 \bar{X}+P_1 (B^j_2)^\top \bar{Y}+P_1 (C^{jj}_2)^\top \bar{Z}^j+P_1 D^j_2 \bar{u}+P_1 E^j_2+P_2 C^j_3 h\\
&\qquad +P_2 C^j_1 m+P_2 C^{jj}_2n^j+P_2 C^{jj}_4 \bar{Z}^j+V^j_1,
\end{aligned}
\end{equation*}
which imply that\vspace{-2mm}
\begin{equation}\label{barZ}
\begin{aligned}
& \bar{Z}^j=(L^j_1)^{-1}\big[P_2^\top A^j_2 \bar{X}+P_2^\top (B^j_2)^\top P_2^\top \bar{X}-P_2^\top (B^j_2)^\top P_3 h+P_2^\top (B^j_2)^\top\varphi_2+P_2^\top D^j_2 \bar{u}\\
&\qquad +P_2^\top E^j_2-P_3 C^j_3 h-P_3 C^j_1 P_1 \bar{X}-P_3 C^j_1 P_2 h-P_3 C^j_1 \vp_1-P_3 C^{jj}_2 n^j+V^j_2\big],
\end{aligned}
\end{equation}
where $L^j_1\triangleq I-P_2^\top (C^{jj}_2)^\top+P_3 C^{jj}_4\in\RR^{Nn\times Nn}$,
\begin{equation}\label{n}\vspace{-2mm}
\begin{aligned}
n^j&= \big[P_1 A^j_2+P_1 (B^j_2)^\top P_2^\top+P_2C^j_1P_1\big] \bar{X}-P_1 (B^j_2)^\top P_3 h+P_1 (B^j_2)^\top \vp_2\\
&\quad +\big[P_1 (C^{jj}_2)^\top+P_2 C^{jj}_4\big] (L^j_1)^{-1}\big[P_2^\top A^j_2 \bar{X}+P_2^\top (B^j_2)^\top P_2^\top \bar{X}\\
&\quad -P_2^\top (B^j_2)^\top P_3 h+P_2^\top (B^j_2)^\top \vp_2+P_2^\top D^j_2 \bar{u}+P_2^\top E^j_2-P_3 C^j_3 h\\
&\quad -P_3 C^j_1 P_1 \bar{X}-P_3 C^j_1 P_2 h-P_3 C^j_1 \varphi_1-P_3 C^{jj}_2 n^j+V^j_2\big]\\
&\quad + P_1 D^j_2 \bar{u}+P_1 E^j_2+P_2 C^j_3 h+P_2 C^j_1P_2 h+P_2 C^j_1 \varphi_1+P_2 C^{jj}_2 n^j+V^j_1,
\end{aligned}
\end{equation}
and\vspace{-2mm}
\begin{equation*}\vspace{-2mm}
\begin{aligned}
&\big[I+\big(P_1 (C^{jj}_2)^\top+P_2 C^{jj}_4\big) (L^j_1)^{-1} P_3 C^{jj}_2-P_2 C^{jj}_2\big] n^j\\
&=\big[P_1 A^j_2+P_1 (B^j_2)^\top P_2^\top+P_2 C^j_1 P_1+\big(P_1 (C^{jj}_2)^\top+P_2 C^{jj}_4\big) (L^j_1)^{-1}\big(P_2^\top A^j_2\\
&\quad +P_2^\top (B^j_2)^\top P^\top_2-P_3 C^j_1 P_1\big)\big] \bar{X}+\big[-P_1 (B^j_2)^\top P_3+\big(P_1 (C^{jj}_2)^\top+P_2 C^{jj}_4\big) \\
&\quad \times(L^j_1)^{-1}\big(-P_2^\top (B^j_2)^\top P_3-P_3 C^j_3-P_3 C^j_1 P_2\big)+P_2 C^j_3+P_2C^j_1P_2\big] h\\
&\quad +\big[P_1 (B^j_2)^\top+\big(P_1 (C^{jj}_2)^\top+P_2 C^{jj}_4\big) (L^j_1)^{-1} P_2^\top (B^j_2)^\top\big] \varphi_{2}\\
&\quad +\big[-\big(P_1 (C^{jj}_2)^\top+P_2 C^{jj}_4\big) (L^j_1)^{-1} P_3 C^j_1+P_2 C^j_1\big] \varphi_1\\
&\quad +\big[\big(P_1 (C^{jj}_2)^\top+P_2 C^{jj}_4\big) (L^j_1)^{-1} P_2^\top D^j_2+P_1 D^j_2\big] \bar{u} +V^j_1\\
&\quad +\big(P_1 (C^{jj}_2)^\top+P_2 C^{jj}_4\big) (L^j_1)^{-1} V^j_2+\big(P_1 (C^{jj}_2)^\top+P_2 C^{jj}_4) (L^j_1)^{-1} P_2^\top E^j_2+P_1 E^j_2.
\end{aligned}
\end{equation*}
Denote
\begin{equation}\label{L1}
\begin{aligned}
L^j_2&\triangleq I+\big(P_1 (C^{jj}_2)^\top+P_2 C^{jj}_4\big) (L^j_1)^{-1} P_3 C^{jj}_2-P_2 C^{jj}_2\in\RR^{n\times n},\\
L^j_3&\triangleq P_1 A^j_2+P_1 (B^j_2)^\top P_2^\top+P_2 C^j_1 P_1+\big(P_1 (C^{jj}_2)^\top+P_2 C^{jj}_4\big) \big(L^j_1)^{-1}(P_2^\top A^j_2\\
&\quad +P_2^\top (B^j_2)^\top P_2^\top-P_3 C^j_1 P_1\big)\in\RR^{n\times n}, \\
L^j_4&\triangleq -P_1 (B^j_2)^\top P_3+\big(P_1 (C^{jj}_2)^\top+P_2 C^{jj}_4\big) (L^j_1)^{-1}\big(-P_2^\top (B^j_2)^\top P_3\\
&\quad -P_3 C^j_3-P_3 C^j_1 P_2\big)+P_2 C^j_3+P_2C^j_1P_2\in\RR^{n\times Nn}, \\
S^j_1&\triangleq \big(P_1 (C^{jj}_2)^\top+P_2 C^{jj}_4\big) (L^j_1)^{-1} P_2^\top+P_1\in\RR^{n\times n},\\
S^j_2&\triangleq \big[P_1 (B^j_2)^\top+\big(P_1 (C^{jj}_2)^\top+P_2 C^{jj}_4\big) (L^j_1)^{-1}P_2^\top (B^j_2)^\top\big] \varphi_{2}\\
&\quad +\big[-\big(P_1 (C^{jj}_2)^\top+P_2 C^{jj}_4\big) (L^j_1)^{-1} P_3 C^j_1+P_2 C^j_1\big] \varphi_{1}\\
&\quad +V^j_1+\big(P_1 (C^{jj}_2)^\top+P_2 C^{jj}_4\big) (L^j_1)^{-1} V^j_2\\
&\quad +\big(P_1 (C^{jj}_2)^\top+P_2 C^{jj}_4\big) (L^j_1)^{-1} P_2^\top E^j_2+P_1 E^j_2\in\RR^n,\\
L_5&\triangleq D_4+\sum_{j=1}^{l_2}(D^j_2)^\top (L^j_2)^{-1} S^j_1D^j_2\in\RR^{m\times m}, \\
L_6&\triangleq -L_5^{-1}\bigg(\sum_{j=1}^{l_2}(D^j_2)^\top(L^j_2)^{-1} L^j_3+D_1^\top P_1\bigg)\in\RR^{m\times n},\\ 
L_7&\triangleq -L_5^{-1}\bigg(D_3^\top+\sum_{j=1}^{l_2}(D^j_2)^\top (L^j_2)^{-1}L^j_4+D_1^\top P_2\bigg)\in\RR^{m\times Nn}, \\
\end{aligned}
\end{equation}
and
\begin{equation}\label{S3}
\begin{aligned}
S_3&\triangleq -L_5^{-1}\bigg(D_1^\top \varphi_1+\sum_{j=1}^{l_2}(D^j_2)^\top (L^j_2)^{-1} S^j_2\bigg)\in\RR^m,
\end{aligned}
\end{equation}
then we have
\begin{equation}\label{u}
\bar{u}=L_6 \bar{X}+L_7 h+S_3.
\end{equation}
Insert $\bar{u}$ of \eqref{u} back to \eqref{barZ} and \eqref{n}, and set
\begin{equation*}
\begin{aligned}
L^j_8&\triangleq (L^j_2)^{-1}\big(L^j_3+S^j_1 D^j_2 L_6\big)\in\RR^{n\times n},\quad L^j_9\triangleq (L^j_2)^{-1}\big(L^j_4+S^j_1 D^j_2 L_7\big)\in\RR^{n\times Nn},\\
L^j_{10}&\triangleq (L^j_1)^{-1}\big[P_2^\top A^j_2+P_2^\top (B^j_2)^\top P_2^\top+P_2^\top D^j_2 L_6-P_3 C^j_1 P_1-P_3 C^{jj}_2 L^j_8\big]\in\RR^{Nn\times n},\\
L^j_{11}&\triangleq (L^j_1)^{-1}\big[-P_2^\top (B^j_2)^\top P_3+P_2^\top D^j_2 L_7-P_3 C^j_3-P_3 C^j_1 P_2-P_3 C^{jj}_2 L^j_9\big]\in\RR^{Nn\times Nn},\\
S^j_4&\triangleq (L^j_2)^{-1}\big(S^j_1 D^j_2 S_3+S^j_2\big)\in\RR^n, \\
S^j_5&\triangleq (L^j_1)^{-1}\big[P_2^\top (B^j_2)^\top \vp_2+P_2^\top D^j_2 S_3+P_2^\top E^j_2-P_3 C^j_1 \varphi_1-P_3 C^{jj}_2 S^j_4+V^j_2\big]\in\RR^{Nn} ,
\end{aligned}
\end{equation*}
then we get
\begin{equation}\label{barZn}
\begin{aligned}
& n^j= L^j_8 \bar{X}+L^j_9 h+S^j_4, \\
& \bar{Z}^j=L^j_{10} \bar{X}+L^j_{11} h+S^j_5.
\end{aligned}
\end{equation}
Insert the pairs $(\bar{Z},n)$ of \eqref{barZn} and $(\bar{Y},m)$ of \eqref{indefrela} back, we have
\begin{equation}\label{P2P3}\left\{
\begin{aligned}
\dot{P}_2^\top&+ P_2^\top\big(A_1+\tilde{B}_1^\top P_2^\top\big)+\sum_{j=1}^{l_2}\big(P_2^\top (C^j_1)^\top+(C^j_3)^\top\big)L^j_{10} +\big(P_2^\top D_1+D_3\big)L_6\\
&\quad -P_3\big(\tilde{B}_1P_1+ B_4P_2^\top\big)- P_3\sum_{j=1}^{l_2}B^j_2 L^j_8+A_3+B_3^\top P_2^\top=0, \\
\dot{P}_3&+P_2^\top \tilde{B}_1^\top P_3-\sum_{j=1}^{l_2}\big(P_2^\top (C^j_1)^\top+(C^j_3)^\top\big) L^j_{11}-\big(P_2^\top D_1+D_3\big) L_7+P_3 B_3\\
&\quad +P_3 \tilde{B}_1 P_2+P_3\sum_{j=1}^{l_2}B^j_2 L^j_9-P_3B_4P_3+B_3^\top P_3=0, \\
P_2&(T)=F^\top,\ P_3(T)=0_{Nn\times Nn},
\end{aligned}
\right.\end{equation}
and
\begin{equation}\label{vp_V_2}
\begin{aligned}
d \varphi_2=&-\bigg[\big(P_2^\top \tilde{B}_1^\top-P_3 B_4+B_3^\top\big) \varphi_2-P_3 \tilde{B}_1 \varphi_1+\sum_{j=1}^{l_2}\big(P_2^\top (C^j_1)^\top+(C^j_3)^\top\big)S^j_{5}\\
&\quad +\big(P_2^\top D_1+D_3\big) S_3+P_2^\top E_1-P_3\sum_{j=1}^{l_2} B^j_2 S^j_4+E_3\bigg] d t +\sum_{j=1}^{l_2}V^j_2 d B^j_t.
\end{aligned}
\end{equation}
Next, we can also obtain
\begin{equation}\label{P1P2}\left\{
\begin{aligned}
\dot{P}_1&+P_1\big(A_1+\tilde{B}_1^\top P_2^\top\big)+P_2B_4P_2^\top+P_1 \sum_{j=1}^{l_2}(C^j_1)^\top L^j_{10}+P_1 D_1 L_6\\
&\quad +\big(P_2 \tilde{B}_1+A_1^\top\big) P_1+\sum_{j=1}^{l_2}\big(P_2 B^j_2+(A^j_2)^\top\big) L^j_8+A_4=0, \\
\dot{P}_2&-\big(P_1 \tilde{B}_1^\top+P_2B_4\big) P_3+P_1 \sum_{j=1}^{l_2}(C^j_1)^\top L^j_{11}+P_1 D_1 L_7+P_2 B_3\\
&\quad +\big(P_2 \tilde{B}_1+A_1^\top\big)P_2+\sum_{j=1}^{l_2}\big(P_2 B^j_2+(A^j_2)^\top\big) L^j_9+A_3^\top=0,\\
P_1&(T)=G,\ P_2(T)=F^\top,
\end{aligned}
\right.\end{equation}
and
\begin{equation}\label{vp_V_1}
\begin{aligned}
d \varphi_1=&-\bigg[\big(P_1 \tilde{B}_1^\top+P_2 B_4\big) \varphi_2+\big(P_2 \tilde{B}_1+A_1^\top\big) \varphi_1+P_1 \sum_{j=1}^{l_2}(C^j_1)^\top S^j_5\\
&\quad +P_1 D_1 S_3+P_1 E_1+\sum_{j=1}^{l_2}\big(P_2 B^j_2+(A^j_2)^\top\big) S^j_4\bigg] d t +\sum_{j=1}^{l_2}V^j_1 d B^j_t.
\end{aligned}
\end{equation}

Next, similar to the key idea in \cite{LZ01}, we can regard the original forward-backward SLQ optimal control problem as a dimension-enlarged forward SLQ one and build the relationship between them. Let
\begin{equation}\label{A25}
\begin{aligned}
\tilde{X}&\equiv \left(\begin{array}{l}
X \\
Y
\end{array}\right),\quad \tilde{u}^j\equiv \left(\begin{array}{l}
u \\
Z^j
\end{array}\right),\quad \tilde{A}\triangleq \left(\begin{array}{cc}
A_1 & \tilde{B}_1^\top \\
-A_3 & -B_3^\top
\end{array}\right),\\
\tilde{B}^j&\triangleq \left(\begin{array}{cc}
\frac{D_1}{l_2} & (C^j_1)^\top \\
-\frac{D_3}{l_2} & -(C^j_3)^\top
\end{array}\right), \quad
\tilde{C}^j\triangleq \left(\begin{array}{cc}
A^j_2 & (B^j_2)^\top \\
0 & 0
\end{array}\right),\\
\tilde{D}^j&\triangleq \left(\begin{array}{cc}
D^j_2 & (C^{jj}_2)^\top \\
0 & I_{Nn\times Nn}
\end{array}\right),\quad
\tilde{E}\triangleq \left(\begin{array}{c}
E_1 \\
-E_3
\end{array}\right),\quad \tilde{\tilde{E}}^j\triangleq \left(\begin{array}{c}
E^j_2 \\
0
\end{array}\right),
\end{aligned}
\end{equation}
then the state equation and cost functional become
\begin{equation}\label{A26}
d \tilde{X}=\bigg[\tilde{A} \tilde{X}+\sum_{j=1}^{l_2}\tilde{B}^j \tilde{u}^j+\tilde{E}\bigg] dt+\sum_{j=1}^{l_2}\big[\tilde{C}^j\tilde{X}+\tilde{D}^j \tilde{u}^j+\tilde{\tilde{E}}^j\big] d B^j_t,
\end{equation}
and
\begin{equation}\label{A27}
J(u(\cdot))=\frac{1}{2} \mathbb{E}\bigg[\int_0^T\bigg(\tilde{X}^\top \tilde{Q} \tilde{X}+\sum_{j=1}^{l_2}(\tilde{u}^j)^\top\tilde{R}^j\tilde{u}^j\bigg)dt +\langle G X_T,X_T\rangle\bigg],
\end{equation}
where
\begin{equation*}
\begin{aligned}
&\tilde{Q}\triangleq \left(\begin{array}{cc}
A_4 & 0 \\
0 & B_4
\end{array}\right)\in\RR^{(N+1)n\times(N+1)n},\quad \tilde{R}^j\triangleq \left(\begin{array}{cc}
D_4 & 0 \\
0 & C^j_4
\end{array}\right)\in\RR^{(Nn+m)\times (Nn+m)}.
\end{aligned}
\end{equation*}
Assume that $d \tilde{\varphi}=-\tilde{\ga} d t+\sum\limits_{j=1}^{l_2}\tilde{V}^j d B^j_t$, and apply It\^o's formula, we achieve
\begin{equation*}
\begin{aligned}
&d\big[(\tilde{X}-\tilde{\vp})^\top\tilde{P}(\tilde{X}-\tilde{\vp})\big] +\bigg[\tilde{X}^\top \tilde{Q} \tilde{X}+\sum_{j=1}^{l_2}(\tilde{u}^j)^\top\tilde{R}^j\tilde{u}^j\bigg]dt\\
&=(\tilde{X}-\tilde{\varphi})^\top \tilde{P}\bigg[\tilde{A} \tilde{X}+\sum_{j=1}^{l_2}\tilde{B}^j \tilde{u}^j+\tilde{E}+\tilde{\ga}\bigg] d t
 +(\tilde{X}-\tilde{\varphi})^\top \tilde{P}\sum_{j=1}^{l_2}\big[\tilde{C}^j \tilde{X}+\tilde{D}^j \tilde{u}^j+\tilde{\tilde{E}}^j-\tilde{V}^j\big] d B^j_t\\
&\quad +(\tilde{X}-\tilde{\varphi})^\top \dot{\tilde{P}}(\tilde{X}-\tilde{\varphi}) d t +\bigg[\tilde{X}^\top \tilde{A}^\top+\sum_{j=1}^{l_2}(\tilde{u}^j)^\top (\tilde{B}^j)^\top
 +\tilde{E}^\top+\tilde{\gamma}^\top\bigg] \tilde{P}(\bar{X}-\tilde{\varphi})dt\\
&\quad +\sum_{j=1}^{l_2}\big[\bar{X}^\top (\tilde{C}^j)^\top+(\tilde{u}^j)^\top(\tilde{D}^j)^\top+(\tilde{\tilde{E}}^j)^\top-(\tilde{V}^j)^\top\big] \tilde{P}(\tilde{X}-\tilde{\varphi}) dB^j_t\\
&\quad +\sum_{j=1}^{l_2}\big[\tilde{X}^\top (\tilde{C}^j)^\top+(\tilde{u}^j)^\top (\tilde{D}^j)^\top+(\tilde{\tilde{E}}^j)^\top-(\tilde{V}^j)^\top\big]
 \tilde{P}\big[\tilde{C}^j\tilde{X}+\tilde{D}^j\tilde{u}^j+\tilde{\tilde{E}}^j-\tilde{V}^j\big]dt \\
&\quad +\bigg[\tilde{X}^\top \tilde{Q} \tilde{X}+\sum_{j=1}^{l_2}(\tilde{u}^j)^\top\tilde{R}^j\tilde{u}^j\bigg]dt\\
&=\bigg\{\tilde{X}^\top\bigg[\tilde{P} \tilde{A}+\dot{\tilde{P}}+\tilde{A}^\top \tilde{P}+\sum_{j=1}^{l_2}(\tilde{C}^j)^\top \tilde{P} \tilde{C}^j+\tilde{Q}\bigg] \tilde{X}
 +\tilde{X}^\top\bigg[\tilde{P} \tilde{E}+\tilde{P} \tilde{\gamma}-\dot{\tilde{P}} \tilde{\vp}-\tilde{A}^\top \tilde{P} \tilde{\varphi}\\
&\qquad +\sum_{j=1}^{l_2}(\tilde{C}^j)^\top \tilde{P} \tilde{\tilde{E}}^j-\sum_{j=1}^{l_2}(\tilde{C}^j)^\top \tilde{P} \tilde{V}^j\bigg]
 +\bigg[-\tilde{\varphi}^\top \tilde{P} \tilde{A}-\tilde{\varphi}^\top \dot{\tilde{P}}+\tilde{E}^\top \tilde{P}+\tilde{\gamma}^\top \tilde{P}\\
&\qquad +\sum_{j=1}^{l_2}(\tilde{\tilde{E}}^j)^\top \tilde{P} \tilde{C}^j-\sum_{j=1}^{l_2}(\tilde{V}^j)^\top \tilde{P} \tilde{C}^j\bigg] \tilde{X}
 +\bigg[\bar{X}^\top\sum_{j=1}^{l_2}\big(\tilde{P}\tilde{B}^j+(\tilde{C}^j)^\top\tilde{P} \tilde{D}^j\big) \tilde{u}^j\\
&\qquad -\tilde{\varphi}^\top \tilde{P} \sum_{j=1}^{l_2}\tilde{B}^j \tilde{u}^j-\sum_{j=1}^{l_2}(\tilde{u}^j)^\top (\tilde{B}^j)^\top \tilde{P} \tilde{\varphi}
 +\sum_{j=1}^{l_2}(\tilde{u}^j)^\top\big((\tilde{B}^j)^\top \tilde{P}+(\tilde{D}^j)^\top \tilde{P} \tilde{C}^j\big)\tilde{X}\\
&\qquad +\sum_{j=1}^{l_2}(\tilde{u}^j)^\top \big(\tilde{R}^j+(\tilde{D}^j)^\top \tilde{P} \tilde{D}^j\big) \tilde{u}^j
 +\sum_{j=1}^{l_2}(\tilde{u}^j)^\top (\tilde{D}^j)^\top\tilde{P}\tilde{\tilde{E}}^j-\sum_{j=1}^{l_2}(\tilde{u}^j)^\top (\tilde{D}^j)^\top\tilde{P} \tilde{V}^j\\
&\qquad +\sum_{j=1}^{l_2}(\tilde{\tilde{E}}^j)^\top \tilde{P} \tilde{D}^j \tilde{u}^j-\sum_{j=1}^{l_2}(\tilde{V}^j)^\top \tilde{P} \tilde{D}^j \tilde{u}^j\bigg]
 -\tilde{\varphi}^\top \tilde{P} \tilde{E}-\tilde{\varphi}^\top \tilde{P} \tilde{\ga}+\tilde{\varphi}^\top \dot{\tilde{P}} \tilde{\varphi}\\
&\qquad -\tilde{E}^\top \tilde{P} \tilde{\varphi}-\tilde{\gamma}^\top \tilde{P} \tilde{\varphi}+\sum_{j=1}^{l_2}(\tilde{\tilde{E}}^j)^\top\tilde{P} \tilde{\tilde{E}}^j
 -\sum_{j=1}^{l_2}(\tilde{\tilde{E}}^j)^\top \tilde{P} \tilde{V}^j \\
&\qquad -\sum_{j=1}^{l_2}(\tilde{V}^j)^\top \tilde{P} \tilde{\tilde{E}}^j+\sum_{j=1}^{l_2}(\tilde{V}^j)^\top \tilde{P} \tilde{V}^j\bigg\}dt +\{\cdots\} dB_t^j\\
&=\bigg\{\tilde{X}^\top\bigg[\tilde{P}\tilde{A}+\dot{\tilde{P}}+\tilde{A}^\top\tilde{P}+\sum_{j=1}^{l_2}(\tilde{C}^j)^\top\tilde{P}\tilde{C}^j+\tilde{Q}
 -\sum_{j=1}^{l_2}(\tilde{P}^\top\tilde{B}^j+(\tilde{C}^j)^\top\tilde{P}^\top\tilde{D}^j)\\
&\qquad \times\big(\tilde{R}^j+(\tilde{D}^j)^\top\tilde{P}\tilde{D}^j\big)^{-1}\big((\tilde{B}^j)^\top \tilde{P}+(\tilde{D}^j)^\top\tilde{P} \tilde{C}^j\big)\bigg] \tilde{X}
 +\tilde{X}^\top\bigg[\tilde{P} \tilde{E}+\tilde{P} \tilde{\gamma}-\dot{\tilde{P}} \tilde{\vp}\\
&\qquad -\tilde{A}^\top \tilde{P} \tilde{\varphi}+\sum_{j=1}^{l_2}(\tilde{C}^j)^\top \tilde{P} \tilde{\tilde{E}}^j-\sum_{j=1}^{l_2}\big(\tilde{C}^j)^\top\tilde{P}\tilde{V}^j
 -\sum_{j=1}^{l_2}(\tilde{P}^\top\tilde{B}^j+(\tilde{C}^j)^\top\tilde{P}^\top\tilde{D}^j\big)\\
&\qquad \times\big(\tilde{R}^j+(\tilde{D}^j)^\top\tilde{P}\tilde{D}^j\big)^{-1}\big((\tilde{D}^j)^\top\tilde{P}\tilde{\tilde{E}}^j-(\tilde{B}^j)^\top\tilde{P}\tilde{\vp}-(\tilde{D}^j)^\top\tilde{P}\tilde{V}^j\big)\bigg]\\
\end{aligned}
\end{equation*}
\begin{equation*}
\begin{aligned}
&\qquad +\bigg[-\tilde{\varphi}^\top\tilde{P}\tilde{A}-\tilde{\vp}^\top \dot{\tilde{P}}+\tilde{E}^\top\tilde{P}+\tilde{\ga}^\top\tilde{P}
 +\sum_{j=1}^{l_2}(\tilde{\tilde{E}}^j)^\top \tilde{P} \tilde{C}^j-\sum_{j=1}^{l_2}(\tilde{V}^j)^\top\tilde{P} \tilde{C}^j\\
&\qquad -\sum_{j=1}^{l_2}\big[(\tilde{\tilde{E}}^j)^\top\tilde{P}^\top \tilde{D}^j -\tilde{\vp}^\top\tilde{P}^\top \tilde{B}^j-(\tilde{V}^j)^\top \tilde{P}^\top\tilde{D}^j\big]
 \big(\tilde{R}^j+(\tilde{D}^j)^\top \tilde{P}^\top \tilde{D}^j)^{-1}\\
&\qquad \times\big((\tilde{B}^j)^\top\tilde{P}+(\tilde{D}^j)^\top \tilde{P}\tilde{C}^j\big)\bigg]\tilde{X}
 -\sum_{j=1}^{l_2}\big[(\tilde{\tilde{E}}^j)^\top\tilde{P}^\top\tilde{D}^j-\tilde{\vp}^\top\tilde{P}^\top\tilde{B}^j\\
&\qquad -(\tilde{V}^j)^\top\tilde{P}^\top\tilde{D}^j\big]\big(\tilde{R}^j+(\tilde{D}^j)^\top\tilde{P}^\top\tilde{D}^j\big)^{-1}
 \big[(\tilde{D}^j)^\top\tilde{P}\tilde{\tilde{E}}^j-(\tilde{B}^j)^\top\tilde{P}\tilde{\vp}-(\tilde{D}^j)^\top\tilde{P}\tilde{V}^j\big]\\
&\qquad -\tilde{\varphi}^\top \tilde{P} \tilde{E}-\tilde{\varphi}^\top \tilde{P} \tilde{\gamma}+\tilde{\varphi}^\top\dot{\tilde{P}}\tilde{\varphi}
 -\tilde{E}^\top \tilde{P} \tilde{\varphi}-\tilde{\gamma}^\top \tilde{P} \tilde{\varphi}
 +\sum_{j=1}^{l_2}(\tilde{\tilde{E}}^j)^\top \tilde{P}\tilde{\tilde{E}}^j-\sum_{j=1}^{l_2}(\tilde{\tilde{E}}^j)^\top \tilde{P} \tilde{V}^j \\
&\qquad -\sum_{j=1}^{l_2}(\tilde{V}^j)^\top \tilde{P} \tilde{\tilde{E}}^j+\sum_{j=1}^{l_2}(\tilde{V}^j)^\top\tilde{P}\tilde{V}^j
 +\sum_{j=1}^{l_2}\Big\{\tilde{u}^j+\big(\tilde{R}^j+(\tilde{D}^j)^\top\tilde{P}\tilde{D}^j\big)^{-1}\big[\big((\tilde{B}^j)^\top\tilde{P}\\
&\qquad +(\tilde{D}^j)^\top\tilde{P}\tilde{C}^j\big)\tilde{X}+(\tilde{D}^j)^\top\tilde{P}\tilde{\tilde{E}}^j-(\tilde{B}^j)^\top\tilde{P}\tilde{\vp}-(\tilde{D}^j)^\top\tilde{P}\tilde{V}^j\big]\Big\}^\top
 \big(\tilde{R}^j+(\tilde{D}^j)^\top\tilde{P}\tilde{D}^j)\\
&\qquad \times\Big\{\tilde{u}^j+\big(\tilde{R}^j+(\tilde{D}^j)^\top\tilde{P}\tilde{D}^j)^{-1}\big[\big((\tilde{B}^j)^\top\tilde{P}+(\tilde{D}^j)^\top\tilde{P}\tilde{C}^j\big)\tilde{X}
 +(\tilde{D}^j)^\top\tilde{P}\tilde{\tilde{E}}^j\\
&\qquad -(\tilde{B}^j)^\top\tilde{P}\tilde{\vp}-(\tilde{D}^j)^\top\tilde{P}\tilde{V}^j\big]\Big\}\bigg\}dt+\{\cdots\}dB_t^j\\
&=\bigg\{\sum_{j=1}^{l_2}\Big\{\tilde{u}^j+(M^j_1)^{-1}\big[M^j_2\tilde{X}-M^j_3\tilde{\vp}-M^j_4\tilde{V}^j+M^j_4\tilde{\tilde{E}}^j\big]\Big\}^\top M^j_1\Big\{\tilde{u}^j+(M_1^j)^{-1}\big[M^j_2\tilde{X}\\
&\qquad -M^j_3\tilde{\vp}-M^j_4\tilde{V}^j+M^j_4\tilde{\tilde{E}}^j\big]\Big\} +\tilde{X}^\top\bigg[\tilde{P} \tilde{A}+\dot{\tilde{P}}+\tilde{A}^\top \tilde{P}+\sum_{j=1}^{l_2}(\tilde{C}^j)^\top \tilde{P} \tilde{C}^j+\tilde{Q}\\
&\qquad -\sum_{j=1}^{l_2}(M^j_2)^{\top} (M^j_1)^{-1} M^j_2\bigg] \tilde{X}+\tilde{X}^\top\bigg[\tilde{P} \tilde{E}+\tilde{P} \tilde{\ga}-\dot{\tilde{P}}\tilde{\varphi}
 -\tilde{A}^\top \tilde{P} \tilde{\varphi}+\sum_{j=1}^{l_2}(\tilde{C}^j)^\top \tilde{P}\tilde{\tilde{E}}^j\\
&\qquad -\sum_{j=1}^{l_2}(\tilde{C}^j)^\top \tilde{P} \tilde{V}^j-\sum_{j=1}^{l_2}(M^j_2)^\top\big(M^j_1)^{-1}(M^j_4 \tilde{\tilde{E}}^j-M^j_3 \tilde{\vp}-M^j_4 \tilde{V}^j\big)\bigg]
 +\bigg[-\tilde{\varphi}^\top \tilde{P}\tilde{A}\\
&\qquad -\tilde{\varphi}^\top \dot{\tilde{P}}+\tilde{E}^\top \tilde{P}+\tilde{\gamma}^\top \tilde{P}+\sum_{j=1}^{l_2}(\tilde{\tilde{E}}^j)^\top \tilde{P} \tilde{C}^j
 -\sum_{j=1}^{l_2}(\tilde{V}^j)^\top \tilde{P}\tilde{C}^j-\sum_{j=1}^{l_2}\big[(\tilde{\tilde{E}}^j)^\top (M^j_{4})^\top\\
&\qquad -\tilde{\varphi}^\top (M^j_{3})^\top-(\tilde{V}^j)^\top (M^j_4)^\top\big] (M^j_1)^{-1} M^j_2\bigg] \tilde{X}+M_5\bigg\}dt+\{ \cdots \} d B_t^j,
\end{aligned}
\end{equation*}
where
\begin{equation*}
\begin{aligned}
M^j_1&\triangleq \tilde{R}^j+(\tilde{D}^j)^\top \tilde{P} \tilde{D}^j\in\RR^{(Nn+m)\times(Nn+m)},\\
M^j_2&\triangleq (\tilde{B}^j)^\top \tilde{P}+(\tilde{D}^j)^\top \tilde{P} \tilde{C}^j\in\RR^{(Nn+m)\times(N+1)n},\\
M^j_3&\triangleq (\tilde{B}^j)^\top \tilde{P}\in\RR^{(Nn+m)\times(N+1)n}, \quad M^j_4\triangleq (\tilde{D}^j)^\top\tilde{P}\in\RR^{(Nn+m)\times(N+1)n},\\
\end{aligned}
\end{equation*}
and
\begin{equation*}
\begin{aligned}
M_5&\triangleq -\sum_{j=1}^{l_2}\big[\tilde{\tilde{E}}^j M_4^\top-\tilde{\varphi}^\top (M^j_3)^\top-(\tilde{V}^j)^\top (M^j_4)^\top\big] (M^j_1)^{-1}\big(M^j_4 \tilde{\tilde{E}}^j-M^j_3 \tilde{\varphi}-M^j_4 \tilde{V}^j\big) \\
 \end{aligned}
\end{equation*}
 \begin{equation}\label{M5}
\begin{aligned}
&\quad -\tilde{\varphi}^\top \tilde{P}\tilde{E}-\tilde{\varphi}^\top \tilde{P}\tilde{\gamma}+\tilde{\varphi}^\top \dot{\tilde{P}} \tilde{\varphi}-\tilde{E}^\top \tilde{P} \tilde{\varphi}
 -\tilde{\gamma}^\top \tilde{P} \tilde{\varphi}+\sum_{j=1}^{l_2}(\tilde{\tilde{E}}^j)^\top \tilde{P}\tilde{\tilde{E}}^j \\
&\quad -\sum_{j=1}^{l_2}(\tilde{\tilde{E}}^j)^\top\tilde{P} \tilde{V}^j-\sum_{j=1}^{l_2}(\tilde{V}^j)^\top \tilde{P}\tilde{\tilde{E}}^j+\sum_{j=1}^{l_2}(\tilde{V}^j)^\top\tilde{P} \tilde{V}^j.
\end{aligned}
\end{equation}
Then we introduce the Riccati equation:
\begin{equation}\label{A30}
\left\{
\begin{aligned}
&\dot{\tilde{P}}=-\tilde{P} \tilde{A}-\tilde{A}^\top \tilde{P}-\sum_{j=1}^{l_2}(\tilde{C}^j)^\top \tilde{P} \tilde{C}^j-\tilde{Q}+\sum_{j=1}^{l_2}(M^j_2)^{\top} (M^j_1)^{-1} M^j_2, \\
&\tilde{R}^j+(\tilde{D}^j)^\top \tilde{P} \tilde{D}^j>0,\\
& \bar{\tilde{u}}^j=-(M^j_1)^{-1}\big[M^j_2 \tilde{X}-M^j_3 \tilde{\varphi}-M^j_4 \tilde{V}^j+M^j_4 \tilde{\tilde{E}}^j\big],
\end{aligned}
\right.
\end{equation}
and
\begin{equation*}
\begin{aligned}
& \tilde{\gamma}=\tilde{P}^{-1}\bigg[\bigg(\dot{\tilde{P}}+\tilde{A}^\top \tilde{P}-\sum_{j=1}^{l_2}(M^j_2)^\top (M^j_1)^{-1} (M^j_3)^\top\bigg) \tilde{\varphi}+\sum_{j=1}^{l_2}\bigg((\tilde{C}^j)^\top \tilde{P}  \\
&\qquad -(M^j_2)^\top (M^j_1)^{-1} M^j_4\bigg)\tilde{V}^j-\tilde{P} \tilde{E}-\sum_{j=1}^{l_2}\big[(\tilde{C}^j)^\top \tilde{P}-(M^j_2)^\top (M_1^j)^{-1} M^j_4\big]\tilde{\tilde{E}}^j\bigg].
\end{aligned}
\end{equation*}
It is known that
\begin{equation*}
\begin{pmatrix}
m\\
-h
\end{pmatrix}\equiv\tilde{P}(\tilde{X}-\tilde{\vp})\equiv\begin{pmatrix}\tilde{P}_1(\bar{X}-\tilde{\vp}_1)+\tilde{P}_2(\bar{Y}-\tilde{\vp}_2)\\ \tilde{P}_2^\top(\bar{X}-\tilde{\vp}_1)+\tilde{P}_3(\bar{Y}-\tilde{\vp}_2)\end{pmatrix},
\end{equation*}
where
\begin{equation*}
\tilde{P}\equiv\begin{pmatrix}
\tilde{P}_1 & \tilde{P}_2 \\
\tilde{P}_2^\top & \tilde{P}_3
\end{pmatrix},\quad \tilde{\vp}=\begin{pmatrix}
\tilde{\vp}_1 \\
\tilde{\vp}_2
\end{pmatrix},\quad \tilde{V}^j=\begin{pmatrix}
\tilde{V}^j_1 \\
\tilde{V}^j_2
\end{pmatrix},
\end{equation*}
and compare with
\begin{equation*}
\begin{aligned}
 m&=\big(P_1+P_2P_3^{-1}P_2^\top \big)\bar{X}-P_2 P_3^{-1} \bar{Y}+P_2 P_3^{-1} \varphi_2+\varphi_1,\\
h&=P_3^{-1} P_2^\top \bar{X}-P_3^{-1} \bar{Y}+P_3^{-1} \vp_2,
\end{aligned}
\end{equation*}
then we can obtain
\begin{equation*}
\left\{
\begin{aligned}
&\tilde{P}_1=P_1+P_2P_3^{-1}P_2^\top, \quad \tilde{P}_2=-P_2 P_3^{-1},\quad \tilde{P}_3=P_3^{-1}, \\
&\tilde{\varphi}_1=-P_1^{-1} \varphi_1,\quad \tilde{V}^j_1=-P_1^{-1} V^j_1, \\
&\tilde{\varphi_2}=-P_2^\top P_1^{-1} \varphi_1+\varphi_2,\quad \tilde{V}^j_2=-P_2^\top P_1^{-1} V^j_1+V^j_2.
\end{aligned}
\right.
\end{equation*}
As a result of $P_3(T)=0$, $P^{-1}_3(T)$ makes no sense. We consider
\begin{equation*}
P^{i}(t)\equiv\begin{pmatrix}
P_{1, i} & P_{2, i} \\
P_{2, i}^\top & P_{3, i}
\end{pmatrix},\quad \varphi^i=\begin{pmatrix}
\varphi_{1, i} \\
\varphi_{2, i}
\end{pmatrix},\quad V^{j,i}=\begin{pmatrix}
V^j_{1, i} \\
V^j_{2, i}
\end{pmatrix}
\end{equation*}
to be the solution to equations corresponding to
\begin{equation}\label{P_vp}
P^i(T)=\begin{pmatrix}
G & F^\top \\
F & \frac{1}{i}I_{Nn\times Nn}
\end{pmatrix}, \quad \vp_T^i=\begin{pmatrix}
0 \\
\xi
\end{pmatrix}.
\end{equation}
It corresponds to a sequence of solutions to Riccati equations for
\begin{equation*}
\tilde{P}^i\equiv\begin{pmatrix}
\tilde{P}_{1,i} & \tilde{P}_{2,i} \\
\tilde{P}^\top_{2,i} & \tilde{P}_{3,i}
\end{pmatrix},\quad \tilde{\varphi}^i=\begin{pmatrix}
\tilde{\varphi}_{1,i} \\
\tilde{\varphi}_{2,i}
\end{pmatrix},\quad \tilde{V}^{j,i}=\begin{pmatrix}
\tilde{V}^j_{1,i} \\
\tilde{V}^j_{2,i}
\end{pmatrix},
\end{equation*}
with terminal conditions
\begin{equation}\label{A35}
\begin{aligned}
& \tilde{P}^i(T)=\begin{pmatrix}
G+iF^\top F& -iF^\top \\
-i F & iI_{Nn\times Nn}
\end{pmatrix}, \quad \tilde{\varphi}^i_T=\begin{pmatrix}
0 \\
\xi
\end{pmatrix},
\end{aligned}
\end{equation}
and similarly
\begin{equation}\label{R1}
\begin{aligned}
&P_{1,i}=\tilde{P}_{1,i}-P_{2,i} P_{3,i}^{-1}P_{2,i}^\top,\quad P_{2,i}=-\tilde{P}_{2,i}\tilde{P}_{3,i}^{-1},\quad P_{3,i}=\tilde{P}_{3,i}^{-1}, \\
& \vp_{1,i}=-P_{1, i} \tilde{\varphi}_{1,i},\quad V^j_{1,i}=-P_{1 i} \tilde{V}^j_{1, i}, \\
& \vp_{2,i}=\tilde{\varphi}_{2,i}-P_{2,i}^\top \tilde{\vp}_{1,i},\quad V^j_{2, i}=\tilde{V}^j_{2,i}-P_{2,i}^\top \tilde{V}^j_{1,i}.
\end{aligned}
\end{equation}
In the following, we have the relationships
\begin{equation}\label{mainrela}
\begin{aligned}
& (M^j_{1,i})^{-1} M^j_{2,i}=-\begin{pmatrix}
L_{6,i}+L_{7,i} P_{3,i}^{-1} P_{2,i}^\top & -L_{7,i} P_{3,i}^{-1} \\
L^j_{10,i}+L^j_{11,i} P_{3,i}^{-1} P_{2,i}^\top & -L^j_{11,i} P_{3,i}^{-1}
\end{pmatrix}\in\RR^{(Nn+m)\times(N+1)n},\\
& (M^j_{1,i})^{-1}\big[M^j_{3,i} \tilde{\varphi}^i+M^j_{4,i} \tilde{V}^{j,i}+M^j_{4,i}\tilde{\tilde{E}}^j\big]
=\begin{pmatrix}
L_{7,i} P_{3,i}^{-1} \varphi_{2,i}+S_{3,i} \\
L^j_{11,i} P_{3,i}^{-1} \varphi_{2,i}+S^j_{5,i}
\end{pmatrix}\in\RR^{Nn+m}.
\end{aligned}
\end{equation}
The proof of relation \eqref{mainrela} is similar to Theorem 4.4 in \cite{HJX23} and we omit it. Moreover, the solvabilities of Riccati equations can be also obtained as those in \cite{HJX23}, but we still illustrate the logic of proof in the following Theorem \ref{Theorem_A1} for the completion.

Next, we derive the feedback optimal control in the indefinite case. First, according to Lemma 4.12 in \cite{HJX23}, it implies that $M_{5,i}$ is uniformly bounded, and assume that $P_{3}(t)=\lim\limits _{i \rightarrow \infty} P_{3, i}(t)>0,\ t\in[0,T)$. Set $\bar{u}=L_6 X^*+L_7 h^*+S_3$, and let
\begin{equation*}
\begin{aligned}
&N_1\triangleq A_1+\tilde{B}_1^\top P_2^\top+\sum_{j=1}^{l_2}(C^j_1)^\top L^j_{10}+D_1 L_6,\quad N_2\triangleq -\tilde{B}_1^\top P_3+\sum_{j=1}^{l_2}(C^j_1)^\top L^j_{11}+D_1 L_7,\\
&N_3\triangleq \tilde{B}_1^\top \varphi_2+\sum_{j=1}^{l_2}(C^j_1)^\top S^j_5+E_{1}+D_1 S_3,\quad N^j_4\triangleq A^j_2+(B^j_2)^\top P_2^\top+(C^{jj}_2)^\top L^j_{10}+D^j_2 L_6,\\
&N^j_5\triangleq -(B^j_2)^\top P_3+(C^{jj}_2)^\top L^j_{11}+D^j_2 L_7,\quad N^j_6\triangleq (B^j_2)^\top \vp_2+(C^{jj}_2)^\top S^j_5+D^j_2 S_3+E^j_2,\\
&N_7\triangleq \tilde{B}_1 P_1+\sum_{j=1}^{l_2}B^j_2 L^j_8+B_4 P_2^\top,\quad N_8\triangleq B_3+\tilde{B}_1 P_2+\sum_{j=1}^{l_2}B^j_2 L^j_9-B_4 P_3,\\
&N_9\triangleq \tilde{B}_1 \varphi_1+\sum_{j=1}^{l_2}B^j_2 S^j_4+B_4 \varphi_2,\quad N^j_{10}\triangleq C^j_1 P_1+C^{jj}_2 L^j_8+C^{jj}_4 L^j_{10},\\
&N^j_{11}\triangleq C^j_3+C^j_1 P_2+C^{jj}_2 L^j_9+C^{jj}_4 L^j_{11},\quad N^j_{12}\triangleq C^j_1 \vp_1+C^{jj}_2 S^j_4+C^{jj}_4 S^j_5, \\
\end{aligned}
\end{equation*}
then we have
\begin{equation*}
\left\{
\begin{aligned}
d X^*&=\big(N_1 X^*+N_2 h^*+N_3\big) d t+\sum_{j=1}^{l_2}\big(N^j_4 X^*+N^j_5 h^*+N^j_6\big) d B^j_t, \\
d h^*&=\big(N_7 X^*+N_8 h^*+N_9\big) d t+\sum_{j=1}^{l_2}\big(N^j_{10} X^*+N^j_{11} h^*+N^j_{12}\big) d B^j_t, \\
X_0^*&=x_0,\ h_0^*=\big(I_{Nn\times Nn}+HP_3(0)\big)^{-1} H\big(P_2^\top(0)X_0+\vp_2(0)\big).
\end{aligned}
\right.
\end{equation*}
Moreover
\begin{equation*}
\begin{aligned}
& \bar{X}_t =X^*_t,\quad h_t=h^*_t, \\
& n^j_t=L^j_8 \bar{X}_t+L^j_9h_t+S^j_4, \quad m_t=P_1(t) \bar{X}_t+P_2(t) h_t+\vp_1(t), \\
& \bar{Z}^j_t=L^j_{10} \bar{X}_t+L^j_{11} h_t+S^j_5, \quad \bar{Y}_t=P^\top_2(t) \bar{X}_t-P_3(t)h_t+\varphi_2(t).
\end{aligned}
\end{equation*}
Set
\begin{equation*}
\begin{aligned}
\tilde{P}=\begin{pmatrix}
P_1+P_2 P_3^{-1}P_2^\top & -P_2 P_3^{-1} \\
-P_3^{-1}P_2^\top & P_3^{-1}
\end{pmatrix}.
\end{aligned}
\end{equation*}
For any given $\varepsilon>0$, $\tilde{P}$ solves the Riccati equations on $[0,T-\varepsilon)$. By the completion of square technique, recalling that
\begin{equation*}
\begin{aligned}
&d\big[(\tilde{X}-\tilde{\vp})^\top\tilde{P}(\tilde{X}-\tilde{\vp})\big] +\bigg[\tilde{X}^\top \tilde{Q} \tilde{X}+\sum_{j=1}^{l_2}(\tilde{u}^j)^\top\tilde{R}^j\tilde{u}^j\bigg]dt\\
&=\bigg\{\sum_{j=1}^{l_2}\Big\{\tilde{u}^j+(M^j_1)^{-1}\big[M^j_2\tilde{X}-M^j_3\tilde{\vp}-M^j_4\tilde{V}^j+M^j_4\tilde{\tilde{E}}^j\big]\Big\}^\top M^j_1\Big\{\tilde{u}^j+(M_1^j)^{-1}\big[M^j_2\tilde{X}\\
&\qquad -M^j_3\tilde{\vp}-M^j_4\tilde{V}^j+M^j_4\tilde{\tilde{E}}^j\big]\Big\} +\tilde{X}^\top\bigg[\tilde{P} \tilde{A}+\dot{\tilde{P}}+\tilde{A}^\top \tilde{P}+\sum_{j=1}^{l_2}(\tilde{C}^j)^\top \tilde{P} \tilde{C}^j+\tilde{Q}\\
&\qquad -\sum_{j=1}^{l_2}(M^j_2)^{\top} (M^j_1)^{-1} M^j_2\bigg] \tilde{X}+\tilde{X}^\top\bigg[\tilde{P} \tilde{E}+\tilde{P} \tilde{\ga}-\dot{\tilde{P}}\tilde{\varphi}
 -\tilde{A}^\top \tilde{P} \tilde{\varphi}+\sum_{j=1}^{l_2}(\tilde{C}^j)^\top \tilde{P}\tilde{\tilde{E}}^j\\
&\qquad -\sum_{j=1}^{l_2}(\tilde{C}^j)^\top \tilde{P} \tilde{V}^j-\sum_{j=1}^{l_2}(M^j_2)^\top\big(M^j_1)^{-1}(M^j_4 \tilde{\tilde{E}}^j-M^j_3 \tilde{\vp}-M^j_4 \tilde{V}^j\big)\bigg]
 +\bigg[-\tilde{\varphi}^\top \tilde{P}\tilde{A}\\
&\qquad -\tilde{\varphi}^\top \dot{\tilde{P}}+\tilde{E}^\top \tilde{P}+\tilde{\gamma}^\top \tilde{P}+\sum_{j=1}^{l_2}(\tilde{\tilde{E}}^j)^\top \tilde{P} \tilde{C}^j
 -\sum_{j=1}^{l_2}(\tilde{V}^j)^\top \tilde{P}\tilde{C}^j-\sum_{j=1}^{l_2}\big[(\tilde{\tilde{E}}^j)^\top (M^j_{4})^\top\\
&\qquad -\tilde{\varphi}^\top (M^j_{3})^\top-(\tilde{V}^j)^\top (M^j_4)^\top\big] (M^j_1)^{-1} M^j_2\bigg] \tilde{X}+M_5\bigg\}dt+\{ \cdots \} d B_t^j,
\end{aligned}
\end{equation*}
we obtain
\begin{equation*}
\begin{aligned}
J(\bar{u})&=\frac{1}{2}\mathbb{E}\bigg[\int_0^{T-\varepsilon}\bigg(\bar{X}A_4 \bar{X}+\bar{Y}^\top B_4 \bar{Y}+\sum_{j=1}^{l_2}(\bar{Z}^j)^\top C^j_4 \bar{Z}^j+\bar{u}D_4 \bar{u}\bigg)d t\\
&\qquad +\int_{T-\varepsilon}^T\bigg(\bar{X}A_4 \bar{X}+\bar{Y}^\top B_4 \bar{Y}+\sum_{j=1}^{l_2}(\bar{Z}^j)^\top C^j_4 \bar{Z}^j+\bar{u}D_4\bar{u}\bigg)dt+\bar{Y}^\top_0H\bar{Y}_0+\bar{X}_TG\bar{X}_T\bigg]\\
&=\frac{1}{2}\mathbb{E}\big[(\tilde{X}_0-\tilde{\vp}_0)^\top\tilde{P}_0(\tilde{X}_0-\tilde{\vp}_0)-(\tilde{X}_{T-\epsilon}-\tilde{\vp}_{T-\epsilon})^\top\tilde{P}_{T-\epsilon}(\tilde{X}_{T-\epsilon}-\tilde{\vp}_{T-\epsilon})\big]\\
&\quad +\frac{1}{2}\mathbb{E}\int_0^{T-\epsilon}\sum_{j=1}^{l_2}\big[\tilde{u}^j+(M^j_1)^{-1}\big(M^j_2\tilde{X}-M^j_3 \tilde{\varphi}-M^j_4 \tilde{V}^j+M^j_4\tilde{\tilde{E}}^j\big)\big]^\top \\
&\qquad \times M^j_1\big[\tilde{u}^j+(M^j_{1})^{-1}\big(M^j_2 \tilde{X}-M^j_3 \tilde{\varphi}-M^j_{4} \tilde{V}+M^j_{4}\tilde{\tilde{E}}^j\big)\big]dt+\frac{1}{2}\mathbb{E}\int_0^{T-\epsilon} M_5 dt\\
&\quad +\frac{1}{2}\mathbb{E} \int_{T-\varepsilon}^T\bigg[\bar{X} A_4 \bar{X}+\bar{Y}^\top B_4 \bar{Y}+\sum_{j=1}^{l_2}(\bar{Z}^j)^\top C^j_4 \bar{Z}^j+\bar{u}D_4 \bar{u}\bigg]dt\\
&\quad +\frac{1}{2}\mathbb{E}\big[\bar{Y}_0^\top H\bar{Y}_0+\bar{X}_TG\bar{X}_T\big],
\end{aligned}
\end{equation*}
where
\begin{equation*}
\begin{aligned}
&\big(X_0^\top-\tilde{\varphi}_1^\top(0), \bar{Y}^\top_0-\tilde{\vp}^\top_{2}(0)\big)\begin{pmatrix}
\tilde{P}_1(0) & \tilde{P}_2(0) \\
\tilde{P}_2^\top(0) & \tilde{P}_3(0)
\end{pmatrix}\begin{pmatrix}
X_0-\vp_1(0) \\
\bar{Y}_0-\tilde{\vp}_2(0)
\end{pmatrix}+\bar{Y}_0^\top H\bar{Y}_0 \\
& =\RR_2+\big(\bar{Y}_0-(I_{Nn\times Nn}+P_3(0) H\big)^{-1}\big(P_2^\top(0) X_0+\varphi_2(0)\big)\big)^\top\big(H+P_3^{-1}(0)\big)\\
&\quad \times\big(\bar{Y}_0-(I_{Nn\times Nn}+P_3(0) H\big)^{-1}\big(P_2^\top(0) X_0+\varphi_2(0)\big)\big),
\end{aligned}
\end{equation*}
with
\begin{equation*}
\begin{aligned}
\RR_2&\triangleq \big(X_0+P_1^{-1}(0) \varphi_1(0)\big)^\top P_1(0)\big(X_0+P_1^{-1}(0) \varphi_1(0)\big)\\
&\quad +\big(P_2^\top(0) X_0+\varphi_2(0)\big)^\top\big(H P_3(0)+I_{Nn\times Nn}\big)^{-1} H\big(P_2^\top(0) X_0+\varphi_2(0)\big),
\end{aligned}
\end{equation*}
and
\begin{equation*}
\begin{aligned}
&\bar{X}_T^\top G\bar{X}_T-\big(\tilde{X}_{T-\varepsilon}-\tilde{\varphi}_{T-\varepsilon}\big)^\top\tilde{P}(T-\varepsilon)\big(\tilde{X}_{T-\vare}-\tilde{\vp}_{T-\epsilon}\big)\\
&=\bar{X}_T^\top G\bar{X}_T-\big(\bar{X}_{T-\vare}+P_1^{-1}(T-\vare)\vp_1(T-\vare)\big)^\top P_1(T-\vare)\big(\bar{X}_{T-\vare}\\
&\quad +P_1^{-1}(T-\vare)\vp_1(T-\vare)\big)-h_{T-\vare}^\top P_3(T-\vare)h_{T-\vare}\rightarrow0, \ \text{when}\ \vare\rightarrow0.
\end{aligned}
\end{equation*}
Similar to \cite{HJX23}, we have
\begin{equation*}
J(\bar{u})=\frac{1}{2} \mathbb{E}\bigg[\RR_2+\int_0^T M_5(t)dt\bigg],
\end{equation*}
and $J(u) \geq J(\bar{u})$, which means $\bar{u}$ is optimal, and finally by \eqref{mainrela}
\begin{equation}\label{A41}
\bar{u}=\big(L_6+L_7 P_3^{-1} P_2^\top\big) \bar{X}-L_7 P_3^{-1} \bar{Y}+L_7 P_3^{-1} \varphi_2+S_3.
\end{equation}
where $P_1(t)=\lim_{i\rightarrow\infty}P_{1,i}, P_2(t)=\lim_{i\rightarrow\infty}P_{2,i}$ and $L_6, L_7, S_3$ are defined in \eqref{L1}, \eqref{S3}.

We conclude the above analysis as the following main theorem.
\begin{theorem}\label{Theorem_A1}
Suppose that (i) the bounded conditions for coefficients with appropriate dimensions similar to \cite[Assumption 2.3 (i)-(ii)]{HJX23} and $E_1,E^j_2,E_{3i}\in L^\infty(0,T;\RR^n), i=1,\dots,N, j=1,\dots,l_2$ hold, and assume that (ii) there exist $i_0$ such that for $i>i_0$, the Riccati equation \eqref{A30} has a positive definite solution $\tilde{P}^i(t)$ satisfying terminal condition \eqref{A35}. Moreover, suppose that (iii) $(L^j_{1,i}(t))^{-1}, (L^j_{2,i}(t))^{-1}$ exist, and \eqref{R1} hold and $L^j_{1,i},L^j_{2,i}$ are defined as above with replacing $P_1(t),P_2(t),P_3(t)$ by $P_{1,i}(t),P_{2,i}(t),P_{3,i}(t)$. (iv) Let for each $i\geq i_0$, $(\tilde{R}^j(\cdot)+(\tilde{D}^j(\cdot))^\top\tilde{P}(\cdot)\tilde{D}^j(\cdot))^{-1}$ is bounded, $j=1,\dots,l_2$ and $P_{1,i}$ have an upper bound and $|P_{2,i}(\cdot)|$, $(L^j_{1,i}(\cdot))^{-1},(L^j_{2,i}(\cdot))^{-1},(L_{5,i}(\cdot))^{-1}$ are uniformly bounded. Assume that (v) $P_3(t)=\lim_{i\rightarrow\infty}P_{3,i}(t)>0$ for $t\in[0,T)$. Then the optimal control problem \eqref{A1} and \eqref{A2} has an optimal control \eqref{A41}.
\end{theorem}
\begin{proof}
The proof is similar to those in \cite{HJX23}, so we mainly sketch the proof framework. Indeed, under assumptions (i), (ii) and (iii), we have $P_{1,i}(\cdot),P_{2,i}(\cdot),P_{3,i}(\cdot)$ solves equations \eqref{P2P3} and \eqref{P1P2} and $(\vp^i,V^{j,i})$ solves the equations \eqref{vp_V_2} and \eqref{vp_V_1} with \eqref{P_vp} by \cite[Theorem 4.3]{HJX23}. Meanwhile, we have the relations \eqref{mainrela} hold by \cite[Theorem 4.4]{HJX23}. Moreover, in addition to (i)-(iii), under assumption (iv), equations \eqref{P2P3} and \eqref{P1P2} admits a unique solution $(P_1(\cdot),P_2(\cdot),P_3(\cdot))$ by \cite[Theorem 4.6]{HJX23}. Therefore, from assumptions (i)-(v), the optimal control problem \eqref{A1} and \eqref{A2} has an optimal control \eqref{A41} by \cite[Theorem 4.13]{HJX23}.
\end{proof}

\begin{remark}
By the above analysis, we can find the optimal control \eqref{A41} actually take the same form as \eqref{u}. However, only the candidate \eqref{u} is obtained (may be optimal in form) based on the guess of relations \eqref{indefrela} (motivated by relations \eqref{defrela}) but not gives the optimal cost in the indefinite case. Therefore, the core of adopting the second idea from \eqref{A25}, \eqref{A26} and \eqref{A27} is to guarantee the final form of optimal control \eqref{A41}(or \eqref{u}) under the classical forward stochastic LQ control problem by the completion of square technique and most importantly to get the optimal cost correspondingly based on the relation \eqref{R1} and \eqref{mainrela} when $i$ converges to infinity.
\end{remark}

\end{document}